\documentclass[11pt, french]{article}
\usepackage{latexsym} 
\usepackage[all]{xy}
\usepackage{amsmath}
\usepackage{amssymb}
\usepackage{amsfonts}
\usepackage[applemac]{inputenc}

\newtheorem{thm}{Th\'eor\`eme}[section]
\newtheorem{prop}[thm]{Proposition}  
\newtheorem{lem}[thm]{Lemme}

\newtheorem{df}[thm]{D\'efinition}
\newtheorem{cor}[thm]{Corollaire}

\newtheorem{rmk}[thm]{Remarque}

\newcommand{\s}{\infty}
\newcommand{\scps}{\infty-\mathbf{Cat}^{pr-\otimes}}
\newcommand{\scms}{\infty-\mathbf{Cat}^{\otimes}}
\newcommand{\uscms}{\infty-\underline{\mathbf{Cat}}^{\otimes}}
\newcommand{\usmon}{\infty-\underline{\mathbf{Mon}}}
\newcommand{\ustop}{\infty-\underline{\mathbf{Top}}}
\newcommand{\uscat}{\infty-\underline{\mathbf{Cat}}}
\newcommand{\uscomm}{\infty-\underline{\mathbf{Comm}}}
\newcommand{\T}{\mathbb{T}}
\newcommand{\rh}{\mathbb{R}\underline{\mathrm{Hom}}}

\newcommand{\Feyn}[1]{#1\kern-0.45em/}

\begin{document}

\title{\textbf{Caract\`eres de Chern, traces \'equivariantes et g\'eom\'etrie alg\'ebrique d\'eriv\'ee}}

\bigskip
\bigskip

\author{\bf{Bertrand To\"en}\\ \small{I3M UMR 5149} \\  \small{Universit\'e de Montpellier2 - France} \\ \small{Montpellier - France} \and \bf{Gabriele Vezzosi}\\ \small{Dipartimento di Sistemi ed Informatica} \\ \small{Sezione di Matematica}\\ \small{Universit\`a di Firenze}\\ \small{Firenze - Italy} \bigskip }

\bigskip
\bigskip

\date{F\'evrier 2011}

\maketitle

\begin{abstract}
L'objectif de ce travail est de donner un traitement d\'etaill\'e de la 
construction du caract\`ere de Chern pour certaines familles en cat\'egories (esquiss\'ee dans l'article 
\cite{tove}). Pour cela nous introduisons, et \'etudions, la notion 
d'$\s$-cat\'egorie mono\"\i dale sym\'etrique rigide. Nous construisons un 
morphisme de trace dans ce cadre, qui est un morphisme de l'$\infty$-groupo\"\i de
des endomorphismes d'objets dans une telle $\s$-cat\'egorie vers celui des
endomorphismes de l'unit\'e. En utilisant le travail r\'ecent d'Hopkins-Lurie sur
\emph{l'hypoth\`ese du cobordisme} (voir \cite{lu2}) nous montrons de plus
que ce morphisme de trace satisfait une propri\'et\'e remarquable 
d'invariance cyclique. Nous utilisons l'existence de cette trace
cyclique afin de construire un caract\`ere de Chern, d\'efini pour tout couple
$(T,\mathcal{A})$ form\'e d'un $\s$-topos $T$ et d'un champ $\mathcal{A}$ en 
$\s$-cat\'egories mono\"\i dales sym\'etriques rigides. Nous pr\'esentons 
deux applications de notre construction g\'en\'erale, obtenues en sp\'ecifiant
le couple $(T,\mathcal{A})$. Nous montrons d'une part comment on peut
retrouver le caract\`ere de Chern des complexes parfaits \`a valeurs dans
l'homologie cyclique et comment notre construction permet de 
l'\'etendre de fa\c{c}on pertinente au cas des complexes parfaits sur
des champs alg\'ebriques d'Artin. Enfin, nous montrons comment 
notre caract\`ere de Chern permet de construire des invariants
de familles alg\'ebriques de dg-cat\'egories. Une  cons\'equence 
de l'existence de ces invariants est la construction d'une connexion de
Gauss-Manin sur le complexe d'homologie cyclique d'une telle
famille g\'en\'eralisant les constructions de \cite{ge,dotats}. Nous montrons aussi
comment on peut construire le \emph{faisceaux des caract\`eres} d'une repr\'esentation
d'un groupe alg\'ebrique dans une dg-cat\'egorie, qui est cat\'egorification 
de la fonction caract\`ere d'une repr\'esentation lin\'eaire ainsi qu'une
extension au cas dg-cat\'egorique de la construction de \cite{gaka}. Pour finir, lorsque
l'on dispose d'une famille de dg-cat\'egories satur\'ees nous construisons un
\emph{caract\`ere de Chern secondaire}, dont l'existence \'etait annonc\'ee
dans \cite{tove}, et \`a valeurs dans une nouvelle th\'eorie cohmologique
que nous appelons \emph{l'homologie cyclique secondaire}.
\end{abstract}

\tableofcontents

\section*{Introduction}

L'objectif de ce travail est de donner un traitement d\'etaill\'e de la 
construction du caract\`ere de Chern pour certaines familles en cat\'egories.
Plus sp\'ecifiquement, nous pr\'esentons ici un formalisme
tr\`es g\'en\'eral dans lequel il est possible de d\'emontrer l'existence d'un caract\`ere de 
Chern d\'efini suivant les grandes lignes de \cite{tove}, mais dont le domaine d'application 
d\'epasse largement le contexte des \emph{faisceaux cat\'egoriques d\'eriv\'es} de l'article en question. 
Avant d'entrer dans des d\'etails plus techniques nous nous permettons de commencer
cette introduction par un bref survol de plusieurs contextes d'applications de notre
construction.
Elle nous permet \emph{d'une part} de retrouver le caract\`ere
de Chern des fibr\'es vectoriels et des complexes parfaits (\`a valeurs dans l'homologie
cyclique), sur des sch\'emas, des espace analytiques, des vari\'et\'es diff\'erentielles etc. 
mais aussi sur des champs alg\'ebriques pour les quels il donne une version du 
\emph{caract\`ere de Chern orbifold} (dans le style de \cite{jkk,to}).  Mais il permet \emph{aussi} de construire certains invariants
de \emph{faisceaux en dg-cat\'egories}, 
comme la construction d'une connexion de Gauss-Manin sur l'homologie cyclique
p\'eriodique d'une famille alg\'ebrique de dg-cat\'egories, ou comme
le \emph{faisceau des caract\`eres} associ\'e \`a une \emph{repr\'esentation 
dg-cat\'egorique} d'un groupe alg\'ebrique, qui est une version
cat\'egorique de la fonction des caract\`eres associ\'ee \`a une repr\'esentation lin\'eaire
ainsi qu'une extension au cadre dg-cat\'egorique de la construction de \cite{gaka}.
Pour les \emph{faisceaux cat\'egoriques d\'eriv\'es parfaits} notre construction permet 
bien \'evidemment de construire le caract\`ere de Chern \`a valeurs dans
\emph{l'homologie cyclique secondaire} (une nouvelle th\'eorie homologique) dont l'existence est 
affirm\'ee, sans arguments r\'eellement solides, dans \cite{tove}. Le degr\'e de g\'en\'eralit\'e
de notre construction permet probablement de trouver de nombreux autres champs 
d'applications, par exemple pour d\'efinir un caract\`ere de Chern 
dans des contextes non-additifs, comme ceux apparaissant dans l'\'etude des
sch\'emas sur le corps \`a un \'el\'ement, ou plus g\'en\'eralement en 
g\'eom\'etrie relative de \cite{tova}; nous n'avons cependant pas
explorer cette voie dans ce travail. \\

\noindent \textbf{La construction dans ses grandes lignes } -- 
Avant d'entrer plus en d\'etails dans la description du contenu de ce travail
nous souhaitons pr\'esenter, en quelques mots, les grandes lignes 
de la construction de notre caract\`ere de Chern. Nous esp\'erons que cela
pourra eclaircir et expliquer le point de vue adopt\'e, \`a savoir celui des
$\infty$-cat\'egories, qu'il nous semble difficile d'\'eviter (m\^eme si l'on 
souhaite ne consid\'erer que le caract\`ere de Chern des fibr\'es vectoriels).

La construction se r\'esume probablement le plus ais\'ement dans le cadre
alg\'ebrique (notons cependant qu'un peu d'imagination permet aussi 
de concevoir des analogues dans le cadre diff\'erentielle, ou encore complexe
analytique). Soit donc $k$ un anneau commutatif, 
$A$ une $k$-alg\`ebre commutative et $M$ un $A$-module projectif de type fini, que nous voyons comme 
le module des sections d'un fibr\'e vectoriel sur le $k$-sch\'ema $X=Spec\, A$. 
Nous souhaitons d\'efinir le caract\`ere de Chern $Ch(M)$ \`a valeurs dans
$HC_{0}^{-}(A/k)$, l'homologie cyclique n\'egative de $A$ (ou de $X$) relativement 
\`a $k$. Pour cela nous consid\'erons la $k$-alg\`ebre commutative \emph{simpliciale} 
$S^{1}\otimes_{k}^{\mathbb{L}}A$, obtenue en tensorisant, de fa\c{c}on
correctement d\'eriv\'ee, $A$ par le groupe simplicial $S^{1}=B\mathbb{Z}$. 
On dispose d'une projection naturelle 
$p : A \longrightarrow S^{1}\otimes_{k}^{\mathbb{L}}A$
(induite par l'inclusion 
du point de base $* \rightarrow S^{1}$).
Ce morphisme poss\`ede de plus une auto-homotopie tautologique, induite
par l'identit\'e de $S^{1}\otimes_{k}^{\mathbb{L}}A$. Cette auto-homotopie 
induit \`a son tour un automorphisme tautologique du foncteur de changement de base
$$(S^{1}\otimes_{k}^{\mathbb{L}}A)\otimes_{A} - :
\mathbf{Proj}^{\textrm{tf}}(A) \longrightarrow \mathrm{Ho}(S^{1}\otimes_{k}^{\mathbb{L}}A\textbf{-Mod}),$$
de la cat\'egorie des $A$-modules projectifs de type fini vers 
la cat\'egorie homotopique des modules \emph{simpliciaux} sur l'anneau simplicial
$S^{1}\otimes_{k}^{\mathbb{L}}A$. Ainsi, on dispose d'un automorphisme
tautologique
$$m : (S^{1}\otimes_{k}^{\mathbb{L}}A)\otimes_{A} M \simeq 
(S^{1}\otimes_{k}^{\mathbb{L}}A)\otimes_{A} M.$$
On observe que $M$ est un objet dualisable de la cat\'egorie
monoidale des $A$-modules (car il est projectif et de type fini), et qu'il en est
donc de m\^eme de $(S^{1}\otimes_{k}^{\mathbb{L}}A)\otimes_{A} M$ en tant 
qu'objet de $\mathrm{Ho}(S^{1}\otimes_{k}^{\mathbb{L}}A-Mod)$. L'automorphisme
$m$ poss\`ede ainsi une trace bien d\'efini, qui est un \'el\'ement
$Tr(m) \in  
\pi_{0}(S^{1}\otimes_{k}^{\mathbb{L}}A)$.
Une des observations principales de ce travail est 
que la trace $Tr(m)$ poss\`ede un rel\`evement naturel en un 
\'el\'ement qui est homotopiquement invariant par l'action de $S^{1}$
$$Tr^{S^{1}}(m) \in \pi_{0}((S^{1}\otimes_{k}^{\mathbb{L}}A)^{hS^{1}}).$$
Nous affirmons que $Tr^{S^{1}}(m)$ \emph{est} le caract\`ere de Chern 
de $M$. Cette affirmation est bas\'ee sur une identification, que nous ne
d\'emontrerons pas dans cet article\footnote{Dans le cas o\`u $A$ est lisse sur $k$, l'\'equivalence entre $S^{1}\otimes_{k}^{\mathbb{L}} A$, muni de l'action naturelle de $S^{1}$, et l'alg\`ebre de de Rham de $A/k$, munie de sa diff\'erentielle de
de Rham, \`a et\'e etablie dans \cite{rhamloop}. Voir aussi \cite{d&d}.}, mais qui semble plus ou moins
connue de fa\c{c}on folklorique
$$\pi_{*}((S^{1}\otimes_{k}^{\mathbb{L}}A)^{hS^{1}}) \simeq HC_{*}^{-}(A/k).$$
Cette construction du caract\`ere de Chern (qui sera rendu 
rigoureuse, voir \S 4.2), sera notre source principale d'inspiration, 
et l'objectif de ce travail consiste essentiellement \`a pr\'esenter un 
cadre abstrait et beaucoup plus general dans lequel elle garde un sens. \\

\noindent \textbf{$\infty$-Topos catannel\'es rigides --}  
G\'eom\'etriquement la construction que nous venons d'esquisser 
s'interpr\`ete de la mani\`ere
suivante. L'anneau simplicial $S^{1}\otimes_{k}^{\mathbb{L}}A$
est \emph{l'anneau des fonctions sur l'espace des lacets
$LX=Map(S^{1},X)$ du sch\'ema $X$}, o\`u $LX$ est construit dans la 
cat\'egorie des $k$-sch\'emas d\'eriv\'es (voir par exemple \cite{hagII,to3}). Le caract\`ere 
de Chern du fibr\'e vectoriel $M$ sur $X$ est alors la fonction sur 
$LX$ qui \`a un \emph{lacet $\gamma$ sur $X$} associe la trace de \emph{la monodromie de $V$} 
le long de $\gamma$, fonction que l'on observe \^etre naturellement invariante par l'action 
de $S^{1}$ sur $LX$ par rotations. L'importance de cette interpr\'etation 
g\'eom\'etrique tient dans le fait remarquable que la construction garde
alors un sens dans un contexte extr\`emement plus g\'en\'eral que celui
des fibr\'es vectoriels sur des sch\'emas. Il suffit en effet 
d'avoir \`a disposition deux ingr\'edients principaux: l'existence
d'un objet des lacets $LX$, et l'existence de la trace
de l'automorphisme tautologique $m$. \\

La donn\'ee de base de notre situation abstraite consistera en deux
\'el\'ements.

\begin{enumerate}

\item Une notion \emph{d'objets g\'eom\'etriques} (e.g. sch\'emas, 
vari\'et\'es diff\'erentielles etc.), que nous organiserons
en une structure cat\'egorique $T$ du type topos. Nous allons
voir que pour que la construction soit pertinente il faut que
$T$ soit un $\s$-topos (au sens de \cite{lu1,tove2}; voir \S 1.5), et non seulement un topos au sens de Grothendieck.

\item Pour tout objet g\'eom\'etrique $X \in T$ du point $(1)$ ci-dessus, 
un notion de \emph{coefficients sur $X$} (e.g. fibr\'es vectoriels
alg\'ebriques ou diff\'erentiels, complexes parfaits \dots), que
nous organiserons en un foncteur $X \mapsto \mathcal{A}(X)$
du type champ en cat\'egories sur $T$. Nous allons
voir que pour que la construction soit pertinente il faut que
$\mathcal{A}(X)$ soit une $\s$-cat\'egorie mono\"\i dale (voir \S 2), 
et non seulement une cat\'egorie mono\"\i dale au sens usuel.

\end{enumerate}

Comme nous l'avons observ\'e il nous sera n\'ecessaire 
de construire pour $X \in T$ un objet $LX$, des lacets sur $X$.
Nous voulons que cet objet soit d\'efini par 
$LX=Map(S^{1},X)$, o\`u $S^{1}=B\mathbb{Z}$ est le cercle simplicial.
Or, \`a \'equivalence faible pr\`es, on a 
$$S^{1}\simeq *\coprod_{*\coprod *}*,$$ 
et donc
$LX\simeq X \times_{X\times X}X$,
o\`u ce produit fibr\'e est construit dans $T$. Nous voyons ici que si 
$T$ est un topos au sens usuel ce produit fibr\'e se simplifie, et on trouve
$LX\simeq X$, ce qui ne rend certainement pas la construction tr\'es int\'eressante.
Cependant, le produit fibr\'e ci-dessus devient non trivial lorsqu'il est
construit dans un contexte homotopique ou encore $\infty$-cat\'egorique
car, contrairement au cas des cat\'egories, le morphisme diagonal dans une cat\'egorie sup\'erieure n'est pas
un monomorphisme en g\'en\'eral. Dans le paragraphe pr\'ec\'edent, rappelant les id\'ees
de la construction,
$LX$ \'etait construit dans la th\'eorie homotopique, ou de mani\`ere
\'equivalente dans l'$\infty$-cat\'egorie, des sch\'emas d\'eriv\'es
de \cite{hagII,to3}, dans laquelle les auto-intersections donnent lieu
\`a des objets non triviaux. Ces sch\'emas d\'eriv\'es ne forment 
pas un topos de Grothendieck mais un $\infty$-topos
au sens de \cite{lu1,tove2} (voir aussi \cite{hagI}). Ainsi, $T$ 
sera pour nous un $\infty$-topos. Cela nous force aussi \`a supposer
que les coefficients  au-dessus d'un objet $X \in T$ 
forment eux-m\^eme une $\infty$-cat\'egorie et non pas seulement une
cat\'egorie. Par exemple, les fibr\'es vectoriels sur le sch\'ema
d\'eriv\'e $LX$ forment naturellent une $\infty$-cat\'egorie qui n'est pas
\'equivalente \`a une cat\'egorie en g\'en\'eral. Comme nous devons de plus pouvoir d\'efinir
la trace d'endomorphismes de coefficients, les $\infty$-cat\'egories
$\mathcal{A}(X)$ devront \^etre munies de \emph{structures mono\"\i dales
sym\'etriques rigides}. Ainsi, $\mathcal{A}$ sera pour nous
un champ sur $T$ et \`a valeurs dans les $\infty$-cat\'egories mono\"\i dales
sym\'etriques rigides. Un tel couple $(T,\mathcal{A})$ sera 
appel\'e un \emph{$\infty$-topos catannel\'e rigide}, et notre r\'esultat
principal est le fait qu'un caract\`ere de Chern peut toujours
\^etre construit pour un tel couple suivant les grandes lignes de
la construction que nous avons esquiss\'ee. Un point crucial
cependant, et qui rend ce travail non formel, est l'existence de morphismes de \emph{trace cylique}
pour les $\s$-cat\'egories mono\"\i dales sym\'etriques rigides (voir plus loin
dans cette introduction), qui comme nous l'expliquerons est un fait relativement
non trivial. \\

\noindent \textbf{Le th\'eor\`eme principal  --}
Le r\'esultat principal de ce travail peut se r\'esumer en l'\'enonc\'e suivant.

\begin{thm}\label{ti}
Soit $T$ un $\infty$-topos et $\mathcal{A}$ un 
champ en $\infty$-cat\'egories mono\"\i dales sym\'etriques rigides. Soit $X\in T$ et 
notons $LX=Map(S^{1},X) \in T$ son objet des lacets, muni de
l'action du groupe $S^{1}=B\mathbb{Z}$ naturelle.
Alors, il existe une application d'$\infty$-groupo\"\i des
$$\mathcal{A}(X)^{iso} \longrightarrow End_{\mathcal{A}(LX)}(\mathbf{1})^{hS^{1}},$$
entre l'$\infty$-groupo\"\i de sous-jacent \`a l'$\s$-cat\'egorie $\mathcal{A}(X)$ et 
l'$\infty$-groupo\"\i de des endomorphismes $S^{1}$-\'equivariants 
de l'objet unit\'e  de $\mathcal{A}(LX)$.
\end{thm}

Le morphisme du th\'eor\`eme pr\'ec\'edent sera appel\'e le \emph{pr\'e-caract\`ere
de Chern} de l'$\s$-topos catannel\'e rigide $(T,\mathcal{A})$, et 
le caract\`ere de Chern en sera une modification obtenue 
par faisceautisation (le pr\'e-caract\`ere de Chern se r\'ev\`ele parfois
plus int\'eressant que le caract\`ere de Chern lui-m\^eme, voir par 
exemple \S 4.2). \\

\noindent\textbf{Cons\'equences et applications --} 
Dans ce travail, nous proposons principalement deux applications du th\'eor\`eme pr\'ec\'edent, 
mais il est certain que l'on pourrait en trouver d'autres, ne serait-ce que 
par exemple en reprenant ces deux applications dans d'autres contextes
g\'eom\'etriques que celui de la g\'eom\'etrie alg\'ebrique d\'eriv\'ee
(e.g. celui de la g\'eom\'etrie diff\'erentielle d\'eriv\'ee de \cite{sp}, 
voir aussi \cite[\S 4]{lu3} pour plusieurs autres exemples
de \emph{g\'eom\'etries d\'eriv\'ees}).

Pour commencer, nous revenons sur la construction du caract\`ere de Chern
des fibr\'es vectoriels en g\'eom\'etrie alg\'ebrique que nous avons
pr\'esent\'ee en d\'ebut de cette introduction. Pour cela, 
nous prendrons pour $T$ le $\s$-topos $dAff_{k}^{\sim,et}$ des 
$k$-champs d\'eriv\'es, qui, par d\'efinition, est le $\s$-topos
des champs sur l'$\s$-site des $k$-alg\`ebres simpliciales
commutatives muni de la topologie \'etale (voir \S 4.1).
Pour $\mathcal{A}$ nous prendrons le champ des complexes parfaits, 
qui \`a une $k$-alg\`ebre commutative simpliciale $B$ associe
l'$\s$-cat\'egorie mono\"\i dale sym\'etrique $Parf(B)$ des 
$N(B)$-dg-modules parfaits (o\`u $N(B)$ est la dg-alg\`ebre 
obtenue \`a patrir de $A$ par normalisation). Pour tout
objet $X\in T$, $End_{\mathcal{A}(X)}(\mathbf{1})$ se d\'ecrit 
comme l'ensemble simplicial des fonctions sur $X$, que nous noterons
$\mathcal{O}(X)$. Le caract\`ere 
de Chern de notre th\'eor\`eme \ref{ti} induit ainsi un morphisme
$$Parf(X)^{iso} \longrightarrow \mathcal{O}(LX)^{hS^{1}},$$
de l'$\s$-groupo\"\i de des complexes parfaits sur le champ d\'eriv\'e $X$
vers l'$\s$-groupo\"\i de des fonctions $S^{1}$-\'equivariantes sur $LX$. 
Lorsque $X$ est un $k$-sch\'ema ceci redonne, modulo l'identification
entre $\mathcal{O}(LX)^{hS^{1}}$ et l'homologie cyclique n\'egative de 
$X$, le caract\`ere de Chern des complexes parfaits (voir Thm. \ref{compfibvect} et Rem. \ref{compparf}). Lorsque $X$ est un 
$k$-champ alg\'ebrique (disons au sens d'Artin) non deriv\'e, $\mathcal{O}(LX)^{hS^{1}}$
est une version \emph{orbifold} de l'homologie cyclique (n\'egative) de $X$
et le caract\`ere de Chern ainsi obtenu est donc une extension au cas
des champs d'Artin 
du caract\`ere de Chern d\'efini par exemple dans \cite{to,jkk}.

Pour notre second contexte d'application $T$ est encore
l'$\s$-topos des $k$-champs d\'eriv\'es mais o\`u 
$\mathcal{A}$ est maintenant un certain champ des dg-cat\'egories
compactement engendr\'ees, que nous noterons $\mathbb{D}g$. 
La construction pr\'ecise du champ
$\mathbb{D}g$ est assez technique et sera donn\'ee en \S 4.3.
Signalons simplement que
pour $X=Spec\, A$, avec $A$ une $k$-alg\`ebre commutative simpliciale, 
$\mathbb{D}g(X)$ est un mod\`ele \`a l'$\s$-cat\'egorie dont
les objets sont des dg-cat\'egories 
$A$-lin\'eaires et petites, et dont les morphismes sont les morphismes Morita
(au sens de \cite{to1}).
Pour $X\in T$ les objets de $\mathbb{D}g(X)$ doivent \^etre vus comme
des familles de dg-cat\'egories cocompl\`etes param\'etr\'ees par $X$ et qui, localement pour
la topologie \'etale sur $X$, poss\`edent des g\'en\'erateurs compacts. Ainsi, 
$End_{\mathbb{D}g(X)}(\mathbf{1})$ n'est autre 
que l'$\s$-groupo\"\i de des complexes quasi-coh\'erents
sur $X$ (qui sont les endomorphismes Morita de la dg-cat\'egorie
unit\'e sur $X$).
Notre caract\`ere de Chern dans ce contexte fournit un morphisme
$$\mathbb{D}g(X)^{iso} \longrightarrow
QCoh^{S^{1}}(LX),$$
o\`u le membre de droite d\'esigne l'$\s$-groupo\"\i de
des complexes quasi-coh\'erents et $S^{1}$-\'equivariants sur 
$LX$.
Ainsi, le caract\`ere de Chern d'une famille de dg-cat\'egories 
param\'etr\'ees par $X$ est un objet dans la cat\'egorie d\'eriv\'ee
$S^{1}$-\'equivariante des quasi-coh\'erents sur $LX$, objet qu'il 
faut voir comme un \emph{faisceau d'anomalies}, au sens
o\`u cette expression est utilis\'ee par exemple
dans \cite[\S 6.2]{bry}. Lorsque 
$X$ est un $k$-sch\'ema lisse de caract\'eristique nulle, il existe une
relation entre $QCoh^{S^{1}}(LX)$ et certains $\mathcal{D}_{X}$-modules
filtr\'es sur $X$ (\cite{d&d}), et \`a travers cette identification 
le caract\`ere de Chern d'une dg-cat\'egorie sur $X$ fournit 
la donn\'ee d'une \emph{connexion de Gauss-Manin}, et d'une
\emph{filtration de Hodge} sur le complexe d'homologie 
p\'eriodique associ\'e (voir \S 4.3). Cette construction est une 
g\'en\'eralisation de \cite{ge} et \cite{dotats}, et parrrait nouvelle dans un tel
degr\'e de g\'en\'eralit\'e. Lorsque maintenant 
$X=BG$ est le champ classifiant d'un $k$-sch\'ema en groupes, 
un objet de $\mathbb{D}g(X)$ peut \^etre consid\'er\'e
comme une \emph{repr\'esentation dg-cat\'egorique de $G$}
dans le style de \cite{gai}. Notre th\'eor\`eme \ref{ti} construit pour 
une telle repr\'esentation son \emph{faisceau des caract\`eres}, qui est 
un complexe quasi-coh\'erent $S^{1}$-\'equivariant sur le champ
$[G/G]$, et donc en particulier un objet dans $D_{qcoh}^{G}(G)$, la
cat\'egorie d\'eriv\'ee $G$-\'equivariante de $G$. La construction
de ce faisceau des caract\`eres nous semble, elle aussi, nouvelle, et 
g\'en\'eralise au cas des dg-cat\'egories la construction du 
caract\`ere $2$-cat\'egorique de \cite{gaka}.

Enfin, dans \S 4.4 nous nous int\'eressons \`a un sous-champ
de $\mathbb{D}g$, form\'e des dg-cat\'egories
\emph{satur\'ees} (aussi appel\'ees propres et lisses), et not\'e
$\mathbb{D}g^{sat}$. Nous observons que le
caract\`ere de Chern pour $\mathbb{D}g$, appliqu\'e \`a un objet
satur\'e, founit un complexe quasi-coh\'erent et $S^{1}$-\'equivariant sur $LX$
dont le complexe sous-jacent est parfait. Ainsi, en composant avec 
le caract\`ere de Chern pour les complexes parfaits (sur $[LX/S^{1}]$)
on obtient une fonction $S^{1}\times S^{1}$-\'equivariante sur
$L^{(2)}X=Map(S^{1}\times S^{1},X)$, l'espace des lacets doubles
sur $X$:
$$\mathbb{D}g^{sat}(X)^{iso} \longrightarrow \mathcal{O}(L^{(2)}X)^{h(S^{1}\times S^{1})}.$$
Nous appelons le membre de droite l'\emph{espace d'homologie cyclique secondaire de $X$} et
nous le proposons comme receptacle pour le caract\`ere de Chern des familles de dg-cat\'egories
satur\'ees. C'est l'existence de ce dernier morphisme qui \'etait propos\'e comme
caract\`ere de Chern des familles de dg-cat\'egories satur\'ees dans \cite{tove}, et dont
la preuve de l'existence n'\'etait qu'esquiss\'ee. \\

\noindent \textbf{Contenu du travail --} 
La premi\`ere section pr\'esente des notions et des r\'esultats
g\'en\'eraux de la th\'eorie des $\s$-cat\'egories (qui seront
pour nous des $1$-cat\'egories de Segal au sens de \cite{be,hisi}) qui sont indispensables
au reste du travail.
Certains sont bien connus et apparaissent d\'ej\`a dans des travaux
ant\'erieurs, comme par exemple les notions d'adjonctions, de limites, 
de localisation que l'on trouve dans \cite{hisi}. D'autres le sont moins, comme les notions 
d'$\s$-cat\'egories (co)fibr\'ees, d'$\s$-cat\'egories localement pr\'esentables et 
d'$\s$-topos, bien que celles-ci se trouvent \^etre trait\'ees dans
le contexte des quasi-cat\'egories dans \cite{lu1}. Nous avons essay\'e
cependant de rendre cette section le plus ind\'ependante possible 
d'autres r\'ef\'erences, et mis \`a part quelques r\'esultats fondamentaux
(comme par exemple l'existence de la structure de mod\`eles sur
les cat\'egories de Segal, ou encore l'\'equivalence avec 
la th\'eorie des cat\'egories simpliciallement enrichies) nous donnons
des arguments complets pour tous les r\'esultats que nous
utiliserons par la suite.

Dans la seconde partie de ce travail nous pr\'esentons une th\'eorie
des $\s$-cat\'egories mono\"\i dales sym\'etriques. Cette th\'eorie
est directement tir\'ee du texte \cite{to0}, et est \'equivalente \`a celle
pr\'esent\'ee dans le cadre des quasi-cat\'egories dans \cite{lu2}. On y introduit
la notion d'$\s$-cat\'egories mono\"\i dales sym\'etriques
rigides, et on montre comment on peut construire la trace 
d'un endomorphisme dans une telle $\s$-cat\'egorie. Une premi\`ere
difficult\'e technique est la construction 
d'une trace avec suffisament de propri\'et\'es de fonctorialit\'e, que l'on
obtient en d\'emontrant l'existence, et l'unicit\'e, d'un mod\`ele universel (Thm. \ref{t1},
Prop. \ref{p6}). Une seconde difficult\'e est de produire une
\emph{trace cyclique}, c'est \`a dire poss\`edant des propri\'et\'es 
d'invariances par rotations dans des suites cycliques d'\'equivalences. 
Cette seconde difficult\'e est r\'esolue en appliquant un th\'eor\`eme
r\'ecent de J. Lurie et M. Hopkins donnant une r\'eponse positive \`a l'\emph{hypoth\`ese du cobordisme}, 
qui pr\'edit le caract\`ere universel de la $\s$-cat\'egorie des bordismes
orient\'es de dimension $1$ (voir \cite{lu2}). Ce th\'eor\`eme permet de montrer 
la contractibilit\'e de l'espace des traces, dont la cons\'equence
imm\'ediate est l'existence et l'unicit\'e d'une trace cyclique (voir
Thm. \ref{t2}). Finalement nous montrons que la trace cyclique 
est additive et multiplicative. 

Notre troisi\`eme section est consacr\'ee \`a la notion d'\emph{$\s$-topos
catannel\'es rigides}, et \`a la construction du caract\`ere de Chern 
pour de tels objets. Il s'agit d'utiliser avec profit les r\'esultats et 
constructions des deux sections pr\'ec\'edentes afin de d\'emontrer 
le th\'eor\`eme \ref{ti}. 

La derni\`ere section pr\'esente quelques
contextes d'applications de la construction dont les principaux ont d\'ej\`a \'et\'e
mentionn\'es dans cette introduction; en particulier, on donne ici la comparaison 
entre notre caract\`ere de Chern et l'usuel, dans le cas des fibr\'es vectoriels sur des vari\'et\'es
projectives lisses sur un corps de caract\'eristique nulle (Th\'eor\`eme \ref{compfibvect}). 

Enfin, un appendice pr\'esente la 
preuve certains r\'esultats que nous utilisons sur les cat\'egories
de simplexes ainsi que la preuve d\'etaill\'ee de la comparaison avec les
caract\`ere de Chern usuel. \\

\noindent \textbf{Questions non trait\'ees --} 
Nous souhaitons signaler que plusieurs questions de comparaison,
qui apparaissent naturellement dans les exemples de la section \S 4,
ont \'et\'e laiss\'ees ouvertes. Par exemple, pour comparer rigoureusement
le caract\`ere de Chern construit en \S 4.2 et le caract\`ere de Chern 
usuel (\`a valeurs dans l'homologie cyclique) il faut commencer par 
comparer l'homologie cyclique n\'egative d'un sch\'ema $X$ et 
les fonctions $S^{1}$-\'equivariantes sur $LX$. M\^eme lorsque 
$X$ est un sch\'ema affine cette comparaison, bien que consid\'er\'ee
comme folklorique (voir par exemple \cite{bena}), ne semble pas
avoir \'et\'e \'ecrite en d\'etails. M. Hoyois a attir\'e notre attention
sur le fait qu'il ne parrait pas clair que la litt\'erature
actuelle suffise pour mener \`a bien cette comparaison et il est possible
qu'un travail non trivial soit n\'ecessaire pour parvenir \`a  \'ecrire une preuve compl\`ete. 
Dans le cas o\`u $X$ est un sch\'ema sur un anneau de caract\'eristique nulle, 
une comparaison directe entre $\mathcal{O}(LX)^{hS^{1}}$ et la 
cohomologie de de Rham de $X$ est cependant possible, sans passer par 
une comparaison avec l'homologie cyclique n\'egative. C'est ce qui
est fait dans \cite{bena} dans le cas lisse, et dans \cite{rhamloop}
en g\'en\'eral. Dans
\S 4.3 nous avons laiss\'e de cot\'e une question de m\^eme nature concernant
la comparaison entre $\mathcal{D}_{X}$-modules sur un sch\'ema lisse 
de caract\'eristique nulle et complexes quasi-coh\'erents et 
$S^{1}$-\'equivariants sur $LX$. Cette comparaison est trait\'ee
dans \cite{bena}, et peut aussi se d\'eduire
de \cite{rhamloop} afin de traiter aussi le cas non-lisse. 

Dans
\S 4.4 nous introduisons 
la notion d'\emph{homologie cyclique secondaire}  qui est le receptacle final du caract\`ere de Chern 
des familles de dg-cat\'egories satur\'ees. Nous n'avons pas encore une d\'escription explicite de cette homologie secondaire, sauf dans le cas du corps de base. Il semble probable que des applications successives du th\'eor\`eme
de Hochschild-Kostant-Rosenberg permettent de la calculer dans le cas lisse en caract\'eristique nulle, mais nous avons pour les moment laiss\'e cette question, \`a la quelle nous pensons cependant qu'il est important de r\'epondre.\\ Enfin, et pour finir avec cette introduction, 
nous avons signal\'e, en d\'ebut d'introduction, un contexte possible d'application de
notre construction, \`a savoir celui de certaines g\'eom\'etries non additives,
comme par exemple celui des sch\'emas au-dessus du corps \`a un \'el\'ement ou plus
g\'en\'eralement ceux apparaissant dans \cite{tova}. Nous n'avons pas inclu
ces contextes comme exemples dans le \S 4, mais nous pensons cependant qu'il 
s'agit d'applications \'eventuelles dignes d'int\'er\^ets. \\

\noindent \textbf{Remerciements.} 
Nous remercions tres particuli\`erement J. Lurie pour nous avoir expliqu\'e les r\'esultats
de \cite{lu2} pendant leur \'elaboration. 
Nous remercions aussi D. Ben-Zvi, M. Hoyois, D. Kaledin, D. Nadler, T. Pantev, D. Sullivan, M. Vaqui\'e pour des 
discussions sur le sujet de ce travail qui nous ont beaucoup apport\'e. Notre reconaissance va aussi \`a tous les participants au ``Cours de Geom\'etrie Alg\'ebrique Deriv\'e dans la Datcha'' (Moscou, Juillet 2009) pour leur int\'er\^et, leur g\'en\'erosit\'e et par les \'echanges fructueux d'id\'ees.

\bigskip
\begin{center}{\textbf{---}} \end{center}
\medskip

\noindent \textbf{Notations --} Nous utiliserons les d\'efinitions de cat\'egories
de mod\`eles de \cite{ho}. Les \emph{mapping spaces}, que nous appellerons
comme il se doit \emph{espaces de morphismes}, d'une cat\'egorie de mod\`eles
$M$ seront not\'es $Map_{M}$, ou plus simplement
$Map$ s'il n'y a pas lieu de sp\'ecifier $M$ (voir \cite[\S 5]{ho}).
Les produits fibr\'es homotopiques dans une cat\'egorie de mod\`eles
seront not\'es $-\times_{-}^{h}-$.

Dans ce travail, le symb\^ole $\Delta^{n}$ ne d\'esignera pas le
simplexe standard mais la cat\'egorie 
$$\xymatrix{
0 \ar[r] & 1 \ar[r] & \dots \ar[r] & n,}$$
classifiant les chaines composables de $n$-morphismes. Le simplexe
standard sera not\'e quand \`a lui $\overline{\Delta^{n}}$. 

\section{$\s$-Cat\'egories}

Cette premi\`ere section contient les d\'efinitions et \'enonc\'es
de la th\'eorie des $\s$-cat\'egories que nous 
utiliserons dans ce travail.  Nous n'aborderons pas, ou tr\`es peu, les
aspects $(2,\infty)$-cat\'egoriques, et dans ce travail
l'expression \emph{$\infty$-cat\'egorie} sera synonyme de 
\emph{$(1,\infty)$-cat\'egorie}, c'est \`a dire d\'esignera 
une $\s$-cat\'egorie dont tous les $i$-morphismes sont 
inversibles pour $i>1$ (voir \cite{be}). Nous avons opt\'e pour la th\'eorie
des cat\'egories de Segal comme mod\`ele \`a celle
des $\s$-cat\'egories, mais tout autre mod\`ele ad\'equat pourrait
faire l'affaire. Une des motivations pour avoir utilis\'e 
l'expression \emph{$\s$-cat\'egorie} au lieu de 
\emph{cat\'egorie de Segal} est pr\'ecis\'ement de permettre
au lecteur d'adpater ce texte \`a d'autre choix de th\'eorie
des $(1,\s)$-cat\'egorie.

Pour ce qui est des r\'ef\'erences, une partie important 
des notions et \'enonc\'es peuvent se trouver dans 
le texte fondateur \cite{hisi}. Certaines notions plus avanc\'ees
sont tir\'ees de \cite{tove2}, et d'autres sont des adaptations
de r\'esultats de \cite{lu1} au cadre des cat\'egories de Segal. 
Aucun des r\'esultats de cette section ne pr\'etend r\'eellement \`a 
l'originalit\'e, m\^eme si certains ne semblent pas
se trouver dans la litt\'erature, comme par exemple le traitement 
des cat\'egories de Segal (co)fibr\'ees. 

\subsection{Th\'eorie homotopique des $\s$-cat\'egories}

Tout au long de ce travail nous adopterons la d\'efinition fondamentale suivante.

\begin{df}\label{d-1}
Une \emph{$\s$-cat\'egorie} est une (1-)cat\'egorie de Segal au sens de \cite[\S 2]{hisi} (voir aussi
\cite[Def. 3.14]{be}).
\end{df}

Pour fixer les notations, rappelons qu'une cat\'egorie de Segal est un 
foncteur $A : \Delta^{op} \longrightarrow \mathbf{SEns}$ v\'erifiant les deux conditions suivantes.
\begin{enumerate}
\item $A_{0}$ est un ensemble simplicial discret (i.e. isomorphe \`a un ensemble
simplicial constant).
\item Pour tout $n$ le morphisme de Segal
$$A_{n} \longrightarrow A_{1}\times _{A_{0}}A_{1} \times_{A_{0}} \dots \times_{A_{0}}A_{1}$$
est une \'equivalence faible d'ensembles simpliciaux.
\end{enumerate}

Rappelons aussi qu'une $\s$-cat\'egorie est \emph{stricte} (nous dirons aussi  
\emph{$\mathbb{S}$-cat\'egorie}), 
si les morphismes du point $(2)$ ci-dessus sont des isomorphismes. Tel est le cas si et seulement 
si $A$ est le nerf d'une cat\'egorie simplicialement enrichie. Nous identifierons toujours
une cat\'egorie simplicialement enrichie \`a son nerf et ainsi \`a la 
cat\'egorie de Segal stricte correspondante. Une
cat\'egorie sera toujours vue comme une cat\'egorie simplicialement enrichie
dont les ensembles simpliciaux de morphismes sont discrets, et 
nous verrons ainisi la th\'eorie des cat\'egories comme une sous-th\'eorie
de celle des $\s$-cat\'egories.

Pour une 
$\s$-cat\'egorie $A$ nous noterons encore $A$ son ensemble d'objets $A_{0}$.
Pour un $n\geq 1$, on dispose d'une d\'ecomposition canonique
$$A_{n}\simeq \coprod_{(x_{0}, \dots, x_{n}) \in A^{n+1}}A(x_{0}, \dots, x_{n}),$$
o\`u $A(x_{0}, \dots, x_{n})$ d\'esigne la fibre de la projection
$A_n \longrightarrow A_{0}^{n+1}$ prise en le point $(x_{0}, \dots, x_{n})$.
Les ensembles simpliciaux $A(x_{0}, \dots, x_{n})$ sont par d\'efinition les espaces
des chaines
composables de morphismes $\xymatrix{x_{1} \ar[r] & x_{2} \ar[r] & \dots \ar[r] & x_{n}}$. 
En particulier, pour $x,y$ deux objets, $A(x,y)$ sera l'espace des morphismes de $x$ vers $y$.
Dans le cas o\`u $A$ est stricte, on a des isomorphismes
$$A(x_{0}, \dots, x_{n}) \simeq A(x_{0},x_{1}) \times A(x_{1},x_{2}) \times \dots \times A(x_{n-1},x_{n}).$$
Il nous arrivera aussi d'utiliser la notation suivante
$$Map_{A}(x,y):=A(x,y),$$
ou encore $Map(x,y):=A(x,y)$
lorsqu'il n'y aura pas d'ambiguit\'e sur $A$.  
Notons que la composition des morphismes dans une $\s$-cat\'egorie $A$ n'est pas un morphisme bien d\'efini dans 
$\mathbf{SEns}$, mais est repr\'esentable par le zig-zag
$$\xymatrix{
A(x,y)\times A(y,z) & \ar[l]_-{\sim} \ar[r] A(x,y,z) & A(x,z).}$$
On dispose donc toujours d'un morphisme
$$-\circ - : A(x,y)\times A(y,z) \longrightarrow A(x,z)$$
bien d\'efini dans $\mathrm{Ho}(\mathbf{SEns})$.  Ceci nous permettra de dire qu'un morphisme $u\in A(x,y)$
(par cela on entend un $0$-simplexe de $A(x,y)$) est une \'equivalence dans $A$ 
si pour tout $z$ le morphisme induit
$$- \circ u : A(y,z) \longrightarrow A(x,z)$$
est un isomorphisme dans $\mathrm{Ho}(\mathbf{SEns})$. 
 
Pour deux $\s$-cat\'egories $A$ et $B$, un $\s$-foncteur $f : A \longrightarrow B$ est par d\'efinition
un morphisme de cat\'egories de Segal (c'est \`a dire un morphisme d'ensembles bisimpliciaux). 
Nous parlerons aussi de morphismes de $\s$-cat\'egories.
Les $\s$-cat\'egories, \'el\'ements d'un univers fix\'e, et les $\s$-foncteurs forment une
cat\'egorie not\'ee $\s -\mathbf{Cat}$. Si $\mathbb{U}$ est un univers nous noterons, lorsque 
une telle pr\'ecision se r\'ev\`ele \^etre n\'ecessaire, $\s -\mathbf{Cat}_{\mathbb{U}}$ la cat\'egorie
des $\s$-cat\'egories \'el\'ements de $\mathbb{U}$ (nous dirons aussi $\mathbb{U}$-petites).

Pour un $\s$-cat\'egorie $A$ nous noterons $[A]$ la cat\'egorie des composantes
connexes de $A$, aussi appel\'ee la cat\'egorie homotopique de $A$. Ses objets
sont les objets de $A$ et ses ensembles de morphismes sont d\'efinis par 
$[A](x,y):=\pi_{0}(A(x,y))$.  La composition des morphismes est donn\'ee par 
le diagramme
$$\xymatrix{
\pi_{0}(A(x,y)) \times \pi_{0}(A(y,z)) & \ar[l]_-{\sim} \ar[r] \pi_{0}(A(x,y,z)) & 
\pi_{0}(A(x,z)),}$$
induit par les identit\'es simpliciales. Les ensembles de morphismes dans $[A]$ seront aussi parfois 
not\'es 
$$[x,y]=[A](x,y).$$
Notons au passage qu'un morphisme de $A$ est une \'equivalence si et seulement si
son image dans $[A]$ est un isomorphisme.
La construction $A \mapsto [A]$ d\'efinit bien \'evidemment un foncteur
$$[-] : \s -\mathbf{Cat} \longrightarrow \mathbf{Cat}.$$
Un $\s$-foncteur $f : A \longrightarrow B$ est pleinement fid\`ele (resp. essentiellement surjectif)
si pour tout $x,y \in A$ le morphisme $A(x,y) \longrightarrow B(f(x),f(y))$ est 
une \'equivalence faible d'ensembles simpliciaux 
(resp. si le foncteur induit $[f] : [A] \longrightarrow [B]$ est essentiellement surjectif). 
Un $\s$-foncteur est une \'equivalence s'il est \`a la fois essentiellement surjectif et pleinement fid\`ele.
La cat\'egorie
homotopique de $\s -\mathbf{Cat}$ obtenue en inversant les \'equivalences sera not\'ee $\mathrm{Ho}(\s -\mathbf{Cat})$. 

La cat\'egorie $\mathrm{Ho}(\s -\mathbf{Cat})$ est en r\'ealit\'e \'equivalente \`a la cat\'egorie homotopique 
d'une cat\'egorie de mod\`eles $\s -\mathbf{Cat}^{pr}$ des \emph{pr\'e-cat\'egories de Segal}, contenant 
$\s -\mathbf{Cat}$ comme sous-cat\'egorie pleine. 
La cat\'egorie $\s -\mathbf{Cat}^{pr}$ est la sous-cat\'egorie pleine des 
foncteurs $\Delta^{op} \longrightarrow \mathbf{SEns}$ qui sont discrets en $[0]$.
 Cette cat\'egorie
est munie d'une structure de mod\`eles pour laquelles les cofibrations sont les
monomorphismes et la restriction des \'equivalences \`a $\s -\mathbf{Cat}$ sont les \'equivalences
d\'ecrites ci-dessus. La d\'efinition g\'en\'erale des \'equivalences est plus subtile et 
ne sera pas utilis\'ee dans ce travail (voir \cite[\S 5]{be} pour les d\'etails). Les objets fibrants 
de $\s -\mathbf{Cat}^{pr}$ sont des cat\'egories de Segal et il existe donc des inclusions
$$(\s -\mathbf{Cat}^{pr})^{fib} \subset \s -\mathbf{Cat} \subset \s -\mathbf{Cat}^{pr}.$$
Ces inclusions sont strictes et une cat\'egorie de Segal est fibrante si et seulement si
elle est fibrante au sens de Reedy comme objet simplicial de $\mathbf{SEns}$ (voir \cite[Cor.  5.13]{be}). 
Une propri\'et\'e fondamentale de $\s -\mathbf{Cat}^{pr}$ est qu'il s'agit d'une
cat\'egorie de mod\`eles interne, c'est \`a dire qu'elle est mono\"\i dale
au sens de \cite[\S 4]{ho} 
pour la structure mono\"\i dale induite par le produit direct
(voir aussi \cite[Def. 2.1]{to2}). Une cons\'equence de cela est l'existence
de Hom internes entre objets fibrants, et en particulier l'existence de Hom internes sur
la cat\'egorie homotopique $\mathrm{Ho}(\s -\mathbf{Cat})$. 
Par la suite nous utiliserons aussi l'expression \emph{$\s$-foncteur} pour d\'esigner un morphisme
dans $\mathrm{Ho}(\s -\mathbf{Cat})$. Cet abus de language est relativement anodin, et signifie que l'on prendra
soins d'effectuer des remplacements fibrants lorsque cela est n\'ecessaire. Cependant, 
ces remplacements seront, la plupart du temps, implicites, et ne seront mentionn\'es que lorsque
l'argument m\'erite un certain degr\'e de pr\'ecision. 

Pour un univers
$\mathbb{U}$ la cat\'egorie de mod\`eles $\s -\mathbf{Cat}_{\mathbb{U}}$ est une cat\'egorie
de mod\`eles $\mathbb{U}$-combinatoire au sens de \cite[App. A]{hagI}. De plus, lorsque $\mathbb{U} \in \mathbb{V}$ 
sont  deux univers fix\'es alors le foncteur d'inclusion induit un foncteur pleinement fid\`ele
$$\mathrm{Ho}(\s -\mathbf{Cat}_{\mathbb{U}}) \hookrightarrow \mathrm{Ho}(\s -\mathbf{Cat}_{\mathbb{V}}).$$
On peut montrer que l'image essentielle de ce foncteur d'inclusion consiste en les $\s$-cat\'egories
$A$ telles que, d'une part l'ensemble des classes d'isomorphismes de $[A]$ est 
en bijection avec un ensemble $\mathbb{U}$-petit, et d'autre part telles que pour tout $x,y\in A$ l'ensemble
simplicial $X$ est \'equivalent \`a un ensemble simplicial $\mathbb{U}$-petit. \\

La structure de la cat\'egorie $\mathrm{Ho}(\s -\mathbf{Cat})$ peut se d\'ecrire avec les techniques d\'evelop\'ees
dans \cite{to1}, correctement modifi\'es pour traiter le cas 
simplicial au lieu du cas \emph{dg}. 
Pour cela on fait appel \`a l'\'equivalence entre la th\'eorie
homotopique des $\s$-cat\'egories et celle des $\mathbb{S}$-cat\'egories (qui ne sont autre
que les $\s$-cat\'egories strictes). Nous
renvoyons \`a \cite{be} pour des d\'etails concernant cette \'equivalence, nous en retiendrons que le
foncteur d'inclusion $\mathbb{S}-\mathbf{Cat} \longrightarrow \s -\mathbf{Cat}$ induit une \'equivalence sur les
cat\'egories homotopiques
$$\mathrm{Ho}(\mathbb{S}-\mathbf{Cat}) \simeq \mathrm{Ho}(\s -\mathbf{Cat}).$$
Pour une $\s$-cat\'egorie $A$, le choix d'une $\mathbb{S}$-cat\'egorie munie $A'$ d'un
isomorphisme $A'\simeq A$ dans $\mathrm{Ho}(\s -\mathbf{Cat})$ sera souvent appel\'e \emph{un mod\`ele
strict de $A$}. Un tel mod\`ele existe et est toujours unique \`a isomorphisme unique pr\`es
dans $\mathrm{Ho}(\s -\mathbf{Cat})$. Les descriptions des espaces de morphimes et des Hom internes
donn\'ees par \cite{to1} fournissent ainsi une fa\c{c}on de d\'ecrire 
la cat\'egorie $\mathrm{Ho}(\s -\mathbf{Cat})$ que nous allons rappeler ci-dessous.

Comme nous l'avons signal\'e 
la cat\'egorie $\mathrm{Ho}(\s -\mathbf{Cat})$ est cart\'esiennement close. Ses $Hom$-internes
seront alors not\'es $\rh$, et sont d\'efinis comme 
$\underline{Hom}(,R(-))$, o\`u $R$ est un foncteur de remplacement fibrant 
dans $\s -\mathbf{Cat}^{pr}$ et $\underline{Hom}$ est le Hom interene de pr\'e-cat\'egories de Segal. 
Ces Hom internes sont compatibles avec la structure
simplicial au sens o\`u l'on dispose d'isomorphismes naturels dans $\mathrm{Ho}(\mathbf{SEns})$
$$Map(A,\mathbb{R}\underline{Hom}(B,C))\simeq Map(A\times B,C),$$
avec $Map$ les espaces de morphismes de la cat\'egorie de mod\`eles $\s -\mathbf{Cat}^{pr}$ (au sens
de \cite[\S 5]{ho}). 

Pour une cat\'egorie de mod\`eles simpliciale $M$ nous notons $Int(M)$ la $\mathbb{S}$-cat\'egorie
des objets fibrants et cofibrants dans $M$. Nous nous empresserons de voir 
$Int(M)$ comme une $\s$-cat\'egorie et donc comme un objet de $\s -\mathbf{Cat}$.
Si $M$ est de plus $\mathbb{U}$-combinatoire et 
$\mathbb{V}$-petite (avec des univers $\mathbb{U}\in \mathbb{V}$), alors pour toute
$\mathbb{S}$-cat\'egorie $\mathbb{U}$-petite $A$ on dispose d'un isomorphisme naturel dans 
$\mathrm{Ho}(\s -\mathbf{Cat}_{\mathbb{V}})$
$$\mathbb{R}\underline{Hom}(A,Int(M)) \simeq Int(M^{A}),$$
o\`u $M^{A}$ est la cat\'egorie des $A$-diagrammes (enrichis sur $\mathbf{SEns}$) 
\`a valeurs dans $M$ munie de sa structure de mod\`eles projective (fibrations et \'equivalences
d\'efinies termes \`a termes). De m\^eme, l'espace des morphismes
$Map(A,Int(M))$ est naturellement \'equivalent \`a $N(WM^{A})$, le nerf de la sous-cat\'egorie des \'equivalences
dans $M^{A}$. En particulier, l'ensemble $[A,Int(M)]$, des morphismes de $A$ vers 
$Int(M)$ dans $\mathrm{Ho}(\s -\mathbf{Cat})$ est en bijection avec l'ensemble des 
classes d'isomorphismes d'objets de $\mathrm{Ho}(M^{A})$. Comme toute $\s$-cat\'egorie est \'equivalente 
\`a une $\mathbb{S}$-cat\'egorie (et ce de fa\c{c}on unique \`a \'equivalence pr\`es) cela 
d\'ecrit les morphismes d'une $\s$-cat\'egorie vers une $\s$-cat\'egorie de la forme $Int(M)$.

Fixons $\Gamma^{*} \rightarrow id$, un foncteur de r\'esolution cofibrante dans la cat\'egorie
de mod\`eles $\s -\mathbf{Cat}^{pr}$ (au sens de \cite[\S 5]{ho}). On peut prendre par exemple le foncteur de r\'esolution 
naturel donn\'e par la structure simpliciale sur $\s -\mathbf{Cat}^{pr}$. On a alors pour une pr\'e-cat\'egorie
de Segal $A$
$\Gamma^{n}(A):=\overline{\Delta([n])}\times A$, o\`u 
$\overline{\Delta([n])}$ est le groupo\"\i de classifiant les chaines
de $n+1$ isomorphismes composables. On d\'efinit alors une adjonction de Quillen
$$\Pi_{\infty} :  \mathbf{SEns} \longleftrightarrow \s -\mathbf{Cat} : Hom(\Gamma^{*}(*),-).$$
L'adjoint \`a droite $\mathcal{I}$ envoie une $\s$-cat\'egorie $A$ sur 
l'ensemble simplicial $Hom(\Gamma^{*}(*),R(A))$, o\`u $R$ est un foncteur de remplacement fibrant. 
L'adjoint \`a gauche
$\Pi_{\infty}$ envoie un ensemble simplicial $K$ sur le coend du diagramme
$\left( (p,q) \mapsto \coprod_{K_{q}}\Gamma^{p}(*) \right)$. L'adjonction d\'eriv\'ee
$$\Pi_{\infty} : \mathrm{Ho}(\mathbf{SEns}) \longleftrightarrow \mathrm{Ho}(\s -\mathbf{Cat}) : \mathbb{R}Hom(\Gamma^{*}(*),-)=:\mathcal{I}$$
est telle que $\Pi_{\infty}$ soit pleinement fid\`ele d'image essentielle form\'ee des $\s$-cat\'egories
$T$ telles que $[T]$ soit un groupo\"\i de (que nous appellerons 
des \emph{$\s$-groupo\"\i des}). 
Cela peut en effet se d\'eduire de l'\'enonc\'e analogue
pour les espaces de Segal complets de (qui devient alors un \'enonc\'e formel) et des r\'esultats
de comparaisons de \cite{be}. Le foncteur $\Pi_{\infty}$ sera appel\'e le foncteur \emph{$\infty$-groupo\"\i de fondamental}, 
et $\mathcal{I}$ le foncteur \emph{espace sous-jacent}. Noter que pour une $\s$-cat\'egorie $A$ le
morphisme d'ajonction $\Pi_{\infty}(\mathcal{I}(A)) \longrightarrow A$ identifie 
le membre de gauche avec le sous-groupo\"\i de maximal de $A$. En d'autres termes on dispose d'un
carr\'e homotopiquement cart\'esien dans $\s -\mathbf{Cat}$
$$\xymatrix{
\Pi_{\infty}(\mathcal{I}(A)) \ar[r] \ar[d] & A \ar[d] \\
[A]^{iso} \ar[r] & [A],}$$
o\`u $[A]^{iso}$ est le sous-groupo\"\i de des isomorphismes dans $[A]$ (c'est aussi la
cat\'egorie homotopique de $\Pi_{\infty}(\mathcal{I}(A))$). Nous utiliserons aussi la notation
$T^{int} \subset T$ pour d\'esigner ce sous-groupo\"\i de maximal, 
o\`u \emph{int} fait r\'ef\'erence \`a \emph{int\'erieur} (voir \cite{hisi} dans lequel cette
terminologie est introduite). 

Le foncteur $\Pi_{\infty} : \mathrm{Ho}(\mathbf{SEns}) \longrightarrow \mathrm{Ho}(\s -\mathbf{Cat})$
est aussi isomorphe au foncteur qui envoie un ensemble simplicial fibrant $K$ sur 
$[n] \mapsto \underline{Hom}_{*}(\Delta^{n},K)$, o\`u 
$\underline{Hom}_{*}(\Delta^{n},K)$ d\'esigne l'ensemble simplicial des morphismes
$\Delta^{n} \rightarrow K$ qui fixent les sommets de $\Delta^{n}$. Ceci montre que le foncteur
$\Pi_{\infty}$ poss\`ede aussi un adjoint \`a gauche
$$|-| : \mathrm{Ho}(\s -\mathbf{Cat}) \longrightarrow \mathrm{Ho}(\mathbf{SEns})$$
qui envoie une $\s$-cat\'egorie sur sa r\'ealisation g\'eom\'etrique (i.e. la diagonale 
de l'ensemble bisimplicial correspondant). Le morphisme d'adjonction 
$A \longrightarrow \Pi_{\infty}(|A|)$ est alors la compl\'etion $\s$-groupo\"\i dale de
$A$.

\subsection{Adjonctions, limites et localisations}

La cat\'egorie homotopique des $\s$-cat\'egorie, $\mathrm{Ho}(\s -\mathbf{Cat})$, peut-\^etre 
promu en une $2$-cat\'egorie $\mathrm{Ho}_{\leq 2}(\s -\mathbf{Cat})$ dont la cat\'egorie
homotopique est naturellement \'equivalente \`a $\mathrm{Ho}(\s -\mathbf{Cat})$. Une m\'ethode
relativement directe consiste \`a consid\'erer $\s -\mathbf{Cat}^{f}$, la cat\'egorie
des $\s$-cat\'egories fibrantes. Cette cat\'egorie en cart\'esiennement close, et 
on peut donc la consid\'er\'ee comme une cat\'egorie enrichie dans
la cat\'egorie $\s -\mathbf{Cat}$. Notons cette $(\s -\mathbf{Cat})$-cat\'egorie
par $\underline{\s -\mathbf{Cat}^{f}}$, qui n'est autre qu'un mod\`ele pour la
$(2,\infty)$-cat\'egorie de Segal des $\s$-cat\'egories (voir \cite{hisi}). En rempla\c{c}ant 
les $\s$-cat\'egories de morphismes par leur cat\'egories homotopiques on trouve une
cat\'egorie enrichie dans $\mathbf{Cat}$, donc une $2$-cat\'egorie (stricte), que nous noterons
$\mathrm{Ho}_{\leq 2}(\s -\mathbf{Cat})$. La cat\'egorie homotopique de cette $2$-cat\'egorie
est naturellement isomorphe \`a la cat\'egories des $\s$-cat\'egories fibrantes 
et des classes d'homotopie de morphismes entre elles. 
On dispose donc bien d'une \'equivalence naturelle
$$\mathrm{Ho}(\s -\mathbf{Cat}) \simeq [\mathrm{Ho}_{\leq 2}(\s -\mathbf{Cat})].$$
Notons que par d\'efinition les objets de $\mathrm{Ho}_{\leq 2}(\s -\mathbf{Cat})$ 
sont les $\s$-cat\'egories fibrantes et la cat\'egorie des morphismes
entre $A$ et $B$  est $[\underline{Hom}(A,B)]$. Il est aussi possible
de transporter la structure $2$-cat\'egorique pour construire une
$2$-cat\'egorie \'equivalente, dont les objets sont toutes les $\s$-cat\'egories
et dont les cat\'egories de morphismes sont $[\underline{Hom}(R(A),R(B))]$, avec
$R$ un foncteur de remplacement fibrant. Dans ce qui suit
$\mathrm{Ho}_{\leq 2}(\s -\mathbf{Cat})$ d\'esignera l'une ou l'autre de ces constructions, et les
cat\'egories de morphismes seront simplement not\'ees $[\rh (A,B)]$. Il s'agit dans les
deux cas d'une $2$-cat\'egorie stricte.

Cet enrichissement de la cat\'egorie homotopique en une $2$-cat\'egorie 
permet de d\'efinir la notion d'adjonction de la fa\c{c}on suivante.
Nous dirons qu'un $\s$-foncteur $f : A \longrightarrow B$ poss\`ede
un adjoint \`a droite s'il existe un morphisme $g : B\longrightarrow A$ dans 
$\mathrm{Ho}_{\leq 2}(\s -\mathbf{Cat})$, et un objet $u \in [\rh (A,A)](id,gf)$, 
v\'erifiant la condition suivante: 
pour toute $\s$-cat\'egorie $C$ et tout objet $k \in [\mathbb{R}\underline{Hom}(C,A)]$,
$h \in [\mathbb{R}\underline{Hom}(C,B)]$ le morphisme induit
$$\xymatrix{
[fk,h] \ar[r]^-{g\circ } & [gfk,gh] \ar[r]^-{\circ u} & [k,gh]}$$
est bijectif (o\`u l'on a not\'e $[-,-]$ les ensembles de morphismes
dans $[\mathbb{R}\underline{Hom}(C,A)]$ et $[\rh (C,B)]$).
En d'autres termes, on demande \`a ce que pour toute $\s$-cat\'egorie
$C$, les foncteurs induits $f_{C} : [\rh (C,A)] \longleftrightarrow [\rh (C,B)] : g_{C}$, munis 
de la transformation naturelle $id \rightarrow g_{C}f_{C}$ induite par $u$,  d\'efinissent
une adjonction de cat\'egories au sens usuel.

Lorsque le morphisme $f$ poss\`ede un adjoint \`a droite le couple
$(g,u)$ est unique \`a isomorphisme unique pr\`es au sens suivant: si $(g',u')$ est un 
second couple satisfaisant \`a la m\^eme condition alors il existe un 
unique isomorphisme $\alpha : g \simeq g'$ dans $[\mathbb{R}\underline{Hom}(B,A)]$
tel que le diagramme suivant commute dans $[\mathbb{R}\underline{Hom}(A,A)]$
$$\xymatrix{
gf \ar[r]^-{u} \ar[d]_-{\alpha.f} & id \\
g'f \ar[ru]_-{u'}.}$$
On dispose dualement d'une notion de $\s$-foncteur poss\`edant un 
adjoint \`a gauche.  De plus si $(g,u)$ est un adjoint \`a droite de 
$f$ comme ci-dessus alors $(f,v)$ est un adjoint \`a gauche de $g$, o\`u 
$v \in [fg,id]$ correspond \`a l'identit\'e dans $[g,g]$
par la bijection $[fg,id] \simeq [g,g]$ d\'ecrite ci-dessus ($h=id$ et $k=g$).
Les couples $(f,u)$ et $(g,v)$ sont alors reli\'es par les relations triangulaires usuelles
$$(v.f)\circ (f.u)=id : f \longrightarrow f \qquad (g.v) \circ (u.g)=id : g \longrightarrow g,$$
qui sont des \'egalit\'es de morphismes dans $[\mathbb{R}\underline{Hom}(A,B)]$
et dans $[\mathbb{R}\underline{Hom}(B,A)]$ respectivement. Ces identit\'es triangulaires 
permettent de montrer que pour toute $\s$-cat\'egorie $C$ et tout 
$\s$-foncteurs $k \in [\mathbb{R}\underline{Hom}(C,A)]$,
$h \in [\mathbb{R}\underline{Hom}(C,B)]$ le morphisme dans $\mathrm{Ho}(\s -\mathbf{Cat})$
$$\xymatrix{
\rh (C,B)(fk,h) \ar[r]^-{g\circ } & \rh (C,A)(gfk,gh) \ar[r]^-{\circ u} & 
\rh (C,A)(k,gh)}$$
est un isomorphisme (alors que la d\'efinition d'adjoint ne donne, \`a priori, qu'une bijection sur les
$\pi_{0}$ de ses ensembles simpliciaux). En effet, un inverse est donn\'e par le morphisme
$$\xymatrix{
\rh (C,A)(k,gh) \ar[r]^-{f\circ } & \mathbb{R}\underline{Hom}(C,B)(fk,fgh) \ar[r]^-{\circ v} & 
\rh (C,A)(fk,h),}$$
les identit\'es triangulaires impliquant que ce morphisme est bien un inverse, dans $Ho(\mathbf{SEns})$, 
au morphisme ci-dessus. Ainsi, en prenant $C=*$ on trouve que pour tout $a\in A$ et tout $b\in B$, le morphisme 
naturel
$$B(f(a),b) \longrightarrow A(a,g(b))$$
est un isomorphisme dans $\mathrm{Ho}(\mathbf{SEns})$. \\

La notion d'ajoint ci-dessus permet de parler de limites et de colimites dans
les $\s$-cat\'egories. Soient $I$ et $A$ deux $\s$-cat\'egories, et consid\'erons
le morphisme diagramme constant
$$c : A \longrightarrow \rh (I,A)$$
adjoint de la projection $A\times I \longrightarrow A$. Nous dirons que $A$ poss\`ede
des limites (resp. des colimites) suivant $I$ si le foncteur $c$ poss\`ede un 
adjoint \`a droite (resp. \`a gauche). Un adjoint \`a droite sera not\'e
$lim_{I} : \rh (I,A) \longrightarrow A$, et un adjoint \`a gauche 
$colim_{I} :  \rh (I,A) \longrightarrow A$. 
Nous dirons plus g\'en\'eralement que $A$
poss\`ede des $\mathbb{U}$-limites (resp. des $\mathbb{U}$-colimites) si elle 
poss\`ede des limites le long de $I$ pour toute $\s$-cat\'egorie $I\in \mathbb{U}$. Si l'univers
$\mathbb{U}$ est implicite nous dirons simplement que $A$ poss\`ede des limites (resp. 
des colimites), sans plus de pr\'ecision. Lorsque $A$ poss\`ede des colimites, 
pour tout objet $a \in A$ et tout ensemble simplicial $K$ il existe un 
objet $K\otimes a \in A$. Cet objet vient avec un morphisme dans $\mathrm{Ho}(\mathbf{SEns})$
$$K \longrightarrow A(a,K\otimes a)$$
de sorte \`a ce que pour tout $b \in A$ le morphisme induit (d\'efini dans $\mathrm{Ho}(\mathbf{SEns})$)
$$K\times A(K\otimes a,b) \longrightarrow A(a,K\otimes a) \times A(K\otimes a,b) \longrightarrow
A(a,b)$$
induise un isomorphisme $A(K\otimes a,b) \simeq Map(K,A(a,b))$
dans $\mathrm{Ho}(\mathbf{SEns})$. Dualement, lorsque $A$ poss\`ede des limites alors 
pour tout $a\in A$ et tout ensemble simplicial $K$ il existe un objet
$a^{K}$. Cet objet vient avec un morphisme $K \longrightarrow A(a^{K},a)$
tel que pour tout $b\in A$ le morphisme induit
$$A(a,b^{K}) \longrightarrow Map(K,A(a,b))$$
soit un isomorphisme dans $\mathrm{Ho}(\mathbf{SEns})$. 

Supposons que $A$ et $B$ soient deux $\s$-cat\'egories poss\'edant des
limites le long d'une $\s$-cat\'egorie fix\'ee $I$. Alors, 
pour un $\s$-foncteur $f : A \longrightarrow B$ on dispose d'un 
morphisme naturel $\phi : f\circ lim_{I} \longrightarrow lim_{I}\circ f^{I}$ dans
la cat\'egorie $[\rh (\rh (I,A),B)]$ (o\`u $f^{I}$ 
d\'esigne le morphisme induit $\rh (I,A) \longrightarrow \rh (I,B)$). 
Ce morphisme est d\'efini de la fa\c{c}on suivante.
Notons $v : c\circ lim_{I} \longrightarrow id$ le morphisme d\'efinissant l'adjonction 
entre $c$ et $lim_{I}$. On dispose alors d'une application naturelle
$$\xymatrix{[c\circ lim_{I},id] \ar[r] & [f^{I}\circ c\circ lim_{I},f^{I}] \simeq
[c\circ f\circ lim_{I},f^{I}] \simeq [f\circ lim_{I},lim_{I}\circ f^{I}]}.$$
L'image de $v$ par cette application fournit le morphisme $\phi$ cherch\'e. Nous dirons alors
que le $\s$-foncteur $f$ commute aux $\mathbb{U}$-limites si le
morphisme ci-dessus $\phi$ est un isomorphisme pour toute $\s$-cat\'egorie $\mathbb{U}$-petite 
$I$. Dualement, nous dirons que $f$ commute aux $\mathbb{U}$-colimites si 
$f^{op} : A^{op} \longrightarrow B^{op}$ commute aux $\mathbb{U}$-limites. 
On v\'erifie qu'un $\s$-foncteur qui poss\`ede un adjoint \`a droite
commute toujours aux colimites. Dualement, un $\s$-foncteur qui poss\`ede un adjoint \`a gauche
commute toujours aux limites. \\

Supposons maintenant que $A$ soit une $\s$-cat\'egorie et que $S$ soit un 
ensemble de morphismes dans $[A]$. Une localisation de $A$ le long de $S$ est 
la donn\'ee d'un couple $(L_{S}A,l)$, form\'e d'une 
$\s$-cat\'egorie $L_{S}A$ et d'un $\s$-foncteur
$l : A \longrightarrow L_{S}A$ tel que pour toute $\s$-cat\'egorie $B$ le morphisme induit
$$l^{*} : \rh (L_{S}A,B) \longrightarrow \rh (A,B)$$
soit pleinement fid\`ele et que son image essentielle consiste en tous
les $\s$-foncteurs $f : A \longrightarrow B$ tel que $[f]$ envoie les morphismes de $S$ sur
des isomorphismes dans $[B]$. On montre qu'une localisation $(L_{S}A,l)$ existe toujours, 
et qu'elle est unique dans la cat\'egorie homotopique $\mathrm{Ho}(\s -\mathbf{Cat})$. 
Un mod\`ele explicite de $L_{S}A$ est donn\'e par la localisation simpliciale de 
Dwyer et Kan (voir \cite{dk}). On peut aussi construire $L_{S}A$ comme le push-out homotopique 
dans $\s -\mathbf{Cat}$ (voir \cite[\S 8.2]{to1})
$$\xymatrix{
\coprod_{S}\Delta^{1} \ar[r] \ar[d] & A \ar[d] \\
\coprod_{S} * \ar[r] & L_{S}A.}$$
Nous mettons en garde que ce push-out, pris dans $\s -\mathbf{Cat}^{pr}$ n'est pas
une cat\'egorie de Segal, et qu'il est n\'ecessaire de composer cette construction avec 
un foncteur de remplacement fibrant pour avoir une construction donnant 
une $\s$-cat\'egorie.
Dans les deux cas on voit que $l : A \longrightarrow L_{S}A$ est toujours
essentiellement surjectif.
Une cons\'equence de la propri\'et\'e universelle des localisations est que 
la localisation commute aux produits: pour $A$ et $B$ deux $\s$-cat\'egories et 
$S_{A}$ (resp. $S_{B}$) un ensemble de morphismes dans $[A]$ (resp. dans $[B]$), le morphisme
naturel $L_{S_{A}\times S_{B}}(A\times B) \longrightarrow L_{S_{A}}A \times L_{S_{B}}B$
est un isomorphisme dans $\mathrm{Ho}(\s -\mathbf{Cat})$. Un autre cons\'equence est 
que si $A'$ est une $\mathbb{S}$-cat\'egorie \'equivalente \`a $A$, alors 
la $\s$-cat\'egorie $L_{S}A$ est toujours \'equivalente 
\`a la sous-$\s$-cat\'egorie pleine de $Int(SPr(A'))$ form\'ee
des pr\'efaisceaux simpliciaux $F : (A')^{op} \longrightarrow \mathbf{SEns}$
qui envoient $S$ sur des \'equivalences.

Pour une cat\'egorie de mod\`eles $M$ nous noterons $L(M):=L_{W}M$
sa localis\'ee le long des \'equivalences faibles. Par d\'efinition
$[L(M)^{op},L(\mathbf{SEns})]$ est en bijection avec les classes d'isomorphismes
d'objets dans la cat\'egorie $\mathrm{Ho}(M^{\wedge})$ (c'est \` a dire avec les
classes d'isomorphismes dans $\mathrm{Ho}(SPr(M))$ des pr\'efaisceaux
qui envoient \'equivalences sur \'equivalences, voir \cite[\S 2.3.2]{hagI}). Ceci implique
que la sous-$\s$-cat\'egorie pleine de $Int(M^{\wedge})$ form\'ee
des objets \'equivalents \`a des repr\'esentables est naturellement \'equivalente
\`a $L(M)$. En particulier, le lemme de Yoneda de \cite[Thm. 4.2.3]{hagI} implique qu'il existe des isomorphismes
naturels dans $\mathrm{Ho}(\mathbf{SEns})$
$$Map_{M}(x,y) \simeq L(M)(x,y)$$
pour $x,y \in M$. Ainsi, si $M$ est une cat\'egorie de mod\`eles 
dont les ensembles de morphismes sont $\mathbb{U}$-petits alors 
$L(M)$ est \'equivalente \`a une $\s$-cat\'egorie dont les ensembles simplicaux
de morphismes sont de m\^emes $\mathbb{U}$-petits. Si maintenant, on suppose de plus
que $M$ est une cat\'egorie de mod\`eles simpliciales alors le plongement de Yoneda
simplicial induit une \'equivalence de $\s$-cat\'egories entre
$Int(M)$ et la sous-$\s$-cat\'egorie de $Int(M^{\wedge})$ form\'ee
des pr\'efaisceaux \'equivalents \`a des repr\'esentables. Ainsi, on a
$Int(M) \simeq L(M)$. Ceci implique que lorsque $M$ est $\mathbb{U}$-engendr\'ee par cofibration 
et simpliciale alors
pour toute
cat\'egorie $I$ qui est $\mathbb{U}$-petite le morphisme naturel
$$L(M^{I}) \longrightarrow \rh (I,L(M))$$
est un isomorphisme dans $\mathrm{Ho}(\s -\mathbf{Cat})$. En d'autres termes, 
la localisation commute aux Homs internes vers une $\s$-cat\'egorie provenant d'une
cat\'egorie de mod\`eles. 
Cet isomorphisme naturel s'\'etend \`a toute cat\'egorie de mod\`eles 
engendr\'ee par cofibrations et Quillen \'equivalente \`a une cat\'egorie de mod\`eles
simpliciale (par exemple \`a toute cat\'egorie de mod\`eles combinatoire). Ceci reste aussi 
vrai pour des sous-cat\'egories pleines de cat\'egories de mod\`eles qui sont 
stables par \'equivalences.

Cette compatibilit\'e de la localisation des cat\'egories de mod\`eles avec 
le passage aux cat\'egories de diagrammes implique aussi une compatibilit\'e 
entre adjonctions de Quillen et adjonctions de $\s$-cat\'egories. Pour 
voir cela, soit $f : M \longleftrightarrow N : g$ une adjonction de Quillen
entre cat\'egorie de mod\`eles simpliciales et $\mathbb{U}$-engendr\'ees par 
cofibrations. On d\'efinit un foncteur $\mathbb{L}f : LM \longrightarrow LN$ par 
le diagramme commutatif suivant
$$\xymatrix{
LM \ar[r]^-{\mathbb{L}f} \ar[d]_-{Q} & LN \\
LM^{c} \ar[r]_-{Lf^{c}} & LN^{c} \ar[u]}$$
o\`u $M^{c}$ est la sous-cat\'egorie pleine des objets cofibrants, 
$Q$ est un foncteur de remplacement cofibrant, et $Lf^{c}$ est le $\s$-foncteur
induit par $f^{c}$ (ce qui a un sens car $f^{c}$ pr\'eserve les \'equivalences). Dualement, 
on d\'efinit $\mathbb{R}g : LN \longrightarrow LM$ par le diagramme commutatif
$$\xymatrix{
LN \ar[r]^-{\mathbb{R}g} \ar[d]_-{R} & LM \\
LN^{f} \ar[r]_-{Lg^{f}} & LM^{f}, \ar[u]}$$
o\`u maintenant $N^{f}$ est la sous-cat\'egorie des objets fibrants et 
$R$ est un foncteur de remplacement fibrant. Ces constructions
d\'efinissent deux morphismes dans $\mathrm{Ho}_{\leq 2}(\s -\mathbf{Cat})$
$$\mathbb{L}f : LM \longrightarrow LN \qquad LM \longleftarrow LN : \mathbb{R}g.$$
Pour un objet $x\in M$, on dispose d'une chaine de morphismes dans $M$
$$gRfQ(x) \longleftarrow gfQ(x) \longleftarrow Q(x) \longrightarrow x,$$
et donc d'un diagramme
$$gRfQ \longleftarrow Q \longrightarrow id$$
d'endofoncteurs de $M$. Ces foncteurs pr\'eservant les \'equivalences on en d\'eduit 
un morphisme bien d\'efini dans $[\rh (LM,LM)]$
$$u : id \longrightarrow \mathbb{R}g \mathbb{L}f.$$
On v\'erifie alors que le couple $(\mathbb{R}g,u)$ est un adjoint \`a droite
de $\mathbb{L}f$ en remarquant que pour toute $\s$-cat\'egorie $C$, le couple induit sur 
$$\mathrm{Ho}(M^{C})\simeq [\rh (C,LM)] \longleftrightarrow [\rh (C,LN)] \simeq \mathrm{Ho}(N^{C})$$
n'est autre que le couple d\'efinissant l'adjonction d\'eriv\'ee
de l'adjonction de Quillen induite sur les cat\'egories de mod\`eles de diagrammes
$$f : M^{C} \longleftrightarrow N^{C} : g.$$
On d\'eduit de cela que pour toute cat\'egorie de mod\`eles simpliciales
et $\mathbb{U}$-engendr\'ee par cofibration $M$ la $\s$-cat\'egorie
$LM$ poss\`ede des $\mathbb{U}$-limites et des 
$\mathbb{U}$-colimites. En effet, pour une $\s$-cat\'egorie
$I$ qui est $\mathbb{U}$-petite, on consid\`ere le foncteur
$Colim_{I} : M^{I} \longrightarrow M$. Ce foncteur \'etant de
Quillen \`a gauche pour la structure de mod\`eles projective sur $M^{I}$
on en d\'eduit que le $\s$-foncteur induit
$$\mathbb{L}Colim_{I} : L(M^{I})\simeq \rh (I,LM) \longrightarrow LM$$
est un adjoint \`a gauche du $\s$-foncteur diagramme constant. De m\^eme, 
si $\widetilde{*}$ est un remplacement de l'objet final de $\mathbf{SEns}^{I}$ alors
le foncteur
$$\underline{Hom}(\widetilde{*},-) : M^{I} \longrightarrow M$$
est de Quillen \`a droite et le $\s$-foncteur induit 
$$L(M^{I})\simeq \rh (I,LM) \longrightarrow LM$$
est un adjoint \`a droite du foncteur diagramme constant (qui est 
le foncteur obtenu par localisation \`a partir de l'adjoint \`a gauche
$x \mapsto \widetilde{*} \otimes x$). Enfin, comme toute cat\'egorie 
de mod\`eles combinatoire est Quillen \'equivalente \`a une
cat\'egorie de mod\`eles simpliciales on voit aussi que 
$LM$ poss\`ede des $\mathbb{U}$-limites et des $\mathbb{U}$-colimites
pour toute cat\'egorie de mod\`eles $\mathbb{U}$-combinatoire $M$.

Terminons ce paragraphe par la remarque suivante: soit $M$ est une
cat\'egorie de mod\`eles combinatoire, $S$ est un ensemble
de morphismes dans $M$ et $L_{S}^{B}M$ la localisation de Bousfield \`a gauche
de $M$ le long de $S$. L'adjonction $id : L_{S}^{B}M \longleftrightarrow M : id$
induit alors une adjonction sur les $\s$-cat\'egories localis\'ees
$$i : L(L_{S}^{B}M) \longleftrightarrow LM : l.$$
Le $\s$-foncteur $l : LM \longrightarrow L(L_{S}^{B}M)$ induit alors un isomorphisme dans $\mathrm{Ho}(\s- \mathbf{Cat})$
$$L_{S}(LM) \simeq L(L_{S}^{B}M).$$

\begin{df}\label{d0}
La \emph{$\s$-cat\'egorie des $\s$-cat\'egories} est $L(\s -\mathbf{Cat})$. Elle sera not\'ee
$\uscat$ (ou encore $\uscat_{\mathbb{U}}$ si l'on se 
restreint aux $\s$-cat\'egories $\mathbb{U}$-petites).
\end{df}

Par d\'efinition on a une \'equivalence de cat\'egories
$$[\uscat]\simeq \mathrm{Ho}(\s -\mathbf{Cat}).$$
De plus, comme $\s -\mathbf{Cat}^{pr}$ est une cat\'egorie de mod\`eles simpliciale on a 
une \'equivalence de $\s$-cat\'egories
$$\uscat \simeq Int(\s -\mathbf{Cat}^{pr}).$$

\subsection{$\s$-Cat\'egories cofibr\'ees}

Pour ce paragraphe nous noterons $I$ une 
cat\'egorie (disons $\mathbb{U}$-petite pour fixer les id\'ees). On d\'efinit 
une adjonction 
$$\int_{I} : (\s -\mathbf{Cat}^{pr})^{I} \longleftrightarrow \s -\mathbf{Cat}^{pr}/I : Se_{I}$$
de la fa\c{c}on suivante. Pour $A \rightarrow I$ un objet 
de $\s -\mathbf{Cat}^{pr}/I$ on pose
$$Se_{I}(A):=\underline{Hom}_{I}(-/I,A) : I \longrightarrow \s -\mathbf{Cat}^{pr}.$$
Dans cette notation $-/I$ est le foncteur $I \longrightarrow \mathbf{Cat}$
qui envoie $i\in I$ sur $i/I$, et $\underline{Hom}_{I}$ sont 
les $Hom$ enrichis de $\s -\mathbf{Cat}^{pr}/I$ \`a valeurs
dans $\s -\mathbf{Cat}^{pr}$ (naturellement induits par ceux 
de $\s -\mathbf{Cat}^{pr}$). Le foncteur
$Se_{I}$ poss\`ede un adjoint \`a gauche
$\int_{I}$. Cet adjoint est caract\'eris\'e par 
le fait qu'il commute aux colimites, \`a l'enrichissement
dans $\s -\mathbf{Cat}^{pr}$, et par le fait que
$$\int_{I}(h^{i})=i/I \rightarrow I,$$
o\`u $h^{i}$ est le pr\'efaisceau d'ensembles corepr\'esent\'e
par l'objet $i$. De mani\`ere plus explicite, pour 
$F : I \longrightarrow \s -\mathbf{Cat}^{pr}$ un pr\'efaisceau
en pr\'e-cat\'egories de Segal, on d\'efinit
une pr\'e-cat\'egorie de Segal $\int_{I}F$ dont les
objets les couples $(i,x)$, avec $i\in I$ et $x$ un objet 
de $F(i)$. Pour $(i_{1},x_{1})$, \dots $(i_{n},x_{n})$ de tels
objets on pose
$$(\int_{I}F)((i_{1},x_{1}), \dots, (i_{n},x_{n})):=
\coprod_{\xymatrix{i_{1} \ar[r] & \dots \ar[r] & i_{n}}}F(i_{n})(x_{1}, \dots, x_{n}),$$
o\`u la somme est prise sur tous les diagrammes $\xymatrix{i_{1} \ar[r] & \dots \ar[r] & i_{n}}$
dans $ I$, et o\`u l'on note encore $x_{j}$ l'image de $x_{j} \in F(i_{j})$
par le morphisme $F(i_{j}) \rightarrow F(i_{n})$. Cette description montre que le foncteur
$\int_{I}$ pr\'eserve les \'equivalences. 

On remarque que le foncteur $Se_{I}$ est de Quillen \`a droite, comme
on peut le d\'eduire du fait que $\s -\mathbf{Cat}^{pr}$ est une
cat\'egorie de mod\`eles interne. Cette adjonction se d\'erive en une adjonction
$$\int_{I} : \mathrm{Ho}(\s -\mathbf{Cat}^{I}) \longleftrightarrow \mathrm{Ho}(\s -\mathbf{Cat}/I) : \mathbb{R}Se_{I}.$$
Soit maintenant $\pi : A \rightarrow I$ un objet de $\s -\mathbf{Cat}/I$. Nous dirons
qu'un morphisme $\alpha \in A(x,y)$ est \emph{cocart\'esien} si pour tout $z \in A$ 
le diagramme suivant
$$\xymatrix{
A(y,z) \ar[r]^-{\circ \alpha} \ar[d]  & A(x,z) \ar[d] \\
I(\pi(y),\pi(z)) \ar[r]_-{\circ \pi(\alpha)} & I(\pi(x),\pi(z))
}$$
est homotopiquement cart\'esien. Dans le cadre des cat\'egories, la notion
ci-dessus de morphismes cocart\'esiens est l\'eg\`erement plus
forte que la notion que l'on trouve
dans \cite[Exp. VI]{sga1}, et elle correspondant en r\'ealit\'e
\`a la condition $FibII'$ de \cite[Exp. VI, Prop. 6.11]{sga1}. Cependant, 
la notion ci-dessus de cat\'egories cofibr\'ees est compatible avec celle
de \cite[Exp. VI]{sga1}. 

\begin{df}\label{dcart}
\begin{enumerate}
\item Un objet $A\rightarrow I$ dans $\s -\mathbf{Cat}/I$ est \emph{cofibr\'e}
si pour tout morphisme $u : i\rightarrow j$ dans $I$, et 
tout objet $x$ de $A$ avec $\pi(x)=i$ il existe $y \in A$ et 
un morphisme cocart\'esien $\alpha \in A(x,y)$ avec $\pi(\alpha)$ isomorphe \`a $u$
dans la cat\'egorie comma $i/I$.
\item Un morphisme dans $\mathrm{Ho}(\s -\mathbf{Cat}/I)$ est \emph{cocart\'esien} s'il pr\'eserve les
morphismes cocart\'esiens.  \\
\end{enumerate}

La sous-cat\'egorie, non pleine, de $\mathrm{Ho}(s- \mathbf{Cat}/I)$ form\'ee des
objets cofibr\'es et des morphismes cocart\'esiens sera not\'ee
$\mathrm{Ho}(\s- \mathbf{Cat}/I)^{cocart}$.
\end{df}

Il existe une version duale de morphismes cart\'esiens, et de $\s$-cat\'egorie fibr\'ee. 
Nous laissons le soin au lecteur de dualiser ces notions ainsi que les \'enonc\'es
que nous allons maintenant d\'emontrer. \\

On remarque qu'un objet $\pi : A \rightarrow I$ est cofibr\'e si et seulement 
si la projection naturelle
$$p : \rh (\Delta^{1},A) \longrightarrow \rh (\{0\},A)\times_{\rh (\{0\},I)}^{h}\rh (\Delta^{1},I)$$
est essentiellement surjective lorsqu'elle est restreinte \`a la sous-$\s$-cat\'egorie pleine
$\rh (\Delta^{1},A)^{cart} \subset \rh (\Delta^{1},A)$, form\'ee
des morphismes cocart\'esiens. Ceci montre que la condition d'\^etre cofibr\'e est stable
par \'equivalences dans $\s -\mathbf{Cat}$. 

\begin{prop}\label{pcart}
Le foncteur $\int_{I} $
induit une \'equivalence de cat\'egorie
$$\int_{I} : \mathrm{Ho}(\s -\mathbf{Cat}^{I}) \longrightarrow \mathrm{Ho}(\s -\mathbf{Cat}/I)^{cart}.$$
\end{prop}

\textbf{Preuve --} Comme\c{c}ons par montrer que 
le foncteur $\int_{I}$ se factorise par la sous-cat\'egorie
$\mathrm{Ho}(\s -\mathbf{Cat}/I)^{cart}$. Pour cela, soit $F : I \longrightarrow \s -\mathbf{Cat}$ 
un foncteur, $\pi : \int_{I}F \rightarrow I$ son int\'egrale, et v\'erifions que
$\pi$ est cofibr\'e. Comme la condition d'\^etre cofibr\'e est visiblement invariante
par \'equivalences dans $\s -\mathbf{Cat}/I$ on peut supposer que $\pi$ est 
de plus une fibration de $\s$-cat\'egories. 
Soit $u : i \rightarrow j$ un morphisme dans $I$, et $a\in A$ 
avec $\pi(a)\simeq i$. Comme $\pi$ est une fibration l'isomorphisme
$\pi(a)\simeq i$ se rel\`eve un une \'equivalence $a \rightarrow a'$ 
dans $A$ avec $\pi(a')=i$. Quitte \`a travailler avec $a$' au lieu de $a$, 
cela montre que l'on peut supposer que $\pi(a)=i$. L'objet $a$ est donc 
de la forme $(i,x)$, avec $x\in F(i)$. Soit $y=u(x) \in F(j)$. Alors, le morphisme $u$ et
l'identit\'e de $y$ d\'efini un \'el\'ement $\alpha$ dans 
$$(\int_{I}F)((i,x),(j,y)):=
\coprod_{v : i \rightarrow j}F(j)(v(x),y).$$
On a $\pi(\alpha)=u$. De plus, si $(k,z)$ est un objet de $A$ le carr\'e suivant
$$\xymatrix{
\coprod_{v : j \rightarrow k}F(k)(v(y),z) \ar[r]^-{\circ \alpha} \ar[d] & 
\coprod_{w : i \rightarrow k}F(k)(w(x),z)
\ar[d] \\
I(j,k) \ar[r] & I(i,k)}$$
est \'evidemment cart\'esien, et donc homotopiquement cart\'esien car les 
bases sont des ensembles simpliciaux discrets. Ceci montre que $\alpha$ est cocart\'esien 
avec $\pi(\alpha)=u$, et donc que $\pi$ est cofibr\'e. 

Soit maintenant $f : F \longrightarrow G$ un morphisme dans $\s -\mathbf{Cat}^{I}$, et montrons que 
le morphisme induit
$$\int_{I}f : \int_{I}F \longrightarrow \int_{I}G$$
est un morphisme cocart\'esien. Pour cela, remarquons qu'un morphisme 
dans 
$$(\int_{I}F)((i,x),(j,y)):=
\coprod_{u : i \rightarrow j}F(j)(u(x),y),$$
de composante $(u,\gamma)$, avec $\gamma \in F(j)(u(x),y)$, est cocart\'esien
si et seulement si $\gamma$ est une \'equivalence dans $A$.  Or, le morphisme
$\int_{I}f$ envoie un morphisme$(u,\gamma)$ sur $(f(u),f(\gamma))$, ce qui montre bien 
que $\int_{I}f$ est toujours un morphisme cocart\'esien. 

Consid\'erons maintenant $\mathbb{R}Se_{I}^{cart}$ le sous-foncteur de 
$\mathbb{R}Se_{I}$ form\'e des sections cocart\'esiennes. En clair, pour un objet fibrant
$\pi : A \rightarrow I$  dans $\s -\mathbf{Cat}^{pr}/I$, et pour $i\in I$, $\mathbb{R}Se_{I}^{cart}(A)(i)$ 
est la sous-$\s$-cat\'egorie pleine de $Se_{I}(A)(i)=\underline{Hom}_{I}(i/I,A)$ 
form\'ee des morphismes cocart\'esiens $i/I \longrightarrow A$. On dispose d'une 
transformation naturelle $\mathbb{R}Se_{I}^{cart} \rightarrow \mathbb{R}Se_{I}$ ce qui permet 
de d\'efinir une transformation naturelle
$$h : \int_{I} \circ \mathbb{R}Se_{I}^{cart}  \longrightarrow id.$$
De m\^eme, la transformation naturelle
$$k : id \longrightarrow \mathbb{R}Se_{I} \int_{I}$$
se factorise par le sous-foncteur $\mathbb{R}Se^{cart}_{I} \int_{I}$, et induit ainsi une
transformation naturelle
$$k : id \longrightarrow \mathbb{R}Se^{cart}_{I} \int_{I}.$$
Les donn\'ees de $h$ et $k$ d\'efinissent de plus une adjonction
$$\int_{I} : \mathrm{Ho}(\s -\mathbf{Cat}^{I}) \longleftrightarrow \mathrm{Ho}(\s -\mathbf{Cat}/I)^{cart} : \mathbb{R}Se_{I}$$
Il nous reste donc \`a montrer que les deux assertions suivantes.

\begin{enumerate}
\item Pour tout objet cofibr\'e $\pi : A \longrightarrow I$ le morphisme d'adjonction
$$h : \int_{I}\circ \mathbb{R}Se_{I}^{cart}(A) \longrightarrow A$$
est un isomorphisme dans $\mathrm{Ho}(\s -\mathbf{Cat})$. 
\item Pour tout objet $F \in \mathrm{Ho}(\s -\mathbf{Cat}^{I})$, le morphisme d'adjonction
$$k : F \longrightarrow \mathbb{R}Se_{I}^{cart}(A)\int_{I}(F)$$
est un isomorphisme. 
\end{enumerate}

Le lemme suivant est le r\'esultat cl\'e.

\begin{lem}\label{lcart}
Soit $\pi : A \rightarrow I$ un objet cofibr\'e de $\mathrm{Ho}(\s -\mathbf{Cat}/I)$ et $i\in I$. Alors, le $\s$-foncteur
d'\'evaluation en $i$
$$ev_{i} : \rh _{I}(i/I,A) \longrightarrow \rh _{I}(\{i\},A)\simeq A\times^{h}_{I}\{i\}$$
est un isomorphisme dans $\mathrm{Ho}(\s -\mathbf{Cat})$.
\end{lem}

\textbf{Preuve --} On supposera que le morphisme $\pi : A \longrightarrow I$ est 
de plus une fibration. Il s'agit donc de montrer que l'\'evaluation en $i$ induit une 
\'equivalence
$$Se_{I}^{cart}(i/I,A) \simeq \pi^{-1}(i).$$
Soit $C$ la cat\'egorie des simplexes du nerf de $i/I$. Rappelons que 
ces objets sont les foncteurs $\Delta^{n} \longrightarrow i/I$, et les morphismes les diagrammes
(strictement) commutatifs
$$\xymatrix{\Delta^{n} \ar[rr] \ar[rd] & & \Delta^{m} \ar[dl] \\
 & i/I. & }$$
On dispose d'une projection naturelle $p : C \longrightarrow i/I$, qui \`a 
$f : \Delta^{n} \longrightarrow i/I$ associe $u(0)$ (o\`u les objets de
$\Delta^{n}$ sont not\'es $\xymatrix{0 \ar[r] & 1 \ar[r] & \dots \ar[r] & n}$). Cette
projection fait de $i/I$ la localisation, au sense de notre \S 1.2, de $C$ le long des morphismes verticaux (i.e. ceux
qui s'envoient sur des identit\'es dans $i/I$, voir par exemple \cite[Lem. 16.1]{hisi} ou encore
\ref{pa1}). En d'autres termes, le foncteur
$$p^{*} : \mathrm{Ho}(\mathbf{SEns}^{i/I}) \longrightarrow \mathrm{Ho}(\mathbf{SEns}^{C})$$
est pleinement fid\`ele et son image essentielle consiste en les diagrammes pour les quels
les morphismes verticaux de $C$ op\`erent par des \'equivalences. De plus, 
le morphisme naturel
$$\mathrm{\mathrm{Hocolim}}_{\Delta^{n} \rightarrow i/I}\Delta^{n} \longrightarrow i/I$$
est un isomorphisme dans $\mathrm{Ho}(\s -\mathbf{Cat})$ (voir \ref{pa2}). On d\'eduit de cela qu'il existe un isomorphisme
naturel dans $\mathrm{Ho}(\s -\mathbf{Cat})$
$$Se_{I}^{cart}(i/I,A)\simeq \mathrm{Holim}_{\Delta^{n} \rightarrow i/I}Se_{I}^{cart}(\Delta^{n},A).$$
De plus, par d\'efinition des morphismes cocart\'esiens et en utilisant le fait que $\pi : A \longrightarrow I$ 
soit cofibr\'e, on voit que le $\s$-foncteur d'\'evaluation en $\{0\} \in \Delta^{1}$
induit, pour tout $f : \Delta^{1} \longrightarrow I$,  une \'equivalence
$$Se_{I}^{cart}(\Delta^{1},A) \longrightarrow \pi^{-1}(f(0)).$$
De m\^eme, on voit que pour tout $f : \Delta^{n} \longrightarrow I$, le foncteur d'\'evaluation
en $\{n\}$ fournit une \'equivalence
$$Se_{I}^{cart}(\Delta^{n},A) \longrightarrow \pi^{-1}(f(0)).$$
Ceci montre que le foncteur
$$(\Delta^{n} \rightarrow i/I) \mapsto Se_{I}^{cart}(\Delta^{n},A)$$
est isomorphe, dans $\mathrm{Ho}(\s -\mathbf{Cat}^{C})$, \`a l'image r\'eciproque par $p : C \longrightarrow i/I$
d'un foncteur $E \in \mathrm{Ho}(\s -\mathbf{Cat}^{i/I})$.  
En particulier, on a 
$$\mathrm{Holim}_{\Delta^{n} \rightarrow i/I}Se_{I}^{cart}(\Delta^{n},A) \simeq 
\mathrm{Holim}_{i \rightarrow j}E(j)\simeq E(i)\simeq Se_{I}(\{i\},A)\simeq \pi^{-1}(i),$$
car $i$ est initial dans $i/I$.  Ceci termine la preuve du lemme. 
\hfill $\Box$ \\

Montrons comment le lemme pr\'ec\'edent implique que pour tout objet cofibr\'e $\pi : A \rightarrow I$
dans $\mathrm{Ho}(\s -\mathbf{Cat})$ le morphisme d'adjonction
$$\int_{I}\mathbb{R}Se_{I}^{cart}(A) \longrightarrow A$$
est un isomorphisme. 
Sans perte de g\'en\'eralit\'e nous pouvons supposer que $\pi$ est une
fibration. Soit alors $x$ et $y$ deux objets de $A$, d'images par $\pi$ not\'ees $i$ et $j$ respectivement.
Pour tout $u \in I(i,j)$ choisissons $\widehat{u} : x \rightarrow x_{u}$ un morphisme cocart\'esien
qui rel\`eve $u$ (comme $\pi$ est une fibration nous pouvons supposer de plus que 
se rel\`evement est strict: $\pi(\widehat{u})=u$). On dispose alors d'un 
morphisme (bien d\'efini dans $\mathrm{Ho}(\mathbf{SEns})$)
$$\coprod_{u\in I(i,j)}A(x_{u},y) \longrightarrow A(x,y),$$
qui consiste \`a pr\'ecomposer avec les morphismes $\widehat{u}$. Notons
$A^{id}(x_{u},y)$ le sous-ensemble simplicial de $A(x_{u},y)$ form\'e
des morphismes dont l'image par $\pi$ est \'egale \`a l'identit\'e de $j=\pi(x_{u})=\pi(y)$.
Alors,  comme les morphismes
$\widehat{u}$ ont \'et\'es choisis cocart\'esiens on voit que le morphisme induit
$$\phi : \coprod_{u\in I(i,j)}A^{id}(x_{u},y) \longrightarrow A(x,y),$$
est un isomorphisme dans $\mathrm{Ho}(\mathbf{SEns})$. Supposons alors que $x$ et $y$ soient les images
de deux objets, $x'$ et $y'$ dans $\int_{I}\mathbb{R}Se_{I}^{cart}(A)$.
Le lemme \ref{lcart} montre alors que le morphisme $\phi$ ci-dessus 
est isomorphe (dans $\mathrm{Ho}(\mathbf{SEns})$) au morphisme naturel
$$\left( \int_{I}\mathbb{R}Se_{I}^{cart}(A) \right)(x',y') \longrightarrow
A(x,y).$$
Ceci montre en particulier que $\int_{I}\mathbb{R}Se_{I}^{cart}(A) \longrightarrow A$ est 
pleinement fid\`ele. L'essentielle surjectivit\'e se voit imm\'ediatement
\`a l'aide du lemme \ref{lcart}, car il implique que 
toutes les fibres $\pi^{-1}(i) \subset A$ sont contenues dans 
l'image essentielle.

Le lemme \ref{lcart} implique d'autre part que le morphisme d'ajonction
$$F \longrightarrow \mathbb{R}Se_{I}^{cart}\int_{I}F$$
est un isomorphisme dans $\mathrm{Ho}(\s -\mathbf{Cat}^{I})$. En effet, il suffit de voir que 
le morphisme naturel
$$F(i) \longrightarrow (\int_{I}F)\times_{I}^{h}\{i\}$$
est une \'equivalence. Ce qui se d\'eduit ais\'ement du fait
que la fibre homotopique de $\int_{I}F \longrightarrow I$ en $i\in I$ 
est naturellement \'equivalente \`a la fibre na\"\i ve, qui elle est 
isomorphe \`a $F(i)$. 
\hfill $\Box$ \\

\subsection{$\s$-Cat\'egories localement pr\'esentables}

Nous utiliserons la notation suivante.

\begin{df}\label{d1}
\begin{enumerate}
\item La \emph{$\s$-cat\'egorie des types d'homotopie} est d\'efinie par
$$\T := L(\mathbf{SEns}).$$
Pour un univers fix\'e $\mathbb{U}$ nous noterons, lorsque cela est n\'ecessaire, 
$$\T _{\mathbb{U}} := L(\mathbf{SEns}_{\mathbb{U}}),$$
o\`u $\mathbf{SEns}_{\mathbb{U}}$ d\'esigne la cat\'egorie des ensembles simpliciaux $\mathbb{U}$-petits.
\item Pour une $\s$-cat\'egorie $A$ qui est $\mathbb{U}$-petite, la \emph{$\s$-cat\'egorie des pr\'efaisceaux sur $A$} est d\'efinie par
$$\widehat{A}:=\rh (A^{op},\T _{\mathbb{U}}).$$
\end{enumerate}
\end{df}

Remarquons que l'on a un isomorphisme naturel $\T \simeq Int(\mathbf{SEns})$ dans $\mathrm{Ho}(\s -\mathbf{Cat})$. De m\^eme, 
lorsque $A$ est une $\s$-cat\'egorie stricte ona 
$$\widehat{A}=\rh (A^{op},\T ) \simeq \rh (A^{op},Int(\mathbf{SEns})) \simeq Int(SPr(A))\simeq L(SPr(A)),$$
o\`u $SPr(A)$ est la cat\'egorie des pr\'efaisceaux simpliciaux sur $A$. 
Pour deux univers $\mathbb{U}\in \mathbb{V}$ le morphisme d'inclusion 
$$Int(\mathbf{SEns}_{\mathbb{U}}) \simeq \T _{\mathbb{U}} \longrightarrow 
Int(\mathbf{SEns}_{\mathbb{V}}) \simeq \T _{\mathbb{V}}$$ 
est ainsi pleinement fid\`ele. Pour une $\s$-cat\'egorie
$\mathbb{U}$-petite $A$ on trouve ainsi un $\s$-foncteur pleinement fid\`ele
$$\rh (A^{op},\T _{\mathbb{U}}) \longrightarrow \rh (A^{op},\T _{\mathbb{V}}).$$
Ainsi, bien que la d\'efinition pr\'ec\'edente de la $\s$-cat\'egorie $\widehat{A}$ d\'epende du 
choix d'un univers, les constructions que nous en tirerons ne d\'ependront pas de ce choix, quitte
\`a choisir un univers de d\'epart suffisemment gros. C'est pour cette raison que la notation
$\widehat{A}$ ne fait pas mention de l'univers $\mathbb{U}$. \\

Pour une $\s$-cat\'egorie stricte $A$ on dispose du plongement de Yoneda $h : A \longrightarrow
\widehat{A}$, adjoint du $\s$ foncteur $A \times A^{op} \longrightarrow \T$ qui envoie
$(a,b)$ sur $A(b,a)$. Pour une $\s$-cat\'egorie quelconque on d\'efinit 
le plongement de Yoneda $h : A \longrightarrow
\widehat{A}$ en choisissant un mod\`ele strict pour $A$ et en appliquand la construction 
pr\'ec\'edente.
Le $\s$-foncteur $h$ est toujours pleinement fid\`ele. Ainsi, 
pour deux $\s$-cat\'egories $A$ et $B$ le morphisme
$$\rh (A,B) \longrightarrow \rh (A,\widehat{B}) \simeq \widehat{A^{op}\times B}$$
est pleinement fid\`ele. Son image consiste en tous les $\s$-foncteurs
$F : A\times B^{op} \longrightarrow \T$ tel que pour tout $a\in A$ le 
$\s$-foncteur $F(a,-) \in \widehat{B}$ est repr\'esentable (i.e. isomorphe 
dans $\widehat{B}$ \`a un objet de l'image essentielle du plongement 
de Yoneda pour $B$). 

\begin{df}\label{d2}
Soit $\mathbb{U}$ un univers fix\'e.
Une $\s$-cat\'egorie $A$ est \emph{presque $\mathbb{U}$-localement pr\'esentable} s'il existe 
une $\s$-cat\'egorie $A_{0}$ qui est $\mathbb{U}$-petite et un $\s$-foncteur pleinement fid\`ele
$$A \longrightarrow \widehat{A_{0}}$$
qui poss\`ede un adjoint \`a gauche. 
\end{df}

La $\s$-cat\'egorie $\widehat{A_{0}}$ \'etant de la forme $Int(SPr(A_{0}'))$, pour $A_{0}'$ un 
mod\`ele strict de $A_{0}$,  elle poss\`ede 
des $\mathbb{U}$-limites et des $\mathbb{U}$-colimites. On d\'eduit formellement de cela
que toute $\s$-cat\'egorie presque $\mathbb{U}$-localement pr\'esentable poss\`ede aussi 
des $\mathbb{U}$-limites et des $\mathbb{U}$-colimites. Si $A \longrightarrow \widehat{A_{0}}$
est comme dans la d\'efinition pr\'ec\'edente, les limites dans $A$ sont simplement 
calcul\'ees dans $\widehat{A_{0}}$ (elles restent dans $A$ can le foncteur d'inclusion est 
adjoint \`a droite et donc commute aux limites). Les colimites sont 
aussi calcul\'ees dans $\widehat{A_{0}}$ puis ramen\'ees dans $A$ \`a l'aide de l'adjoint 
\`a gauche $\widehat{A_{0}} \longrightarrow A$. 

\begin{df}\label{d3}
Nous dirons qu'une $\s$-cat\'egorie $A$ est \emph{$\mathbb{U}$-localement pr\'esentable}
s'il existe une cat\'egorie de mod\`eles $M$ qui est $\mathbb{U}$-combinatoire
et un isomorphisme dans $\mathrm{Ho}(\s -\mathbf{Cat})$
$$A \simeq LM.$$
\end{df}

Une remarque importante est que $\uscat_{\mathbb{U}}$ est 
une $\s$-cat\'egorie $\mathbb{U}$-localement pr\'esentable. \\

D'apr\`es \cite{du} toute cat\'egorie de mod\`eles $\mathbb{U}$-combinatoire est Quillen \'equivalente
\`a une localisation de Bousfield d'une cat\'egorie de pr\'efaisceaux simpliciaux sur
une cat\'egorie $\mathbb{U}$-petite. Ceci implique qu'une $\s$-cat\'egorie 
$\mathbb{U}$-localement pr\'esentable est presque $\mathbb{U}$-localement pr\'esentable. R\'eciproquement
on peut montrer qu'une $\s$-cat\'egorie presque $\mathbb{U}$-localement pr\'esentable est 
$\mathbb{U}$-localement pr\'esentable si et seulement si tous ses objets
sont $\kappa$-compact (pour un cardinal r\'egulier $\kappa$ qui d\'epend de l'objet consid\'er\'e).

Une des propri\'et\'es fondamentales des $\s$-cat\'egories (presque) $\mathbb{U}$-localement pr\'esentables
que nous utiliserons est le crit\`ere d'existence d'adjoint suivant. 

\begin{prop}\label{p1}
Soit $f : A \longrightarrow B$ un $\s$-foncteur.
\begin{enumerate}
\item Si $A$ est presque $\mathbb{U}$-localement pr\'esentable et si 
$f$ commute aux $\mathbb{U}$-colimites alors $f$ poss\`ede un adjoint \`a droite.

\item Si $A$ est $\mathbb{U}$-compactement engendr\'ee et si 
$F$ commute aux $\mathbb{U}$-limites et aux colimites $\kappa$-filtrantes pour un cardinal
r\'egulier $\kappa$, alors $f$ poss\`ede un adjoint \`a gauche.
\end{enumerate}
\end{prop}

\textbf{Preuve --} Quitte \`a choisir des mod\`eles stricts nous supposerons que 
les $\s$-cat\'egories en jeu sont des $\mathbb{S}$-cat\'egories. 

$(1)$ On consid\`ere l'adjonction de Quillen
$$f_{!} : \mathbf{SEns}_{\mathbb{V}}^{A} \longleftrightarrow \mathbf{SEns}_{\mathbb{V}}^{B} : f^{*}$$
o\`u $\mathbb{U}\in \mathbb{V}$ avec $A$ et $B$ \'el\'ements de l'univers $\mathbb{V}$. 
Cette adjonction de Quillen donne lieu \`a une adjonction de $\s$-cat\'egories
$$\mathbb{L}f_{!} : \widehat{A} \simeq L(\mathbf{SEns}_{\mathbb{V}}^{A}) \longleftrightarrow L(\mathbf{SEns}_{\mathbb{V}}^{B})\simeq \widehat{B} : \mathbb{R}f^{*}.$$
On dispose de plus d'un diagramme commutatif dans $\mathrm{Ho}(\s -\mathbf{Cat})$
$$\xymatrix{
A \ar[r]^-{f} \ar[d] & B \ar[d] \\
\widehat{A} \ar[r]_-{\mathbb{L}f_{!}} & \widehat{B},}$$
o\`u les morphismes verticaux sont les plongements de Yoneda. On voit ainsi que pour montrer que 
$f$ poss\`ede un 
adjoint \`a droite il suffit de montrer que le $\s$-foncteur $\mathbb{R}f^{*}$ 
pr\'eserve les objets repr\'esentables. Pour tout $b\in B$ le $\s$-foncteur
$f^{*}(h_{b}) : A \longrightarrow \mathbb{T}^{op}$ commute aux $\mathbb{U}$-colimites par hypoth\`eses.
On voit ainsi qu'il suffit de montrer que tout objet $F \in \widehat{A}$ qui envoie colimites dans
$A$ sur des limites dans $\mathbb{T}$ est repr\'esentables. Pour cela, on \'ecrite $A$ comme
une sous-$\s$-cat\'egorie reflexive d'une $\s$-cat\'egorie de pr\'efaisceaux
$$i : A \longleftrightarrow \widehat{A_{0}} : l.$$
Le foncteur $l\circ F : \widehat{A_{0}}^{op} \longrightarrow \mathbb{T}$ 
commute aux $\mathbb{U}$-limites et est donc repr\'esentable: l'objet de $\widehat{A_{0}}$
qui le repr\'esentable n'est autre que sa propre restriction comme $\s$-foncteur $A_{0}^{op} \longrightarrow \mathbb{T}$. En effet, cela se d\'eduit du fait que l'adjonction de Quillen 
induite par le plongement de Yoneda $h : A_{0} \longrightarrow \widehat{A_{0}}$
$$h_{!} : \mathbf{SEns}_{\mathbb{V}}^{A_{0}^{op}} \longleftrightarrow \mathbf{SEns}_{\mathbb{V}}^{\widehat{A_{0}}^{op}} : h^{*}$$
induise un foncteur pleinement fid\`ele 
$$\mathbb{L}h_{!} : \mathrm{Ho}(\mathbf{SEns}_{\mathbb{V}}^{{A_{0}}^{op}}) \longrightarrow
\mathrm{Ho}(\mathbf{SEns}_{\mathbb{V}}^{\widehat{A_{0}}^{op}})$$
dont l'image essentielle consiste en les pr\'efaisceaux simpliciaux sur $\widehat{A_{0}}$ qui 
envoient colimites sur limites (voir \cite[Lem. 7.3]{to1}). 

$(2)$ Vue les compatibilit\'es entre limites et colimites homotopiques dans une cat\'egorie
de mod\`eles $M$ et les limites et colimites dans la $\s$-cat\'egorie $LM$ on voit qu'il 
faut montrer l'\'enonc\'e suivant: si $F : M \longrightarrow \mathbf{SEns}_{\mathbb{V}}$ est un foncteur
qui v\'erifie les conditions suivantes

\begin{enumerate}
\item $F$ pr\'eserve les \'equivalences.

\item $F$ commute aux limites homotopiques.

\item $F$ commute aux colimites homotopiques $\kappa$-filtrantes.

\end{enumerate}

Alors il existe un objet $x\in M$ tel que $F$ et $Map(x,-)$ soient isomorphes
comme objets dans $SPr(M^{op})$ (ici $M$ est $\mathbb{V}$-petite et les pr\'efaisceaux
sont des $\mathbb{V}$-pr\'efaisceaux). 

D'apr\`es les r\'esultats de \cite{du} on peut supposer que $M$ est de plus simpliciale et 
avec tous les objets cofibrants. Les Homs simpliciaux seront alors
not\'es $\underline{Hom}$.
Quitte \`a choisir $\kappa$ suffisemment grand on peut faire en sorte 
que $\kappa$ soit un cardinal suffisament grand (comme
pour les preuves des r\'esultats de  \cite{du}). Notons alors
$C\subset M$ la sous-cat\'egorie pleine form\'ee des objets $\kappa$-compacts. \\
Notons $\underline{h}^{x} : M \longrightarrow \mathbf{SEns}$ le foncteur qui 
\`a $y$ associe $Hom(x,R(y))$, o\`u $R$ est un 
foncteur de remplacement fibrant. De m\^eme, notons 
$h^{x} : M \longrightarrow Ens$ le foncteur qui \`a $y$ associe $Hom(x,y)$, que nous
verrons comme un foncteur en ensembles simpliciaux discrets. On dispose d'un 
morphisme naturel $h^{x} \longrightarrow \underline{h}^{x}$. 

On note 
$\Delta(C/F)$ la cat\'egorie des simplexes de $F$ au-dessus de $C$. Ses objets sont 
des triplets $(x,n,a)$, avec $x\in C$, $[n] \in \Delta$ et $a\in F(x)_{n}$. 
Les morphismes $(x,n,a) \rightarrow (y,m,b)$ sont donn\'es par
un morphisme $x \rightarrow y$ dans $C$ et un morphisme
$[m] \rightarrow [n]$ dans $\Delta$ qui envoient le simplexe $a \in F(x)_{n}$ sur $b \in F(y)_{m}$.
 On dispose d'un 
foncteur naturel $\phi : \Delta(C/F) \longrightarrow C$ qui \`a $(x,n,a)$
associe $x\in C$. Notons $X:=holim_{\Delta(C/F)}\phi \in M$ la limite
homotopique de ce foncteur. On dispose d'un diagramme naturel dans $\mathrm{Ho}(SPr(M^{op}))$
$$\xymatrix{
\mathrm{\mathrm{Hocolim}}_{(x,n,a) \in \Delta(C/F)}h^{x} \ar[r] \ar[d]  & F \\
h^{X}. & & }$$
De plus, comme $F$ commute aux limites homotopiques, il existe un unique
morphisme $h^{X} \longrightarrow F$ qui rende le diagramme
$$\xymatrix{
\mathrm{\mathrm{Hocolim}}_{(x,n,a)\in \Delta(C/F)}h^{x} \ar[r] \ar[d]  & F \\
h^{X} \ar[ru] &  }$$
commutatif dans $\mathrm{Ho}(SPr(M^{op}))$. Comme $F$ pr\'eserve les \'equivalences, ce diagramme
induit un nouveau diagramme commutatif (voir \cite[Lem. 4.2.2]{hagI})
$$\xymatrix{
\mathrm{Hocolim}_{(x,n,a) \in \Delta(C/F)}\underline{h}^{x} \ar[r] \ar[d]  & F \\
\underline{h}^{X}. \ar[ru] &  }$$
Il est clair que le morphisme horizontal devient une \'equivalence lorsque l'on 
restreint les foncteurs \`a la sous-cat\'egorie $C$ (voir \cite[Prop. 4.7]{du}). De plus, comme
les deux foncteurs en question commutent aux colimites $\kappa$-filtrantes, et que
de plus tout objet de $M$ est \'equivalent \`a une colimite homotopique $\kappa$-filtrante
d'objets de $C$, ceci implique que le morphisme
$$\mathrm{Hocolim}_{(x,n,a) \in \Delta(C/F)}\underline{h}^{x} \longrightarrow F$$
est une \'equivalence. Ainsi, le morphisme 
$\underline{h}^{X} \longrightarrow F$ poss\`ede une section dans $\mathrm{Ho}(SPr(M^{op}))$. Le foncteur 
$F$ est donc un r\'etracte d'un repr\'esentable et donc est lui-m\^eme repr\'esentable
par un facteur direct (dans $\mathrm{Ho}(M)$) de $X$. \hfill $\Box$ \\

Notons au passage un r\'esultat que nous avons utilis\'e lors de la preuve du la
partie $(1)$ de la proposition pr\'ec\'edente: pour une $\s$-cat\'egorie $A_{0}$ qui est 
$\mathbb{U}$-petite et pour $B$ toute $\s$-cat\'egorie qui poss\`ede des
$\mathbb{U}$-colimites, le morphisme de restriction
$$\rh (\widehat{A_{0}},B) \longrightarrow \rh (A_{0},B)$$
induit un isomorphisme dans $\mathrm{Ho}(\s -\mathbf{Cat})$
$$\rh^{c} (\widehat{A_{0}},B) \longrightarrow \rh (A_{0},B),$$
o\`u $\rh^{c}(\widehat{A_{0}},B)$ d\'esigne la sous-$\s$-cat\'egorie pleine
de $\rh (\widehat{A_{0}},B)$ form\'ee des $\s$-foncteurs qui commutent aux
colimites.

\subsection{$\s$-Topos}

\begin{df}\label{d4}
\begin{enumerate}
\item
Un \emph{$\mathbb{U}-\s$-topos} est une $\s$-cat\'egorie $T$ tel qu'il existe une $\s$-cat\'egorie 
$\mathbb{U}$-petite $T_{0}$ et un $\s$-foncteur pleinement fid\`ele
$$i : T \longrightarrow \widehat{T_{0}}$$
poss\`edant un adjoint \`a gauche qui commute aux limites finies (c'est \`a dire aux
produits fibr\'es et aux produits finis). 
\item
Un \emph{$\s$-topos} est un $\mathbb{U}-\s$-topos pour un univers $\mathbb{U}$.
\item Si $T$ et $T'$ sont deux $\s$-topos, un \emph{morphisme g\'eom\'etrique} $f : T \longrightarrow T'$, est 
un $\s$-foncteur qui commute aux colimites et dont l'adjoint \`a gauche commute aux limites finies. 
\end{enumerate}
La \emph{$\s$-cat\'egorie des $\mathbb{U}-\s$-topos} est la sous-$\s$-cat\'egorie (non pleine) 
de $\uscat$ form\'ee des $\mathbb{U}-\s$-topos et des morphismes g\'eom\'etriques. Elle
sera not\'ee $\ustop^{\mathbb{U}}$ (ou $\ustop$ si l'on ne souhaite pas mentioner univers $\mathbb{U}$).
\end{df}

Par d\'efinition les $\mathbb{U}-\s$-topos sont des $\s$-cat\'egories presques $\mathbb{U}$-localement
pr\'esentables. Nous verrons que sous une hypoth\`ese additionnelle de compl\`etude 
ils sont en fait $\mathbb{U}$-localement pr\'esentables. \\

Tout comme le cas des topos les $\s$-topos sont caract\'erisables par 
des axiomes de type Giraud (voir \cite{lu1,hagI,tove2}). Pour un $\s$-topos $T$ fix\'e 
nous retiendrons les propri\'et\'es caract\'eristiques suivantes.

\begin{enumerate}
\item Il existe une sous-$\s$-cat\'egorie $T_{0}\subset T$ qui est $\mathbb{U}$-petite et telle que
le $\s$-foncteur induit 
$$T \longrightarrow \widehat{T_{0}}$$
soit pleinement fid\`ele (on dit aussi que \emph{$T_{0}$ engendre $T$ par \'epimorphismes
stricts}).
\item La $\s$-cat\'egorie $T$ poss\`ede des $\mathbb{U}$-limites et des $\mathbb{U}$-colimites.
\item Les $\mathbb{U}$-colimites dans $T$ sont universelles: si 
$x : I \longrightarrow T$ est un diagramme d'objets de $T$, augment\'e au-dessus
d'un objet $z\in T$, alors pour tout $y \rightarrow z$ dans $T$, le morphisme naturel
$$colim_{I}(y\times_{z}x_{i}) \rightarrow y\times_{z}(colim_{I}x_{i})$$
est un isomorphisme dans $[T]$. 
\item Les sommes sont disjointes dans $T$: si $\{x_{i}\}_{i\in I}$ est une famille d'objets
dans $T$, de somme $\coprod_{i\in I}x_{i}=:x$, alors pour tout $i\neq j$ les morphimes 
$$x_{i} \longrightarrow x_{i}\times_{x}x_{i} \qquad \emptyset \longrightarrow x_{i} \times_{x}x_{j}$$
sont des isomorphismes dans $[T]$. 
\item Les relations d'\'equivalences sont effectives: si $X_{*} : \Delta^{op} \longrightarrow T$
est un \emph{groupo\"\i de de Segal dans $T$} (voir \cite{hagI,tove2}), alors le morphisme naturel
$$X_{1} \longrightarrow X_{0}\times_{colim_{\Delta} X_{*}}X_{0}$$
est un isomorphisme dans $[T]$. 
\end{enumerate}

Notons aussi que les propri\'et\'es $(3)-(5)$ peuvent aussi s'exprimer de mani\`ere uniforme par 
la \emph{propri\'et\'e de conservation}. Cette propri\'et\'e apparait dans \cite{rez} et peut 
se traduire de la fa\c{c}on suivante. Pour l'\'enoncer nous aurons besoin 
de la notion de $\s$-cat\'egorie comma: si $T$ est une $\s$-cat\'egorie
et $x\in T$ est un objet nous noterons $T/x$ la $\s$-cat\'egorie d\'efinie comme la fibre homotopique
en $x$ du morphisme
$$\rh (\Delta^{1},T) \longrightarrow T,$$
o\`u $\Delta^{1}$ est la cat\'egorie classifiant les morphismes ($\Delta^{1}:=\left( 0 \rightarrow 1 \right)$), 
et le morphisme ci-dessus est l'\'evaluation 
en l'objet $1 \in \Delta^{1}$. Soit maintenant $x : I \longrightarrow T$ un diagramme d'objet 
de $T$ avec $I$ une cat\'egorie $\mathbb{U}$-petite. On consid\`ere
$\rh (I,T)$ la $\s$-cat\'egorie des $I$-diagrammes dans $T$ et $\rh (I,T)/x$ la $\s$-cat\'egorie
des diagrammes au-dessus de $x$. Lorsque $T$ poss\`ede des colimites on dispose du $\s$-foncteur
$$colim_{I} : \rh (I,T)/x \longrightarrow T/|x|,$$
o\`u $|x|$ est la colimite du diagramme $x$. Lorsque $T$ poss\`ede des limites finies ce foncteur
poss\`ede un adjoint \`a droite
$$-\times_{|x|}x : T/|x| \longrightarrow \rh (I,T)/x.$$
La propri\'et\'e de conservation affirme alors que le foncteur 
$-\times_{|x|}x$ est pleinement fid\`ele et que son image essentiellement consiste
en des diagrammes $y : I \longrightarrow T$, augment\'es sur $x$, et tel que pour tout
$i\rightarrow j$, morphisme dans $I$, le morphisme
$$y_{i} \longrightarrow (x_{i})\times_{(x_{j})}(y_{j})$$
est un isomorphisme dans $[T]$.  En d'autres termes, le foncteur
$$colim_{I} : \rh (I,T)/x \longrightarrow T/|x|,$$
restreint \`a la sous-cat\'egorie des diagrammes $y \rightarrow x$ v\'erifiant la propri\'et\'e ci-dessus, 
est une \'equivalence. \\

Soit $T$ un $\s$-topos. Un objet $x\in T$ est \emph{$n$-tronqu\'e} si pour tout 
$y\in T$ l'ensemble simplicial $T(y,x)$ est $n$-tronqu\'e. De mani\`ere \'equivalente, 
$x$ est $n$-tronqu\'e si pour tout $i>n$ le morphisme naturel
$x \longrightarrow x^{S^{i}}$ est un isomorphisme dans $[T]$. La sous-$\s$-cat\'egorie pleine
des objets $n$-tronqu\'es de $T$ sera not\'ee $T_{\leq n} \subset T$. Les $\s$-foncteurs
d'inclusion $T_{\leq n} \longrightarrow T$ poss\`ede des adjoints \`a gauche
$t_{n} : T \longrightarrow T_{\leq n}$. Nous dirons alors que $T$ est \emph{t-complet}
si et seulement si tout objet $x\in T$ tel que $t_{n}(x) \simeq *$ pour tout $n$ est tel que
$x\simeq *$. En d'aurtes termes, $T$ est t-complet s'il n'existe aucun objet non trivial qui soit
$n$-connexe pour tout $n$. Les $\s$-topos t-complets se d\'ecrivent alors \`a l'aide
d'une notion de $\s$-topologie. On montre plus pr\'ecis\'ement que 
pour une $\s$-cat\'egorie $T_{0}$ fix\'ee, il existe une bijection entre 
les classes d'\'equivalences de sous-$\s$-topos t-complets de $\widehat{T_{0}}$ et les
topologie de Grothendieck sur la cat\'egorie $[T_{0}]$ (voir \cite[Thm. 3.8.3]{hagI}). Ceci implique en particulier
un proc\'ed\'e de construction des $\s$-topos \`a l'aide d'une notion de 
\emph{$\s$-site}, c'est \`a dire \`a l'aide de couple $(T_{0},\tau)$, form\'es
d'une $\s$-cat\'egorie $T_{0}$ est d'une topologie sur $[T_{0}]$. En effet, pour un tel 
$\s$-site, avec $T_{0}$ strict, on peut construire une cat\'egorie de mod\`eles $SPr_{\tau}(T_{0})$ qui est 
une localisation de Bousfield \`a gauche pour une notion d'\'equivalences
locale dans $SPr(T_{0})$ (voir \cite[\S 3.4]{hagI}). La $\s$-cat\'egorie $L(SPr_{\tau}(T_{0}))$
est alors un $\s$-topos t-complet, que l'on notera $T_{0}^{\sim,\tau}$. Une cons\'equence
de cela est que tout $\s$-topos t-complet est de la forme $LM$ pour $M$ un 
topos de mod\`eles t-complet. En particulier tout $\s$-topos
t-complet est $\mathbb{U}$-localement pr\'esentable.\\

\noindent \textbf{Champs \`a valeurs dans une $\s$-cat\'egorie.} Pour finir, soit $(T_{0},\tau)$ un $\s$-site qui est de plus $\mathbb{U}$-petit et consid\'erons
le $\s$-topos $T:=T_{0}^{\sim,\tau}$ associ\'e. C'est la sous-$\s$-cat\'egorie
pleine de $\widehat{T_{0}}$ form\'ee des $\s$-foncteurs $T_{0}^{op} \longrightarrow \mathbb{T}_{\mathbb{U}}$
qui poss\`ede la propri\'et\'e de descente pour les hyper-recouvrements (voir \cite[Cor. 3.4.7]{hagI}). 
Le $\s$-foncteur d'inclusion $i : T \hookrightarrow \widehat{T_{0}}$
poss\`ede un adjoint \`a gauche exact
$$a : \widehat{T_{0}} \longrightarrow T$$
(c'est le $\s$-foncteur de \emph{faisceautisation} aussi appel\'e \emph{champ associ\'e}). 
Le compos\'e du plongement de Yoneda
$h: T_{0} \longrightarrow \widehat{T_{0}}$ et de $a$ d\'efinit un $\s$-foncteur
$$h^{\sim} : T_{0} \longrightarrow T$$
qui n'est plus pleinement fid\`ele en g\'en\'eral (par d\'efinition il l'est si 
la topologie $\tau$ est sous-canonique). Soit maintenant $A$ une $\s$-cat\'egorie qui poss\`ede 
des $\mathbb{U}$-limites. On dispose alors d'un morphisme de restriction
$$\rh ^{l}(T^{op},A) \longrightarrow \rh (T_{0},A),$$
o\`u $\rh ^{l}(T^{op},A)$ est la sous-$\s$-cat\'egorie pleine de 
$\rh (T^{op},A)$ form\'ee des $\s$-foncteurs qui commutents aux limites. On v\'erifie que ce $\s$-foncteur
est pleinement fid\`ele. Son image dans $\rh (T_{0},A)$ est appel\'e par d\'efinition
la $\s$-cat\'egorie des \emph{champs sur $T_{0}$ \`a valeurs dans $A$}. Cette image se caract\'erise de
la fa\c{c}on suivante: un $\s$-foncteur $F : T_{0}^{op} \longrightarrow A$ est
un champs s'il poss\`ede la propri\'et\'e de descente par rapport aux
hyper-recouvrements. Nous posons la d\'efinition suivante.

\begin{df}\label{d5}
Pour un $\mathbb{U}-\s$-topos $T$ et une $\s$-cat\'egorie $A$ poss\'edant des $\mathbb{U}$-limites, 
la \emph{$\s$-cat\'egorie des champs sur $T$ \`a valeurs dans $A$} (aussi appel\'ee la 
\emph{$\s$-cat\'egorie des $A$-champs sur  $T$}) est $\rh ^{l}(T^{op},A)$. Elle sera not\'e
$Ch(T,A)$. 
\end{df}

Notons pour finir que le foncteur d'inclusion 
$$\rh ^{l}(T^{op},A) \longrightarrow \rh (T_{0}^{op},A),$$
poss\`ede lui aussi un adjoint \`a gauche exact. Cet adjoint \`a gauche est 
construit en composant l'\'equivalence 
$$\rh (T_ {0}^{op},A)^{op}\simeq \rh (T_{0},A^{op})\simeq \rh^{c}(\widehat{T_{0}},A^{op})$$
avec le $\s$-foncteur 
$$a_{!} : \rh^{c}(\widehat{T_{0}},A^{op}) \longrightarrow \rh^{c}(T,A^{op})\simeq \rh^{l}(T^{op},A)^{op},$$
adjoint \`a gauche de la restriction suivant $a$.
o\`u $a : \widehat{T_{0}} \longrightarrow T$ est le $\s$-foncteur champ associ\'e.

\section{Structures monoidales}

Dans cette seconde section nous pr\'esentons la th\'eorie des
$\s$-cat\'egories mono\"\i dales sym\'etriques. Nous adopterons 
le point de vue qui \'etait utilis\'e par le premier auteur dans
\cite{to0}, et de nombreuses notions et constructions sont tir\'es de ce texte.
L'\'equivalence entre la th\'eorie des $\s$-cat\'egories cofibr\'ees sur
$\Gamma$ et celle des $\s$-foncteurs de $\Gamma \longrightarrow \s -\mathbf{Cat}$
implique que notre approche est aussi \'equivalente \`a celle 
utilis\'ee r\'ecemment par J. Lurie (voir \cite{lu2} par exemple).

En dehors des aspects d\'efinitionels cette section contient essentiellement
deux r\'esultats fondamentaux. Tout d'abord le Th\'eor\`eme \ref{t1} qui affirme
l'existence de $\s$-cat\'egories mono\"\i dales sym\'etriques rigides 
librement engendr\'ee par des $\s$-cat\'egories. Une cons\'equence
de ce fait est la Proposition \ref{p6} affirmant l'existence
et l'unicit\'e des morphismes de traces dans les $\s$-cat\'egories mono\"\i dales sym\'etriques rigides. Le second r\'esultat (Cor. \ref{c4} et Def. \ref{d13}) affirme de plus que cette trace
poss\`ede une propri\'et\'e d'invariance cyclique qui sera pour nous
cruciale pour la construction du caract\`ere de Chern. Nous n'avons pas
trouv\'e de preuve directe raisonnable \`a ce second
r\'esultat, et nous montrons qu'il s'agit d'une cons\'equence 
de la solution \`a l'\emph{hypoth\`ese du cobordisme} (\cite{lu2}). Ce second r\'esultat est de toute \'evidence le th\'eor\`eme
le plus important de cette section, voire de cet article. Nous terminerons
cette partie en montrant que, comme il se doit, la trace
est multiplicative par rapport au la structure mono\"\i dale.

\subsection{Th\'eorie homotopique des $\s$-cat\'egories mono\"\i dales sym\'etriques}

Nous consid\'erons la cat\'egorie des ensembles finis point\'es et applications
point\'ees. Par soucis de petitesse nous nous restreindrons \`a un squelette de 
de cette cat\'egorie, not\'e $\Gamma$, form\'e des ensembles $\{0, \dots, n\}$, point\'es
en $0$, pour tout $n\geq 0$. L'objet $\{0, \dots, n\}$ sera not\'e
$[n]$. 

Pour un tel objet $[n] \in \Gamma$, nous disposons de $n$-morphismes canoniques
$$s_{i} : [n] \longrightarrow [1],$$
pour $0<i\leq n$, 
d\'efinis par $s_{i}(j)=\delta_{ij}$ (i.e. $s_{i}$ envoie tout sur $0$ sauf $i$). Les morphismes
$s_{i}$ sont g\'en\'eralement appel\'es les \emph{morphismes de Segal}.

\begin{df}\label{d6}
\begin{enumerate}
\item Une \emph{$\s$-cat\'egorie pr\'emono\"\i dale symm\'etrique} (nous dirons aussi \emph{$\s$-CPS}) 
est la donn\'ee d'un foncteur
$$A : \Gamma \longrightarrow \s -\mathbf{Cat}.$$
\item Une $\s$-CPS $A$ est une \emph{$\s$-cat\'egorie mono\"\i dale symm\'etrique} (nous dirons 
aussi \emph{$\s$-CMS}) si 
pour tout entier $n\geq 0$ le morphisme
$$\prod_{1\leq i\leq n}s_{i} : A([n]) \longrightarrow A([1])^{n}$$
est une \'equivalence de $\s$-cat\'egories. 
\end{enumerate}
Un morphisme de $\s$-CPS (resp. de $\s$-CMS) est une transformation naturelle
de foncteurs $\Gamma \longrightarrow \s -\mathbf{Cat}$.
\end{df}

Notons que la condition d'\^etre une $\s$-CMS implique que 
$A([0])$ est \'equivalente \`a la $\s$-cat\'egorie ponctuelle $*$ (c'est le cas
$n=0$ dans la condition de la d\'efinition). \\

La cat\'egorie des $\s$-CPS est par d\'efinition la sous-cat\'egorie pleine
de la cat\'egorie des foncteurs $\s- \mathbf{Cat}^{\Gamma}$ form\'ee des $\s$-CPS. De m\^eme, 
la cat\'egorie des $\s$-CMS est la sous-cat\'egorie pleine de la cat\'egorie
des $\s$-CPS form\'ee des $\s$-CMS. Ces deux cat\'egories seront not\'ees
respectivement 
$$\scps \qquad \scms$$ 
(ou encore $\scps_{\mathbb{U}}$,  $\scms_{\mathbb{U}}$
si fixer des univers il nous faut). Les morphismes de $\s$-CPS et de $\s$-CMS sont par
d\'efinitions les transformations naturelles de $\Gamma$-diagrammes. Les morphismes
entre $\s$-CMS seront
appel\'es des \emph{$\otimes-\s$-foncteurs}.

On dispose sur $\scps$ de la structure de mod\`eles projective pour laquelle 
les \'equivalences (resp. les fibrations) sont les morphismes $A \longrightarrow B$ tels que
$A([n]) \longrightarrow B([n])$ soit une \'equivalence (resp. une fibration) pour tout $n$.
Nous noterons $\mathrm{Ho}(\scps)$ la cat\'egorie homotopique correspondante, et 
$\mathrm{Ho}(\scms) \subset \mathrm{Ho}(\scps)$ la sous-cat\'egorie pleine des $\s$-CMS. Par extension
les morphismes de $\mathrm{Ho}(\scms)$ seront aussi appel\'es des $\otimes-\s$-foncteurs.
Nous verrons par la
suite que $\mathrm{Ho}(\scms)$ s'identifie naturellement \`a la cat\'egorie homotopique d'une
structure de mod\`eles localis\'ee sur $\scps$ dont les objets locaux sont 
exactement les $\s$-CMS. Mais avant cela nous pr\'esentons quelques exemples.\\

\noindent \textbf{-- Premiers exemples de $\s$-CMS --} 
\begin{enumerate}

\item \textbf{Mono\"\i des ab\'eliens.} La sous-cat\'egorie pleine de $\scms$ form\'ee 
des foncteurs $A : \Gamma \longrightarrow \s -\mathbf{Cat}$ tels que 
$A([n])$ soit un ensemble (vu comme $\s$-cat\'egorie discr\`ete) pour tout $n$ est \'equivalente \`a la cat\'egorie
des mono\"\i des commutatifs (associatifs et unitaires). Pour voir cela on construit 
un foncteur $\Gamma^{op} \longrightarrow CMon$, de $\Gamma^{op}$
vers la cat\'egorie $CMon$ des mono\"\i des commutatifs de la fa\c{c}on suivante.
L'objet $[n]$ est envoy\'e sur $Hom^{0}([n],\mathbb{N})$, l'ensemble des applications
$f : [n] \longrightarrow \mathbb{N}$ telles que $f(0)=0$, muni de sa structure de mono\"\i de
induite par la structure additive de $\mathbb{N}$. Ceci d\'efinit un foncteur 
de mani\`ere \'evidente $\Gamma^{op} \longrightarrow CMon$, en envoyant un morphisme $[n] \rightarrow [m]$
sur le morphisme $Hom^{0}([m],\mathbb{N}) \longrightarrow Hom^{0}([n],\mathbb{N})$ obtenu par 
composition. Pour un mono\"\i de commutatif $E \in CMon$ on d\'efinit alors un foncteur
$N_{\Gamma}(E) : \Gamma \longrightarrow Ens$ en pr\'ecomposant le foncteur $Hom(-,E)$ avec
$\Gamma^{op} \longrightarrow CMon$. En d'autres termes on a
$N_{\Gamma}(E)([n])=E^{n}$, et pour $u : [n] \rightarrow [m]$ dans $\Gamma$ le morphisme induit
$$u_{!} : E^{n} \longrightarrow E^{m}$$
est donn\'e pour $x=(x_{i}) \in E^{n}$ par la formule
$$u_{!}(x)_{j}:=\Sigma_{i\in u^{-1}(j)}x_{i}.$$
On v\'erifie alors que le foncteur $E \mapsto N_{\Gamma}(E)$ est une \'equivalence
de $CMon$ vers la sous-cat\'egorie de $\s$-CMS form\'ee des foncteurs qui prennent leurs valeurs
dans les ensembles. 

\item \textbf{Cat\'egories mono\"\i dales symm\'etriques.} (voir aussi \cite[Def. 3.3.7]{lei}). 
L'exemple pr\'ec\'edent se g\'en\'eralise au cas des cat\'egories mono\"\i dales 
symm\'etriques de la fa\c{c}on suivante. Notons $\mathbf{Cat}^{\otimes}$ la cat\'egorie
des cat\'egories mono\"\i dales symm\'etriques et des foncteurs mono\"\i daux
symm\'etriques. Notons $Fin$ la cat\'egorie dont les objets sont les
ensembles $\{1, \dots, n\}$ et dont les morphimes sont les applications bijectives. La cat\'egorie
$Fin$ est munie de sa structure mono\"\i dale symm\'etrique usuelle induite par 
la somme disjointe d'ensembles (en identifiant $\{1, \dots, n\} \coprod \{1, \dots, m\}$ avec 
$\{1, \dots, n+m\}$ par concat\'enation).
On construit un foncteur $\Gamma^{op} \longrightarrow \mathbf{Cat}^{\otimes}$ 
en envoyant $[n]$ sur $Hom^{0}([n],Fin)$, la cat\'egorie mono\"\i dale symm\'etrique
des foncteurs $[n] \longrightarrow Fin$ qui envoient $0$ sur $\{1\}$. Ainsi, pour 
une cat\'egorie mono\"\i dale symm\'etrique $\mathcal{A}$, on peut pr\'ecomposer le foncteur
$\underline{Hom}^{\otimes}(-,\mathcal{A})$, qui envoie $\mathcal{B}$ sur la cat\'egorie
des foncteurs mono\"\i daux symm\'etriques de $\mathcal{B}$ dans $\mathcal{A}$, pour obtenir
$$N_{\Gamma}(\mathcal{A}) : \Gamma \longrightarrow \mathbf{Cat}.$$ 
On identifie alors $\mathbf{Cat}$ \`a la sous-cat\'egorie pleine de $\s -\mathbf{Cat}$ form\'ee
des $\s$-cat\'egories dont les ensembles simpliciaux de morphismes sont discrets, et on 
voit ainsi $N_{\Gamma}(\mathcal{A})$ comme un objet de $\scps$.
On v\'erifie alors que $N_{\Gamma}(\mathcal{A})$ est une
$\s$-CMS: en effet, pour tout $[n]$ la cat\'egorie
mono\"\i dale symm\'etrique 
$Hom^{0}([n],Fin)$ s'\"ecrit comme une somme, dans la $2$-cat\'egorie des cat\'egories mono\"\i dales
symm\'etriques
$$\coprod_{s_{i} : [n] \rightarrow [1]} Hom^{0}([1],Fin) \simeq Hom^{0}([n],Fin).$$
Ceci implique donc que pour tout $n$ le morphisme
$$N_{\Gamma}(\mathcal{A})([n]) \longrightarrow N_{\Gamma}(\mathcal{A})([1])^{n}$$
est une \'equivalence de cat\'egories. Nous avons donc construit ainsi 
un foncteur $N_{\Gamma}$, de la cat\'egorie des cat\'egories mono\"\i dales symm\'etriques vers
la sous-cat\'egorie de $\scms$ form\'ee des foncteurs qui sont niveaux par niveaux des cat\'egories. 
Ce foncteur induit alors un foncteur sur les cat\'egories homotopiques
$$N_{\Gamma} : \mathrm{Ho}(\mathbf{Cat}^{\otimes}) \longrightarrow \mathrm{Ho}(\scms).$$
On peut montrer (nous le ferons pas) que ce foncteur est pleinement fid\`ele et que son image 
essentielle consiste en toute les $\s$-CMS $A$ telles que 
$A([1])$ soit \'equivalente \`a une cat\'egorie. Une construction 
dans l'autre sens, fournissant un inverse du foncteur $N_{\Gamma}$, sera pr\'esent\'ee plus loin.

\item \textbf{Spectres connectifs.} Soit $M$ une cat\'egorie de mod\`eles point\'ee et engendr\'ee
par cofibrations et $X\in M$ un objet. On d\'efinit un foncteur
$$\Gamma^{op} \longrightarrow M$$
en envoyant $[n]$ sur l'objet $X^{[n]}$, o\`u $X^{[n]}$ d\'esigne l'exponentielle 
de $X$ par l'ensemble point\'e $[n]$. En d'autre termes, 
ce foncteur envoie $[n]$ sur $X^{n}$, et un morphisme $u : [n] \longrightarrow [m]$
sur le morphisme $u^{*} : X^{m} \longrightarrow X^{n}$
obtenu en envoyant la composante $i$ sur la composante $u(i)$ (en prenant comme
convention que cette composante est $*$ lorsque $u(i)=0$). Fixons un foncteur de r\'esolution
fibrante $id \rightarrow R_{*}$ sur $M$ (au sens de \cite[\S 5]{ho}). On d\'efinit alors, pour $Y \in M$ un 
foncteur 
$$\begin{array}{ccc}
\Gamma & \longrightarrow & \mathbf{SEns} \\
 n  & \mapsto & Hom(Q(X^{[n]}),R_{*}(Y)),
\end{array}$$
 o\`u $Q$ est un foncteur de remplacement cofibrant dans $M$.
On composant avec le foncteur $\Pi_{\infty} : \mathbf{SEns} \longrightarrow \s -\mathbf{Cat}$ on trouve
un foncteur 
$$\begin{array}{cccc}
M(X,Y) : &  \Gamma & \longrightarrow & \mathbf{SEns} \\
 & [n]  & \mapsto & \Pi_{\infty}(Hom(X^{[n]},R_{*}(Y))).
 \end{array}$$
Cette construction d\'efinit un foncteur 
$$M(X,-) : M \longrightarrow \scps.$$
Lorsque $M$ est une cat\'egorie de mod\`eles stable alors le foncteur ci-dessus
$M(X,-)$ se factorise par la sous-cat\'egorie $\scms$. En effet, pour $Y\in M$, sa valeur en $[n]\in \Gamma$
est alors
$$M(X,Y)([n])=Hom(Q(X^{[n]}),R_{*}(Y)) \simeq \Pi_{\infty}(Map([n]\vee^{\mathbb{L}} X,Y)) 
\simeq \Pi_{\infty}(Map(X,Y))^{n},$$
comme cela se voit en utilisant le fait que les sommes homotopiques finies sont aussi des
produits homotopiques finis dans $M$. On trouve ainsi un foncteur bien d\'efini
$$M(X,-) : \mathrm{Ho}(M) \longrightarrow \mathrm{Ho}(\scms).$$
Appliquons cela \`a $M=Sp^{\Sigma}$ la cat\'egorie de mod\`eles des spectres sym\'etriques, 
et $X=S$ l'unit\'e pour structure monod\"\i dale $\wedge$ sur $M$
(voir \cite{hss}). On trouve ainsi un foncteur
$$\mathrm{Ho}(Sp^{\Sigma}) \longrightarrow \mathrm{Ho}(\scms).$$
On peut v\'erifier que ce foncteur est pleinement fid\`ele lorsqu'on le restreint \`a la sous-cat\'egorie
pleine des spectres connectifs. De plus, son image essentielle consiste en toutes les $\s$-SMS
$A$ v\'erifiant les deux conditions suivantes:
\begin{enumerate}
\item La cat\'egorie $[A([1])]$ est un groupo\"\i de.
\item Le foncteur $\Gamma \longrightarrow Ens$ qui \`a 
$[n]$ associe l'ensemble des classes d'isomorphismes de $[A([n])]$
est le nerf d'un groupe ab\'elien (i.e. de la forme $N_{\Gamma}(E)$ pour $E$ un groupe ab\'elien, 
voir notre exemple $1$).
\end{enumerate}
Tout ceci d\'ecoule de la relation bien connue entre $\Gamma$-espaces tr\`es sp\'eciaux
et spectres connectifs (voir \cite{mmss}), ainsi que de la relation entre 
ensembles simpliciaux et $\s$-groupo\"\i des rappel\'ee \`a la fin de notre \S 1.1. 

\end{enumerate}

Une $\s$-CMS $A$ poss\`ede deux objets sous-jacents fondamentaux. Tout d'abord
la $\s$-cat\'egorie $A([1])$, qui est appel\'ee la $\s$-cat\'egorie \emph{sous-jacente} 
\`a $A$. Notons que $A([1])$ est naturellement muni d'une structure
de mono\"\i de commutatif dans $\mathrm{Ho}(\s -\mathbf{Cat})$. En effet, la loi de mono\"\i de
est induite par le diagramme suivant
$$\xymatrix{A([1])^{2} & A([2]) \ar[r] \ar[l]_-{\sim} & A([1]),}$$
o\`u le morphisme de gauche est $s_{1}\times s_{2}$, et le morphisme de droite est 
induit par $p : [2] \rightarrow [1]$ qui envoie $1$ et $2$ sur $1$. Ce diagramme
induit un morphisme bien d\'efini dans $\mathrm{Ho}(\s -\mathbf{Cat})$
$$A([1]) \times A([1]) \longrightarrow A([1]).$$
On v\'erifie, \`a l'aide des valeurs de $A$ sur $[3]$ et de ses fonctorialit\'es, que 
cette loi de mono\"\i de est associative, unitaire et commutative. 

Le second objet sous-jacent \`a $A$ est le foncteur $[A] : \Gamma \longrightarrow \mathbf{Cat}$. 
Ce foncteur est encore une $\s$-CMS car $A \mapsto [A]$ commute aux produits finis. 
Ainsi, d'apr\`es le second exemple ci-dessus $[A]$ s'\'ecrit comme
$N_{\Gamma}(\mathcal{A})$ pour une cat\'egorie mono\"\i dale sym\'etrique
$\mathcal{A}$. Donnons une description plus d\'etaill\'ee de cette structure. 
Pour cela on peut clairement supposer que $A=[A]$ (cela simplifie les notations).
On choisit des adjoints \`a droite des \'equivalences
$$c : A([2]) \longrightarrow A([1])^{2} \qquad a : A([3]) \longrightarrow A([1]).$$
Ces adjoints seront not\'es respectivement $d$ et $b$, et leurs counit\'es
$$k : c\circ d \Rightarrow id \qquad h : b\circ a \Rightarrow id.$$
On d\'efinit le foncteur $-\otimes - : A([1]) \times A([1]) \longrightarrow A([1])$
par la composition
$$\xymatrix{
A([1])\times A([1]) \ar[r]^-{d} & A([2]) \ar[r]^-{p} & A([1]),}$$
o\`u $p : [2] \rightarrow [1]$ est tel que $p(1)=p(2)=1$.
Par d\'efinition le morphisme $p$ est $\Sigma_{2}$-invariant, et $d$ est l'adjoint \`a droite
de $A([2]) \longrightarrow A([1])\times A([1])$ qui est lui aussi $\Sigma_{2}$-invariant
(au sens strict). 
Ainsi, $d$ est muni d'une unique structure de foncteur $\Sigma_{2}$-\'equivariant de sorte
\`a ce que la counit\'e $k$ soit $\Sigma_{2}$-invariante. Ceci muni 
la loi $\-\otimes -$ d'une contrainte de commutativit\'e bien d\'etermin\'ee et ne d\'ependant 
que des choix de $d$ et de $k$. 
Pour la contrainte d'associativit\'e on consid\`ere le diagramme commutatif suivant
$$\xymatrix{
A([3]) \ar[r]^-{q} \ar[d]_{a'} & A([2]) \ar[d]^-{c}\ar[r]^-{p} & A([1]) \\
A([2])\times A([1]) \ar[r]_-{p\times id} \ar[d]_-{c\times id} & A([1]) \times A([1]) & \\
A([1])\times A([1]) \times A([1]). & & }
$$
Le morphisme $q$ est induit par $q : [3] \rightarrow [2]$ avec 
$q(1)=1$, $q(2)=q(3)=2$, et le morphisme $a'$ est induit par 
$[3] \rightarrow [2]$ qui envoie $3$ sur $0$ et $[3] \rightarrow [1]$ 
qui envoie $1$ et $2$ sur $0$.
Le foncteur $(-\otimes -)\otimes - : A([1])^{3} \longrightarrow A([1])$
est par d\'efinition $p\circ d\circ (p\times id)\circ (d\times id)$. 
De plus, l'adjoint $b$ d\'efinit un unique adjoint \`a droite $b$' de $a'$, donn\'e par la composition
$b':=b\circ (c\times id)$, et avec l'unique counit\'e $l : a'\circ b' \rightarrow id$
telle que $(c\times id).l = h.(c\times id)$. 
De cette fa\c{c}on, on obtient un 
isomorphisme
$$m : d\circ (p\times id) \longrightarrow q\circ b',$$
adjoint de $p\times id \longrightarrow c\circ q\circ b'=(p\times id)\circ a'\circ b'$, lui-m\^eme
induit par l'unit\'e $l^{-1} : id \rightarrow a'\circ b'$. On construit de cette
fa\c{c}on un isomorphisme
$$\xymatrix{
p\circ d\circ (p\times id)\circ (d\times id) \ar[r]^-{m} &  p\circ q\circ b' \circ (d\times id)
=p\circ q\circ b\circ (c\times id)\circ (d\times id) \ar[r]^-{k} & 
p\circ q\circ b}.
$$
En d'autre termes, nous avons construit un isomorphisme
$$(-\otimes -)\otimes -  \simeq p\circ q\circ b,$$
qui ne d\'epend cette fois que du choix de $b$, $d$, $k$ et $h$. De fa\c{c}on duale on construit un
isomorphisme
$$-\otimes (-\otimes -) \simeq p\circ q\circ b.$$
En composant ces deux isomorphismes on trouve une contrainte d'associativit\'e
$$(-\otimes -)\otimes -  \simeq -\otimes (-\otimes - )$$
uniquement d\'etermin\'ee par les choix de $b$, $d$, $h$ et $k$. On v\'erifie alors que 
ces contraintes de commutativit\'e et d'associativit\'e munissent 
$A([1])$ d'une structure mono\"\i dale sym\'etrique. Une unit\'e pour cette structure
est alors donn\'ee par le choix d'un objet de $A([1])$ dans l'image du foncteur
$A([0]) \rightarrow A([1])$. Nous laissons les d\'etails restant aux lecteurs
(c'est \`a dire de v\'erifier l'axiome du pentagone et les compatibilit\'es entre
toutes les contraintes ainsi construites). Nous laissons aussi le soin au lecteur
de montrer qu'un $\otimes-\s$-foncteur $f : A \longrightarrow B$ induit une structure mono\"\i dale sym\'etrique
sur le foncteur induit $[f] : [A([1])] \longrightarrow [B([1])]$. \\

Venons-en maintenant \`a la structure de mod\`eles localis\'ee
sur $\scps$. Nous appellerons cette structure la \emph{structure de mod\`eles sp\'eciale sur 
$\scps$}, pour rappeler que ces objets locaux ne sont autre que l'analogue des 
$\Gamma$-espaces sp\'eciaux de \cite{sch}. 

\begin{prop}\label{p2}
Il existe sur $\scps$ une unique structure de mod\`eles, appel\'ee 
la \emph{structure sp\'eciale}, v\'erifiant les conditions suivantes.
\begin{enumerate}
\item Les cofibrations sont celles de la structure projective niveaux par niveaux de $\scps$. 
\item Un morphisme $A \longrightarrow B$ est une \'equivalence pour la structure 
sp\'eciale si et seulement si pour toute $\s$-CMS $C$ le morphisme induit
$$Map(B,C) \longrightarrow Map(A,C)$$
est un isomorphisme dans $\mathrm{Ho}(\mathbf{SEns})$ (o\`u les espaces de morphismes sont ici 
calcul\'es pour la structure niveaux par niveaux).
\end{enumerate}
\end{prop}

\textbf{Preuve --} Il s'agit d'un argument standard de localisation de Bousfield \`a gauche.
Pour cela, notons $I=\{X \rightarrow Y\}$ un ensemble de cofibrations g\'en\'eratrices
pour la cat\'egorie de mod\`eles $\s -Ca^{pr}t$. Pour $[n]\in \Gamma$ nous notons
$h^{n} : \Gamma \longrightarrow Ens$ le foncteur corepr\'esent\'e par $[n]$
($h^{n}([m]):=\Gamma([n],[m])$). On dispose d'un morphisme dans $\scps$
$$\coprod_{1\leq i\leq n}h^{1} \longrightarrow h^{n}$$
qui est induit par les $n$ morphismes $s_{i} : [n] \longrightarrow [1]$. Pour $n=0$ 
ce morphisme n'est autre que $\emptyset \rightarrow h^{0}=*$. 

On consid\`ere alors
l'ensemble de morphismes dans $\scps$
$$S:=\lbrace 
\left( X\rightarrow Y \right) \Box \left( \coprod_{1\leq i\leq n}h^{1} \longrightarrow h^{n} \right) \rbrace,$$
pour tout $n\geq 0$ et $X \rightarrow Y$ dans $I$ (voir \cite{ho} pour la d\'efinition du $\Box$ de deux
morphismes). La cat\'egorie de mod\`eles $\s -\mathbf{Cat}$ \'etant combinatoire il en est de m\^eme de 
$\scps$. On peut donc consid\'erer la structure de mod\`eles localis\'ee \`a gauche $L_{S}^{B}\scps$. 
Il nous reste \`a montrer que $L_{S}^{B}\scps$ v\'erifie les deux conditions
de la proposition. La premi\`ere est claire par d\'efinition de la structure localis\'ee. 
Quant \`a la seconde, il suffit de montrer que les objets $S$-locaux dans $\mathrm{Ho}(\scps)$
ne sont autre que les $\s$-CMS. Or, un objet $C\in \mathrm{Ho}(\scps)$ est $S$-local
si et seulement si pour tout $f : A \rightarrow B$ \'el\'ement de $S$ le morphisme
$Map(B,C) \longrightarrow Map(A,C)$ est une \'equivalence. Mais par d\'efinition 
des morphismes de $S$ les objets $S$-locaux sont donc les $C$ tels que pour tout $n$ et 
toute cofibration g\'en\'eratrice $X\rightarrow Y$ de $\s -\mathbf{Cat}^{pr}$
le morphisme
$$Map_{\s -\mathbf{Cat}^{pr}}(Y,C([n])) \longrightarrow Map_{\s -\mathbf{Cat}^{pr}}(Y,C([1])^{n}) 
\times^{h}_{Map_{\s -\mathbf{Cat}^{pr}}(X,C([1])^{n}}Map_{\s -\mathbf{Cat}^{pr}}(X,C([n]))$$
est une \'equivalence. 
Comme les cofibrations $X\rightarrow Y$ engendrent $\s -\mathbf{Cat}^{pr}$ ceci est bien \'equivalent au fait que
pour tout $n$ le morphisme
$$C([n]) \longrightarrow C([1])^{n}$$
soit une \'equivalence de $\s$-cat\'egories.
\hfill $\Box$ \\

Remarquons que les objets fibrants de la structure de mod\`eles localis\'ee de la proposition 
pr\'ec\'edente sont les foncteurs $\Gamma \longrightarrow \s -\mathbf{Cat}$ qui sont d'une part
des $\s$-CMS et d'autre part fibrant niveau par niveau.

\begin{df}\label{d7}
La cat\'egorie $\scps$, munie de sa structure de mod\`eles sp\'eciale, sera not\'ee
$\scps_{sp}$. Pour les distinguer, les espaces de morphismes de $\scps$ seront 
not\'es $Map^{pr-\otimes}$, et ceux de $\scps_{sp}$ seront not\'es $Map^{\otimes}$. 
\end{df}

En d'autres termes, la d\'efinition pr\'ec\'edente nous dit 
$$Map^{\otimes}(A,B):=Map^{pr-\otimes}(A,RB)$$
o\`u $RB$ est un remplacement fibrant de $B$ pour la structure sp\'eciale. \\

L'adjonction de Quillen $id : \scps \longleftrightarrow \scps_{sp} : id$
induit un foncteur pleinement fid\`ele
$$\mathrm{Ho}(\scps_{sp}) \longrightarrow \mathrm{Ho}(\scps)$$
qui identifie le membre de gauche \`a la sous-cat\'egorie pleine
de $\mathrm{Ho}(\scps)$ form\'ee des $\s$-CMS. Nous identifierons toujours 
$\mathrm{Ho}(\scps_{sp})$ \`a cette sous-cat\'egorie pleine, et ce de fa\c{c}on implicite par 
la suite.
De plus, le foncteur de remplacement 
fibrant $R$ de $\scps_{sp}$ induit un adjoint \`a gauche au foncteur d'inclusion
$$R : \mathrm{Ho}(\scps) \longrightarrow \mathrm{Ho}(\scps_{sp}).$$
Ainsi, pour tout foncteur $A : \Gamma \longrightarrow \s -\mathbf{Cat}$, le morphisme
d'adjonction $A \longrightarrow R(A)$, vu comme morphisme dans $\mathrm{Ho}(\scps)$, 
exhibe $R(A)$ comme la \emph{$\s$-CMS engendr\'ee par $A$}, au sens o\`u 
elle est universelle recevant un morphisme de $A$. Ceci permet de construire 
de nombreux exemples de $\s$-CMS \`a partir de la simple donn\'e d'un 
diagramme de $\s$-cat\'egories parametr\'e par $\Gamma$. 

\begin{df}\label{d8}
La \emph{$\s$-cat\'egorie des $\s$-CMS} est $L(\scps_{sp})$. Elle sera not\'ee
$\uscms$ (ou encore $\uscms_{\mathbb{U}}$ si l'on se restreint
aux $\s$-CMS $\mathbb{U}$-petites).
\end{df}

Notons imm\'ediatement que $\uscms_{\mathbb{U}}$ est une $\s$-cat\'egorie
$\mathbb{V}$-petite (pour $\mathbb{U}\in \mathbb{V}$) 
et $\mathbb{U}$-localement pr\'esentable. Les crit\`eres
d'existence d'adjoints de notre proposition \ref{p1} s'appliquent donc. De plus, on sait que 
$\uscat_{\mathbb{U}}$ poss\`ede toutes les 
$\mathbb{U}$-limites et $\mathbb{U}$-colimites. 

A titre indicatif, le foncteur d'\'evaluation en $[1]$ induit un $\s$-foncteur
$$\uscms \longrightarrow \uscat.$$
Ce foncteur \'etant induit par le foncteur de Quillen \`a droite
$\scps \longrightarrow \s -\mathbf{Cat}$ (\'evaluation en $[1]$), on voit que
le $\s$-foncteur ci-dessus poss\`ede un adjoint \`a gauche
$$Fr^{\otimes}(-) : \uscat \longrightarrow \uscms$$
appel\'e \emph{$\s$-CMS libre}. Explicitement ce $\s$-foncteur
est induit en localisation par la construction qui prend une $\s$-cat\'egorie
$A$ et qui l'envoie sur $R(h^{1}\times A)$, o\`u $h^{1} : \Gamma \longrightarrow Ens$
est corepr\'esent\'e par $[1]$, et o\`u $R$ est un foncteur de remplacement fibrant 
dans $\scps_{sp}$. \\

\textbf{Localisation des cat\'egories mono\"\i dales sym\'etriques.}
Terminons cette premi\`ere partie en donnant un proc\'ed\'e relativement efficace de
construction de $\s$-CMS. 
On se donne une cat\'egorie mono\"\i dale sym\'etrique
$C$, munie d'un ensemble de morphismes $W$ (et contenant les isomorphismes). 
On suppose que la structure mono\"\i dale
$\otimes : C\times C \longrightarrow C$ envoie $W\times W$ dans $W$. On construit alors
une $\s$-CMS $LC$ de la fa\c{c}on suivante. Par le foncteur $N_{\Gamma} : \mathbf{Cat}^{\otimes} \longrightarrow
\scms$ la cat\'egorie mono\"\i dale sym\'etrique $C$ donne lieu \`a un 
foncteur $\Gamma \longrightarrow \mathbf{Cat}$, qui est de plus une $\s$-CMS. Le sous-ensemble
de morphismes $W$ en lui-m\^eme envoy\'e sur un sous-foncteur de $N_{\Gamma}$ (c'est \`a dire
pour tout $[n]$ on a $W_{n}$ en ensemble de morphismes dans $N_{\Gamma}$, fonctoriel en 
$[n]$). En appliquant une version fonctorielle de la localisation (par exemple
celle d\'efinie par push-out, voir \S 1.2), on trouve
un nouveau foncteur
$$L_{W}^{\otimes}C : \Gamma \longrightarrow \s -\mathbf{Cat}$$
qui en $[n]$ vaut $L_{W_{n}}N_{\Gamma}(C)$. Le fait que $L$ soit compatible aux produits finis (voir \S 1.2)
implique de plus que $L_{W}^{\otimes}C$ est un $\s$-CMS, qui est de plus
$\mathbb{U}$-petite si $C$ l'est. Par construction on dispose d'un
morphisme $l : C \longrightarrow L_{W}^{\otimes}C$ dans $\mathrm{Ho}(\scms)$, et il n'est pas difficile de voir que 
ce morphisme v\'erifie la propri\'et\'e universelle suivante:
pour tout $A\in \mathrm{Ho}(\scms)$, le morphisme 
$$[L_{W}^{\otimes}C,A] \longrightarrow [C,A]$$
est injectif et son image consiste en tous les $\otimes-\s$-foncteurs $f : C \longrightarrow A$ tels que
$f([1]) : C([1])\simeq C \longrightarrow A([1])$ envoie $W$ dans les isomorphismes de $[A([1])]$. 

La construction $(C,W) \mapsto L_{W}^{\otimes}C$ fournit donc un 
foncteur de la cat\'egorie homotopique des cat\'egories mono\"\i dales sym\'etriques
munies d'un ensemble de morphismes $W$ compatible \`a la structure mono\"\i dale, vers 
la cat\'egorie homotopique des $\s$-CMS. Nous appliquerons cette construction \`a la situation
suivante. Soit $M$ une cat\'egorie de mod\`eles mono\"\i dale sym\'etrique. 
Soit $C_{0} \subset M^{c}$ une sous-cat\'egorie pleine d'objets cofibrants qui est de plus stable
par la structure mono\"\i dale ainsi que par \'equivalences dans $M$,  
mais qui ne contient pas forc\'ement l'unit\'e. Soit 
alors $C$ la sous-cat\'egorie pleine de $M$ form\'ee de $C_{0}$ et de l'unit\'e
(noter que cette derni\`ere n'est pas forc\'ement cofibrante). Notons
$W$ la trace des \'equivalences de $M$ sur $C$. Alors $W \otimes W \subset W$ et on peut donc
appliquer la construction pr\'ec\'edente (voir \cite{koto} pour des d\'etails sur le traitement
de l'unit\'e). On obtient ainsi une $\s$-CMS $L_{W}^{\otimes}C$. De cette fa\c{c}on, toute cat\'egorie
de mod\`eles mono\"\i dale sym\'etrique $M$ donne lieu \`a une $\s$-CMS $L_{W}^{\otimes}M^{c}$. 
La $\s$-cat\'egorie
sous-jacente de $L_{W}^{\otimes}M^{c}$ est naturellement \'equivalente \`a la localis\'ee
de $M^{c}$, la sous-cat\'egorie des objets cofibrants dans $M$. Comme le foncteur 
d'inclusion et le remplacement cofibrant induisent des isomorphismes inverses l'un de
l'autre dans $\mathrm{Ho}(\s -\mathbf{Cat})$ entre $LM^{c}$ et $LM$, on voit que la $\s$-cat\'egorie
sous-jacente \`a $L_{W}^{\otimes}M^{c}$ est naturellement \'equivalente \`a
$LM$. En d'autres termes, $L_{W}^{\otimes}M^{c}$ fournit une
structure de $\s$-CMS sur $LM$. 

\begin{df}\label{d9}
Pour une cat\'egorie de mo\`eles mono\"\i dale sym\'etrique $M$ nous noterons
$L^{\otimes}M:=L^{\otimes}_{W}M^{c} \in \scms$ (en gardant les notations
pr\'ec\'edentes).
\end{df}

\subsection{$\s$-Cat\'egories mono\"\i dales rigides}

Pour une $\s$-CMS $A$ nous noterons $[A]$ le foncteur 
induit $\Gamma \longrightarrow \mathbf{Cat}$ qui envoie $[n]$ sur $[A([n])]$. Comme nous l'avons
expliqu\'e au paragraphe pr\'ec\'edent ce foncteur permet de construire une structure
mono\"\i dale sym\'etrique sur la cat\'egorie $[A([1])]$. Par abus de notation cette
cat\'egorie mono\"\i dale sym\'etrique sera not\'ee $[A]$. Nous utiliserons
en particulier les contraintes d'associativit\'e, de commutativit\'e et 
d'unit\'e pour $[A]$, et souvent de fa\c{c}on implicite. Notons aussi 
que le foncteur $\otimes$ poss\`ede un rel\`evement naturel 
en un $\s$-foncteur 
$$-\otimes - : A([1]) \times A([1]) \longrightarrow A([1]),$$
d\'efini comme le compos\'e
$$\xymatrix{
A([1]) \times A([1]) \ar[r]^-{d} & A([2]) \ar[r] & A([1]),}$$
avec $d$ un adjoint \`a droite de l'\'equivalence $A([2]) \longrightarrow A([1])^{2}$.
Par la suite nous utiliserons \`a plusieurs reprises le $\s$-foncteur 
$\otimes$ d\'efini sur $A([1])$. 

\begin{prop}\label{p3}
Soit $A$ une $\s$-cat\'egorie mono\"\i dale sym\'etrique et $x\in A$ un objet de $A([1])$.
Les conditions suivantes sont \'equivalentes.
\begin{enumerate}
\item Il existe un objet $y\in A([1])$, et un morphisme $t : x\otimes y \rightarrow \mathbf{1}$
dans $[A([1])]$, tel que pour tout objet $z,z'\in A([1])$ le morphisme
$$\xymatrix{
A([1])(z,y\otimes z') \ar[r]^-{x\otimes -} & A([1])(x\otimes z,x\otimes (y\otimes z')) \ar[r]^-{\sim} & 
A([1])(x\otimes z,(x\otimes y)\otimes z')
\ar[r]^-{t_{*}}
& A([1])(x\otimes z,z')}$$
soit un isomorphisme dans $\mathrm{Ho}(\mathbf{SEns})$.
\item 
Il existe un objet $y\in A([1])$, et un morphisme $t : x\otimes y \rightarrow \mathbf{1}$
dans $[A([1])]$, tel que pour tout objet $z,z'\in A([1])$ l'application
$$\xymatrix{
[z,y\otimes z'] \ar[r]^-{x\otimes -} & [x\otimes z,x\otimes (y\otimes z')] \ar[r]^-{\sim} & 
[x\otimes z,(x\otimes y)\otimes z']
 \ar[r]^-{t_{*}}
& [x\otimes z,z']}$$
soit une bijection.
\item Il existe un objet $y\in A([1])$, et deux morphismes $t : x\otimes y \rightarrow \mathbf{1}$, 
et $u : \mathbf{1} \longrightarrow y\otimes x$, tels que les morphismes compos\'es suivants
$$\xymatrix{
y \ar[r]^-{u\otimes y} & (y\otimes x)\otimes y \ar[r]^-{\sim} & y\otimes (x\otimes y) 
\ar[r]^-{y\otimes t} & y} \qquad
\xymatrix{x \ar[r]^-{x\otimes u} & x\otimes (y\otimes x) \ar[r]^-{\sim} &  
(x\otimes y)\otimes x \ar[r]^-{t\otimes x} & x}$$
soient des identit\'es dans $[A([1])]$.
\item Il existe un objet $y\in A([1])$, et un morphisme $u : \mathbf{1} \rightarrow y\otimes x$
dans $[A([1])]$, tel que pour tout objet $z,z'\in A([1])$ le morphisme
$$\xymatrix{
A([1])(x\otimes z,z') \ar[r]^-{y\otimes -} & A([1])(y\otimes (x\otimes z),y\otimes z') 
\ar[r]^-{\sim} & A([1])((y\otimes x)\otimes z,y\otimes z') 
\ar[r]^-{u^{*}}
& A([1])(z,y\otimes z')}$$
soit un isomorphisme dans $\mathrm{Ho}(\mathbf{SEns})$.
\item L'objet $x$ est rigide dans la cat\'egorie mono\"\i dale sym\'etrique $[A]$.
\end{enumerate}
\end{prop}

\textbf{Preuve --} $(1) \Rightarrow (2)$ est clair, il suffit de prendre 
les $\pi_{0}$ des ensembles simpliciaux de morphismes dans $A([1])$. 
Pour voir que $(2)$ implique $(3)$ on applique la bijection 
$[z,y\otimes z'] \simeq [x\otimes z,z']$ de l'\'enonc\'e au cas
$z=\mathbf{1}$ et $z'=x$. L'identit\'e dans $[x,x]$ induit alors
un \'el\'ement $u \in [\mathbf{1},y\otimes x]$. Il est facile de voir par construction que
$u$ v\'erifie les conditions de $(3)$. Pour voir que $(3)$ implique $(4)$ on construit un morphisme
en sens inverse en consid\'erant le morphisme compos\'e
$$\xymatrix{
A([1])(z,y\otimes z') \ar[r]^-{x\otimes -} & A([1])(x\otimes z,x\otimes y\otimes z') \ar[r]^-{t_{*}}
& A([1])(x\otimes z,z').}$$
Les conditions sur $u$ et $t$ impliquent que ces deux morphismes sont inverses l'un de l'autre
dans $\mathrm{Ho}(\mathbf{SEns})$. Un argument sym\'etrique  \`a celui 
que l'on vient de donner pour $(1) \Rightarrow (4)$ montre aussi que $(1) \Rightarrow (4)$.
Enfin, $(5)$ est par d\'efinition une r\'e\'ecriture de la condition $(3)$.
\hfill $\Box$ \\

\begin{df}\label{d10}
Soit $A$ une $\s$-CMS. Un objet $x\in A([1])$ est \emph{rigide}, ou encore
\emph{dualisable}, s'il v\'erifie les conditions \'equivalentes de la proposition
\ref{p3}.
\end{df}

Les objets rigides poss\`edent les propri\'et\'es formelles suivantes.

\begin{prop}\label{p4}
Soit $A$ une $\s$-CMS.
\begin{enumerate}
\item Les objets rigides sont stables par isomorphismes
dans $[A([1])]$. 
\item L'unit\'e $\mathbf{1}$ est un objet rigide. 
\item Tout objet 
inversible pour la structure mono\"\i dale
$\otimes$ est un objet rigide.
\item Les objets rigides sont stables par la structure mono\"\i dale
$\otimes$.
\item Si $f : A \longrightarrow B$ est un $\otimes-\s$-foncteur,
alors l'image par $f$ d'un objet rigide est rigide.
\item Pour un objet rigide fix\'e $x$ de $A$, les donn\'ees de $y$, $t$ et $u$ comme dans la proposition 
\ref{p3} $(3)$ sont uniques \`a isomorphisme unique pr\`es dans $[A]$. 
\end{enumerate}
\end{prop}

\textbf{Preuve --} Pour voir ces propri\'et\'es on utilise le point 
$(5)$ de la proposition \ref{p3}. La proposition \ref{p4} sont alors
des faits bien connus pour les objets rigides dans une cat\'egorie
mono\"\i dale sym\'etrique. \hfill $\Box$ \\

Comme il en est l'usage, pour un objet $x$ rigide fix\'e dans une $\s$-CMS
$A$, nous noterons $x^{\vee}:=y$ un objet dual comme dans 
la proposition \ref{p3}. Lorsque nous utiliserons $x^{\vee}$ il est sous-entendu 
que les unit\'es et counit\'es $t$ et $u$ satisfaisant \ref{p3} $(3)$ ont aussi
\'et\'e fix\'ees. En d'autres termes, l'expression \emph{soit $x^{\vee}$ un dual 
de $x$} sera synonyme de \emph{soit $(y,u,t)$ un triplet comme dans 
la proposition \ref{p3} $(3)$}. \\

D'apr\`es la proposition \ref{p4} les objets rigides dans une 
$\s$-CMS $A$ forment \`a leur tour une sous-$\s$-SMC $A^{rig} \subset A$. 
La d\'efinition pr\'ecise de $A^{rig}$ est la suivante: pour $[n] \in \Gamma$, 
$A^{rig}([n])$ est la sous-$\s$-cat\'egorie pleine de $A([n])$ form\'ee
des objets dont toutes les images dans $A([1])$ par les morphismes
$s_{i} : [n] \rightarrow [1]$ sont des objets rigides. Ceci d\'efinit un sous-foncteur
$A^{rig}$ de $A$ qui est encore une $\s$-CMS. De plus, cette construction est 
fonctorielle pour les $\otimes-\s$-foncteurs et d\'efinit donc un endofoncteur
$A \mapsto A^{rig}$ de la cat\'egorie $\scms$. Cet endofoncteur est muni d'une
transformation naturelle vers l'identit\'e. De plus, le morphisme naturel
$(A^{rig})^{rig} \rightarrow A^{rig}$ est l'identit\'e.

\begin{df}\label{d11}
\begin{enumerate}
\item Une $\s$-CMS est \emph{rigide} si tous ses objets sont rigides. 
\item La sous-cat\'egorie pleine des $\s$-CMS rigides sera not\'ee
$\scms_{rig}$.
\item La sous-$\s$-cat\'egorie pleine de $\uscms$ form\'ee des $\s$-CMS rigides  
sera not\'ee $\uscms_{rig}$. 
\end{enumerate}
\end{df}

Le principal r\'esultat de ce paragraphe est le suivant. Il sera 
fondamental dans l'\'etude et la construction de morphismes de traces
du paragraphe suivant.

\begin{thm}\label{t1}
Le $\s$-foncteur
$$\uscms \hookrightarrow \uscms$$
qui \`a $A$ associe $A^{rig}$ poss\`ede un adjoint \`a gauche.
\end{thm}

\textbf{Preuve --} Fixons nous l'univers ambiant $\mathbb{U}$, et consid\'erons
le $\s$-foncteur
$$(-)^{rig} : \uscms_{\mathbb{U}} \hookrightarrow \uscms_{\mathbb{U}}.$$
Comme $\underline{scms}_{\mathbb{U}}\simeq L(\scps_{sp,\mathbb{U}})$, $(-)^{rig}$
est un $\s$-foncteur entre deux $\s$-cat\'egorie $\mathbb{U}$-localement pr\'esentables. 
Nous pouvons donc appliquer le crit\`ere de la proposition \ref{p1}. Pour d\'emontrer le th\'eor\`eme
il nous suffit donc de montrer que le $\s$-foncteur $A \mapsto A^{rig}$ commute aux 
$\mathbb{U}$-limites. Pour cela il suffit de montrer ind\'ependamment qu'il commute
aux produits (\'eventuellement infinis) et aux produits fibr\'es. De mani\`ere \'equivalente, 
il faut montrer que le foncteur $A \mapsto A^{rig}$, commute aux produits et aux 
produits fibr\'es homotopiques dans $\scms$. 

Pour voir que $A \mapsto A^{rig}$ commute aux produits il suffit de v\'erifier que pour
une famille de $\s$-CMS $\{A_{\alpha}\}$, un objet $(x_{\alpha})\in \prod_{\alpha}A_{\alpha}$
est rigide si et seulement si $x_{\alpha}\in A_{\alpha}$ est rigide pour tout $\alpha$, et 
la proposition \ref{p3} $(3)$ montre que tel est le cas. Il nous reste donc \`a montrer que 
$A \rightarrow A^{rig}$ commute aux produits fibr\'es homotopiques dans $\scms$. 
Soit donc $\xymatrix{A \ar[r] & C & \ar[l] B}$ un diagramme dans $\scms$, et supposons que
$A \rightarrow C$ soit une fibration dans $\scps$ (c'est \`a dire une fibration niveau 
par niveau). Notons $D:=A\times_{C}B$ le produit fibr\'e dans
$\scms$ du diagramme ci-dessus. On commence par remarquer que le morphisme induit
$A^{rig} \rightarrow C^{rig}$ est encore une fibration, ce qui se d\'eduit ais\'ement
du fait que pour tout $[n]$ la sous-$\s$-cta\'egorie pleine $A([n])^{rig} \subset A([n])$
soit stable par \'equivalences dans $A([n])$. Il nous faut donc montrer que le morphisme induit
$$D^{rig} \longrightarrow A^{rig}\times_{C^{rig}}B^{rig}$$
est une \'equivalence. Comme il s'agit d'un morphisme de 
$\s$-CMS il suffit de montrer que le $\s$-foncteur de $\s$-cat\'egories
$$\phi : D([1])^{rig} \longrightarrow A([1])^{rig}\times_{C([1])^{rig}}B([1])^{rig}$$
est une \'equivalence.
Or, le morphisme $\phi$
est un $\s$-foncteur pleinement fid\`ele, comme cela se voit en consid\'erant le diagramme
commutatif suivant
$$\xymatrix{
D([1])^{rig} \ar[dr] \ar[r] & A([1])^{rig}\times_{C([1])^{rig}}B([1])^{rig} \ar[d] \\
& D([1])=A([1])\times_{C([1])}B([1]).}$$
Il nous reste donc \`a montrer que $\phi$ est essentiellement surjectif. 
Pour cela on consid\`ere le
diagramme commutatif suivant de cat\'egories mono\"\i dales sym\'etriques
$$\xymatrix{
[D^{rig}]\ar[r] \ar[rd] & [A^{rig} \times_{C^{rig}}B^{rig}] \ar[d] \\
 & [A^{rig}] \times_{[C^{rig}]}[B^{rig}].}$$
D'apr\`es la proposition \ref{p3} $(5)$ ce diagramme est \'equivalent au diagramme
$$\xymatrix{
[D]^{rig}\ar[r] \ar[rd] & [A^{rig} \times_{C^{rig}}B^{rig}] \ar[d] \\
 & [A]^{rig} \times_{[C]^{rig}}[B]^{rig}.}$$
On remarque alors que le foncteur $[A^{rig} \times_{C^{rig}}B^{rig}] \longrightarrow  [A]^{rig} \times_{[C]^{rig}}[B]^{rig}$ induit une bijection sur les ensembles de classes d'isomorphismes d'objets. 
Ainsi, il nous suffit de montrer que le foncteur induit
$$[D]^{rig} \longrightarrow [A]^{rig} \times_{[C]^{rig}}[B]^{rig}$$
est essentiellement surjectif. Ce foncteur se factorise par le foncteur naturel
$$[A\times_{C}B]^{rig}=[D]^{rig}  \longrightarrow ([A]\times_{[C]}[B])^{rig}.$$
On voit donc qu'il nous reste \`a  montrer que le foncteur 
$$([A]\times_{[C]}[B])^{rig} \longrightarrow [A]^{rig} \times_{[C]^{rig}}[B]^{rig}$$
est essentiellement surjectif. En d'autres termes on s'est ramen\'es au cas
o\`u $A$, $B$ et $C$ sont des cat\'egories mono\"\i dales sym\'etriques. Dans ce
cas notons $p : A \longrightarrow C$ et $q : B \longrightarrow C$ les deux projections.
Pour deux objets $x \in A$ et $y\in B$, poss\'edant la m\^eme image, $z \in C$, 
choisissons des duaux $x^{\vee}$, $y^{\vee}$, et des
unit\'es 
$$u : \mathbf{1} \rightarrow x^{\vee}\otimes x \qquad v :  \mathbf{1} \rightarrow y^{\vee}\otimes y.$$
L'unicit\'e des duaux implique l'existence d'un unique isomorphisme 
$\alpha : p(x^{\vee}) \simeq q(y^{\vee})$, dans $C$, compatible
avec $p(u)$ et $q(v)$. Comme $A \longrightarrow C$ est une fibration
on peut trouver un isomorphisme $\beta : x^{\vee} \simeq x'$ dans $A$ avec $p(\beta)=\alpha$.
Dans ce cas, $x'$, muni de l'unit\'e compos\'ee
$$u' : \xymatrix{
\mathbf{1} \ar[r]^-{u} & x^{\vee}\otimes x  \ar[r]^-{\beta} & x'\otimes x,}$$
est un dual de $x$ dans $A$. Comme on a $p(x')=q(y^{\vee})$, les objets 
$x'$ et $y^{\vee}$ se recollent en un objet $(x',y^{\vee}) \in C$. On v\'erifie facilement que cet objet, 
et le morphisme
$$(u',v) : (\mathbf{1},\mathbf{1}) \longrightarrow
(x'\otimes x,y^{\vee},\otimes y) = (x',y^{\vee}) \otimes (x,y),$$
est un dual de $(x,y) \in C$. Ceci montre que $(x,y) \in C^{rig}$ est un ant\'ec\'edent de
$(x,y) \in A^{rig}\times_{C^{rig}}B^{rig}$, et termine donc la preuve du th\'eor\`eme. 
\hfill $\Box$ \\

L'adjoint \`a gauche du foncteur $A \mapsto A^{rig}$ sera not\'e
$A \mapsto Fr^{rig}(A)$. La $\s$-CMS $Fr^{rig}(A)$ est la $\s$-CMS rigide
engendr\'ee $A$. On dispose d'un morphisme bien d\'efini dans $\mathrm{Ho}(\scms)$
$$A \longrightarrow (Fr^{rig}(A))^{rig}$$
qui est tel que pour toute $\s$-CMS $B$ le morphisme induit 
$$Map^{\otimes}(Fr^{rig}(A),B) \longrightarrow 
Map^{\otimes}(Fr^{rig}(A)^{rig},B^{rig}) \longrightarrow Map^{\otimes}(A,B^{rig})$$
doit un isomorphisme dans $\mathrm{Ho}(\mathbf{SEns})$. 

\begin{lem}\label{l1}
Pour toute $\s$-CMS $A$, la $\s$-CMS $Fr^{rig}(A)$ est rigide. 
\end{lem}

\textbf{Preuve --} La propri\'et\'e universelle du morphisme
$A \longrightarrow Fr^{rig}(A)^{rig}$ rappel\'ee ci-dessus
implique en particulier le fait suivant: pour toute $\s$-CMS $B$, 
et $f : A \longrightarrow B$ un morphisme dans $\mathrm{Ho}(\scms)$, 
le morphisme $f$ se factorise par $B^{rig} \subset B$ si et seulement si $f$ se factorise, 
dans $\mathrm{Ho}(\scms)$, par le morphisme $A \rightarrow Fr^{rig}(A)^{rig} \subset Fr^{rig}(A)$
$$\xymatrix{
Fr^{rig}(A) \ar[r] & B \\
A. \ar[u] \ar[ru]_-{f} & }$$
Lorsque cette factorisation existe elle est alors unique, car le morphisme induit
$$[Fr^{rig}(A),B] \simeq [A,B^{rig}] \subset [A,B]$$
est injectif.
Appliqu\'e \`a $B=Fr^{rig}(A)^{rig}$ et au morphisme d'adjonction $A \rightarrow Fr^{rig}(A)^{rig}$
on trouve l'existence d'un diagramme commutatif
$$\xymatrix{
Fr^{rig}(A) \ar[r]^-{r} & Fr^{rig}(A)^{rig} \\
A. \ar[u] \ar[ru] & }$$
Le diagramme commutatif
$$\xymatrix{
Fr^{rig}(A) \ar[r]^-{i\circ r} & Fr^{rig}(A) \\
A \ar[u] \ar[ru] & }$$
et l'unicit\'e de la factorisation pour $B=Fr^{rig}(A)$ implique
que $i\circ r=id$. Ceci montre que l'inclucion $Fr^{rig}(A)^{rig} \subset Fr^{rig}(A)$
poss\`ede une section et est donc une \'equivalence. Ceci montre que $Fr^{rig}(A)$ est rigide. \hfill $\Box$

Le lemme pr\'ec\'edent et la propri\'et\'e universelle implique donc que le $\s$-foncteur
$$Fr^{rig} : \uscms \longrightarrow \uscms$$
se factorise par la sous-$\s$-cat\'egorie des $\s$-CMS rigides, et fournit alors un 
adjoint \`a gauche du $\s$-foncteur d'inclusion $\uscms_{rig} \hookrightarrow \uscms$. 

\begin{cor}\label{c1}
La sous-$\s$-cat\'egorie $\uscms_{rig} \subset \uscms$ est stable par 
limites et colimites. 
\end{cor}

\textbf{Preuve --} En effet le $\s$-foncteur d'inclusion
$\uscms_{rig} \subset \uscms$ poss\`ede un adjoint \`a gauche $Fr^{rig}$. Il poss\`ede 
de plus un adjoint \`a droite $(-)^{rig}$. \hfill $\Box$ \\

\begin{cor}\label{c2}
La $\s$-cat\'egorie $\uscms_{rig,\mathbb{U}}$ est 
presque $\mathbb{U}$-localement pr\'esentable. 
\end{cor}

\textbf{Preuve --} En effet le $\s$-foncteur d'inclusion
$\uscms_{rig,\mathbb{U}} \hookrightarrow \uscms_{\mathbb{U}}$ poss\`ede
un adjoint \`a gauche et $\uscms_{\mathbb{U}}$ est une elle m\^eme 
$\mathbb{U}$-localement pr\'esentable. \hfill $\Box$ \\

Nous allons maintenant voir que $\uscms_{rig,\mathbb{U}}$ est en r\'ealit\'e
une $\s$-cat\'egorie $\mathbb{U}$-localement pr\'esentable. Pour cela, notons
$\mathbb{F}:=Fr^{\otimes}(*)$ la $\s$-cat\'egorie libre engendr\'ee par un g\'en\'erateur. 
On consid\`ere alors le morphisme 
$$u : \mathbb{F} \longrightarrow Fr^{rig}(\mathbb{F}).$$
La cat\'egorie de mod\`eles sp\'eciale $\scps_{sp}$ \'etant 
$\mathbb{U}$-combinatoire, nous pouvons la localiser \`a gauche 
le long du seul morphisme $u$. Cette localisation sera not\'ee
$$\scps_{sp,rig}:=L^{B}_{u}\scps_{sp}.$$

\begin{prop}\label{p5}
Il existe un isomorphisme dans $\mathrm{Ho}(\s -\mathbf{Cat})$
$$\uscms_{rig}\simeq L(\scps_{sp,rig}).$$
\end{prop}

\textbf{Preuve --} A la vue des d\'efinitions des $\s$-cat\'egories
en jeu il suffit de montrer que les $\s$-CMS $A$ qui sont 
$u$-locales sont exactement les $\s$-cat\'egories rigides. Or, une telle $A$ est 
$u$-locale si et seulement si le morphisme
$$Map^{\otimes}(Fr^{rig}(\mathbb{F}),A) \longrightarrow Map^{\otimes}(\mathbb{F},A)$$
est un isomorphisme dans $\mathrm{Ho}(\mathbf{SEns})$. Par adjonction ce morphisme est lui-m\^eme
isomorphe au morphisme
$$Map_{\s -\mathbf{Cat}}(*,A^{rig}([1])) \longrightarrow Map_{\s -\mathbf{Cat}}(*,A([1])).$$
Ceci implique en particulier que l'inclusion 
$[*,A^{rig}([1])] \hookrightarrow [*,A([1])]$ est bijective, o\`u 
encore que $[A^{rig}([1])] \hookrightarrow [A([1])]$ est essentiellement surjectif. 
En d'autres termes $A$ est rigide si et seulement si elle est $u$-locale. \hfill $\Box$ \\

\begin{cor}\label{c3}
La $\s$-cat\'egorie $\uscms_{rig,\mathbb{U}}$ est 
$\mathbb{U}$-localement pr\'esentable. 
\end{cor}

\subsection{Traces et traces cyliques}

Soit $A$ une $\s$-CMS et $x$ un objet rigide de $A$. Notons $x^{\vee}$ un dual 
de $x$ dans $A$, muni de morphismes dans $[A]$
$$u : \mathbf{1} \longrightarrow x\otimes x^{\vee} \qquad
t : x^{\vee} \otimes x \longrightarrow \mathbf{1},$$
v\'erifiant les identit\'es triangulaires de la proposition \ref{p3} $(3)$. 
On dispose alors d'un morphisme
dans $\mathrm{Ho}(\mathbf{SEns})$
$$\xymatrix{
A([1])(x,x) \ar[r]^-{sim} & A([1])(\mathbf{1},x^{\vee}\otimes x) \ar[r]^-{t} & A([1])(\mathbf{1},\mathbf{1}).}$$
Ce morphisme est bien d\'efini dans $\mathrm{Ho}(\mathbf{SEns})$, et ne d\'epend pas des choix
de $x^{\vee}$, $u$ et $t$ (car un tel choix est unique \`a isomorphisme
unique pr\`es dans $[A]$).

\begin{df}\label{d12}
Pour $x$ un objet rigide dans une $\s$-CMS $A$ le morphisme d\'efini ci-dessus est 
le \emph{morphisme trace en $x$}. Il sera not\'e
$$Tr_{x} : A([1])(x,x) \longrightarrow A([1])(\mathbf{1},\mathbf{1}).$$
\end{df}

Il est facile de voir que l'image du dual $x^{\vee}$ d'un objet 
$x$ par un $\otimes-\s$-foncteur $f$ est un dual de $f(x)$. Cela implique en particulier
que pour tout $\otimes-\s$-foncteur $f : A \longrightarrow B$, et tout objet 
$x$ de $A([1])$ le diagramme suivant
$$\xymatrix{
A([1])(x,x) \ar[r]^-{Tr_{x}} \ar[d]_-{f} & A([1])(\mathbf{1},\mathbf{1}) \ar[d]^-{f} \\
B([1])(f(x),f(x)) \ar[r]_-{Tr_{f(x)}} & B([1])(\mathbf{1},\mathbf{1})}$$
commute dans $\mathrm{Ho}(\mathbf{SEns})$. \\

\noindent\textbf{Fonctorialit\'e des traces.} Dans le reste de ce paragraphe nous nous attacherons \`a montrer comment
les morphismes traces de la d\'efinition ci-dessus peuvent \^etre rendus
fonctoriels en $x$ et en $A$. Nous allons commencer par montrer qu'il existe
une unique fa\c{c}on de construire de tels morphismes, fonctoriel en $x$ et $A$
\`a la fois. Cependant, notre premi\`ere approche montrera que cette unicit\'e
est elle-m\^eme d\'etermin\'ee \`a des isomorphismes non-uniques pr\`es. En d'autres termes
nous montrerons que l'espace des traces est connexe sans pour autant montrer qu'il est 
contractile. Pour montrer que cet espace est effectivement contractile nous aurons besoin
d'une cons\'equence d'un th\'eor\`eme de Hopkins-Lurie concernant les $\s$-cat\'egories
de bordismes orient\'es. La contactibilit\'e de l'espace des traces sera alors utilis\'ee
pour construire des \emph{traces cycliques}, construction fondamentale dans
le caract\`ere de Chern que nous pr\'esenterons dans la prochaine section. \\

Nous d\'efinissons deux $\s$-foncteurs 
$$\scms_{rig} \longrightarrow \T.$$
Le premier envoie une $\s$-CMS rigide $A$ sur $Map(S^{1},\mathcal{I}(A([1]))$, 
o\`u $\mathcal{I}$ est le foncteur espace sous-jacent d\'ecrit dans \S 1.1. Le second 
envoie une $\s$-CMS rigide $A$ sur l'ensemble simplicial $A([1])(\mathbf{1},\mathbf{1})$
des endomorphismes de l'unit\'e. Comme l'unit\'e $\mathbf{1}$ n'est pas bien d\'efini
dans $A([1])$ il nous faut pr\'eciser cette construction. Pour cela on identifiera, 
\`a \'equivalences naturelles pr\`es, $A([1])(\mathbf{1},\mathbf{1})$ avec le produit
fibr\'e homotopique
$$Map_{\s -\mathbf{Cat}}(\Delta^{1},A([1])) \times_{Map_{\s -\mathbf{Cat}}(0\coprod 1,A([1]))}^{h}Map_{\s -\mathbf{Cat}}(0\coprod 1,A([0])),$$
o\`u $0\coprod 1 \subset \Delta^{1}$ est l'inclusion des deux objets de la cat\'egorie
$\Delta^{1}$. En prenant des mod\`eles fonctoriels des $Map_{\s -\mathbf{Cat}}$ et 
du produit fibr\'e homotopique on construit ainsi un $\s$-foncteur
$\uscms \longrightarrow \T$, que nous noterons
symboliquement $A \mapsto A([1])(\mathbf{1},\mathbf{1})$. 

On dispose ainsi de deux $\s$-foncteurs
$$\uscms_{rig} \longrightarrow \T,$$
 not\'es
respectivement 
$$A \mapsto \mathcal{LI}(A) \qquad A \mapsto End_{A}(\mathbf{1}).$$
Notons que $\mathcal{LI}(A)$ s'identifie \`a l'espace des
objets de $A([1])$ munis d'auto\'equivalences. En particulier, on dispose pour tout 
$x\in A([1])$ d'un morphisme bien d\'efini dans $\mathrm{Ho}(\mathbf{SEns})$
$$A([1])(x,x)^{eq} \longrightarrow \mathcal{LI}(A)$$
de l'espace des auto\'equivalences de $x$ dans $\mathcal{LI}(A)$. 

\begin{prop}\label{p6}
Il existe une transformation naturelle 
$$Tr : \mathcal{LI} \longrightarrow End_{-}(\mathbf{1})$$
telle que pour tout $A\in \s -\mathbf{Cat}_{rig}$ et tout objet $x\in A([1])$, le morphisme induit
$$A([1])(x,x)^{eq} \longrightarrow End_{A}(\mathbf{1})$$
sur l'espace des auto\'equivalences de $x$ soit \'egal, dans $\mathrm{Ho}(\mathbf{SEns})$, au morphisme
$Tr_{x}$ d\'efini dans la d\'efinition \ref{d12}. De plus, 
un tel morphisme $Tr$ est unique dans la cat\'egorie
$[\rh (\uscms_{rig},\mathbb{T})]$. 
\end{prop}

\textbf{Preuve --} Pour $X\in \mathbf{SEns}$
on pose 
$$\mathbb{B}(X):=Fr^{rig}(Fr^{\otimes}(\Pi_{\infty}(X))) \in \scms_{rig},$$
la $\s$-CMS rigide libre sur le $\s$-groupo\"\i de fondamental de $X$. Par d\'efinition, 
la $\s$-CMS rigide $\mathbb{B}(S^{1})$ corepr\'esente le $\s$-foncteur $\mathcal{LI}$.
Ainsi, le lemme de Yoneda implique qu'il existe un isomorphisme dans $\mathrm{Ho}(\mathbf{SEns})$ 
$$Map(\mathcal{LI},End_{-}(\mathbf{1}))\simeq
End_{\mathbb{B}(S^{1})}(\mathbf{1}).$$
Le morphisme d'adjonction $S^{1} \longrightarrow \mathcal{LI}(\mathbb{B}(S^{1}))$
d\'efini un $\s$-foncteur $B\mathbb{Z} \longrightarrow \mathbb{B}(S^{1})([1])$, et donc
un couple $(x,u)$, o\`u 
$x$ est un objet de $\mathbb{B}(S^{1})([1])$ et $u : x \simeq x$ est un 
auto\'equivalence de $x$. Notons $Z:=Tr_{x}(u) \in End_{\mathbb{B}(S^{1})}(\mathbf{1})$. 
On v\'erifie alors que la transformation naturelle d\'efinie par $Z$ (\`a travers
l'\'equivalence rappel\'ee ci-dessus) satisfait aux conditions de la proposition. 
\hfill $\Box$ \\

\begin{thm}\label{t2}
Notons $Tr \in Map(\mathcal{LI},End_{-}(\mathbf{1}))$
un morphisme comme dans la proposition \ref{p6}. Alors
$$\pi_{i}(Map(\mathcal{LI},End_{-}(\mathbf{1})),Tr)=0 \;  \forall \; i>0.$$
\end{thm}

\textbf{Preuve --} Comme dans la preuve de la proposition \ref{p6}, on a 
$$Map(\mathcal{LI},End_{-}(\mathbf{1}))\simeq
End_{\mathbb{B}(S^{1})}(\mathbf{1}).$$
Il s'agit donc de montrer que 
$$\pi_{i}(End_{\mathbb{B}(S^{1})}(\mathbf{1}),Z)=0,$$
o\`u $Z$ est d\'efini par $Tr_{x}(u)$ comme dans la preuve de la proposition \ref{p6}. 
Ceci est alors une cons\'equence du th\'eor\`eme de Hopkins-Lurie sur la
propri\'et\'e universelle des $\s$-cat\'egories mono\"\i dales sym\'etriques
de bordismes (voir \cite{lu2}). Nous ne reproduirons pas ces \'enonc\'es ici, signalons 
juste comment la propri\'et\'e universelle de la $\s$-cat\'egorie
$\mathbf{Bord}_{1}^{X}$, des $1$-bordismes orient\'es au-dessus d'un espace $X$, implique
notre r\'esultat. 

Le th\'eor\`eme \cite[2.4.18]{lu2} nous dit que la $\s$-cat\'egorie $\mathbb{B}(X)([1])$ 
est naturellement \'equivalente \`a la $\s$-cat\'egorie 
$\mathbf{Bord}_{1}^{X}$, des $1$-bordismes
orient\'es au-dessus de $X$ (o\`u l'on prend $\zeta$ le fibr\'e trivial de rang $1$ sur $X$). 
L'unit\'e
$\mathbf{1}$ de $\mathbb{B}(X)([1])$ correspond \`a travers cette \'equivalence 
au bordisme vide $\emptyset \rightarrow X$. De plus, l'espace des endomorphismes de cet objet 
$Map(\emptyset,\emptyset)$ s'identifie naturellement \`a l'espace classifiant des
vari\'et\'es compactes, orient\'ees de dimension $1$, munis de morphismes vers $X$. 
En d'autres termes on dispose d'une projection
$$\pi : Map(\emptyset,\emptyset) \longrightarrow \coprod_{n\geq 0}(B(\Sigma_{n}\ltimes (S^{1})^{n}),$$
o\`u $B(\Sigma_{n}\ltimes (S^{1})^{n})$ est l'espace classifiant des 
vari\'et\'es orient\'ees isomorphes \`a $\{1, \dots, n\}\times S^{1}$, 
une r\'eunion disjointe de $n$ copies de $S^{1}$ (pou $n=0$ cet espace
est un point). De plus, 
la fibre homotopique de $\pi$ en un point correspondant \`a un entier $n\geq 1$, est naturellement 
\'equivalente \`a l'espace des morphismes $Map(\{1, \dots, n\}\times S^{1},X)$, et 
l'action induit de $\Sigma_{n}\ltimes (S^{1})^{n}$ est celle induite par 
l'action naturelle sur $\{1, \dots, n\}\times S^{1}$. 

La cat\'egorie mono\"\i dale 
$[\mathbf{Bord}_{1}^{X}]\simeq [\mathbb{B}(X)]$ se d\'ecrit comme suit. Ses objets sont 
des vari\'et\'es compactes orient\'es de dimension $0$ et munies d'applications continues vers $|X|$
(o\`u l'on note ici $|X|$ la r\'ealisation g\'eom\'etrique de $X$). 
Il s'agit en d'autres termes d'une suite finie non ordonn\'ee
$M=(x_{1}^{\pm}, \dots, x_{n}^{\pm})$, \'eventuellement vide, de
points de $X$ d\'ecor\'es d'un signe $+$ ou $-$. L'ensemble des morphismes entre
deux tels objets $M$ et $M'$, est l'ensembles des classes de diff\'eomorphismes orient\'es
de vari\'et\'es compactes orient\'ees $B$ de dimension $1$, munies 
d'applications continues $B \longrightarrow |X|$ ainsi que d'identification 
$\partial B \longrightarrow |X|$ avec le morphisme $\bar{M}\coprod M \longrightarrow |X|$
(o\`u $\bar{M}$ est la vari\'et\'e $M$ munie de son orientation oppos\'ee). La structure 
mono\"\i dale sur $[\mathbf{Bord}_{1}^{X}]$ est induite naturellement par
la somme disjointe des vari\'et\'es orient\'ees. 
Dans le cas o\`u $X=S^{1}=B\mathbb{Z}$, le morphisme d'adjonction 
$B\mathbb{Z} \longrightarrow \mathbb{B}(S^{1})$ classifie 
un objet $x \in [\mathbf{Bord}_{1}^{S^{1}}]$ muni d'un automorphisme $u$. Cet objet 
est le point de base $x \in S^{1}$ (disons l'unit\'e de la structure de groupe sur $S^{1}$), 
positivement orient\'e, et $u$ est l'automorphisme correspondant au morphisme
$$u : [0,1] \longrightarrow |S^{1}|=|B\mathbb{Z}|$$
qui est un repr\'esentant du g\'en\'erateur de $\pi_{1}(|B\mathbb{Z}|,x)\simeq \mathbb{Z}$
(o\`u $[0,1]$ est muni de son orientation positive usuelle). Pour calculer la trace de
cette automorphisme $u : x^{+} \simeq x^{+}$, on commence par d\'ecrire 
le morphisme induit
$$v : x^{+} \coprod  x^{-} \rightarrow \emptyset$$
par la bijection
$$[x^{+}\coprod x^{-},y]\simeq [x^{+},x^{+}\coprod y].$$
Le morphisme $v$ est quand \`a lui toujours repr\'esent\'e
par le morphisme $u : [0,1] \longrightarrow |S^{1}|$, mais maintenant consid\'er\'e
comme un morphisme de $x^{+}\coprod x^{-}$ vers $\emptyset$. En pr\'ecomposant 
$v$ par le morphisme $\emptyset \rightarrow x^{+} \coprod x^{-}$, correspondant 
maintenant \`a l'identit\'e de $x^{+}$, on trouve le morphisme
$$Z:=Tr_{x^{+}}(u) : \emptyset \longrightarrow \emptyset.$$
Ce morphisme correspond \`a l'application continue 
$$[0,1]/0=1 \longrightarrow |S^{1}|$$
induite par $u$. 

Notons maintenant $Map_{Z}(\emptyset,\emptyset)$ la composante connexe
de $Map(\emptyset,\emptyset)$, pour $\emptyset \in \mathbf{Bord}_{1}^{S^{1}}$, 
contenant $Z$. On a $Map_{Z}(\emptyset,\emptyset) \subset Map_{1}(\emptyset,\emptyset)$, 
o\`u $Map_{1}(\emptyset,\emptyset)$ est la r\'eunion des composantes connexes
de $Map(\emptyset,\emptyset)$ form\'ee des vari\'et\'es $B \rightarrow X$, avec 
$B$ isomorphe \`a un cercle. 
On dispose d'une suite exacte de fibration
$$Map(S^{1},S^{1}) \longrightarrow Map_{1}(\emptyset,\emptyset) \longrightarrow BSO(2),$$
qui montre que $\pi_{0}(Map(S^{1},S^{1})) \longrightarrow 
\pi_{0}(Map_{1}(\emptyset,\emptyset))$ est bijectif. En particulier, la composante
connexe contenant $Z$ correspond \`a la composante connexe contenant l'identit\'e
de $S^{1}$. Ceci implique que l'on dispose d'une suite exacte de fibration 
$$Map_{id}(S^{1},S^{1}) \longrightarrow Map_{Z}(\emptyset,\emptyset) \longrightarrow BSO(2)\simeq BS^{1}.$$
Plus pr\'ecis\'ement, $Map_{Z}(\emptyset,\emptyset)$ est 
\'equivalent au quotient homotopique de $Map_{id}(S^{1},S^{1})$ par l'action 
naturelle de $S^{1}$. Comme $Map_{id}(S^{1},S^{1})$ est \'equivalent \`a $S^{1}$, 
en tant qu'ensemble simplicial $S^{1}$-\'equivariant, on trouve finalement que 
$Map_{Z}(\emptyset,\emptyset)$ est contractile. Ceci implique bien que 
$$\pi_{i}(End_{\mathbb{B}(X)}(\mathbf{1}),Z)=0 \qquad \forall \; i>0,$$
et termine donc la preuve du th\'eor\`eme. \hfill $\Box$ \\

La cons\'equence du th\'eor\`eme \ref{t2} que nous retiendrons par la suite est la suivante.
Afin de l'\'enoncer, nous noterons $S^{1}-\T$ la $\s$-cat\'egorie
des ensembles simpliciaux $S^{1}$-\'equivariant. Elle peut-\^etre d\'efinie
par exemple par la formule suivante
$$S^{1}-\mathbb{T}:=\rh (BS^{1},\T),$$
o\`u $BS^{1}$ d\'esigne ici le $\s$-groupo\"\i de
avec un unique objet ayant $S^{1}$ comme groupe simplicial d'endomorphismes. On a une \'equivalence
$$S^{1}-\mathbb{T} \simeq Int(\mathbf{SEns}^{S^{1}})\simeq L(\mathbf{SEns}^{S^{1}})\simeq \rh (BS^{1},\T),$$
o\`u $\mathbf{SEns}^{S^{1}}$ est la cat\'egorie de mod\`eles des ensembles simpliciaux munis d'une
action du groupe simplicial $S^{1}=B\mathbb{Z}$, et $BS^{1}$ est la $\s$-cat\'egorie
avec un unique objet et $S^{1}$ comme groupe simplicial d'endomorphismes. \\

\noindent\textbf{Traces cycliques.} On construit alors deux $\s$-foncteurs $\widetilde{\mathcal{LI}}$ et $\widetilde{End}_{-}(\mathbf{1})$. 
Le $\s$-foncteur $\widetilde{End}_{-}(\mathbf{1})$ est simplement le compos\'e
$$\xymatrix{
\uscms_{rig} \ar[r]^-{End_{-}(\mathbf{1})} & \T \ar[r]^{c}& S^{1}-\T},$$
o\`u $c$ est le $\s$-foncteur \emph{action triviale}. Le second $\s$-foncteur
$\widetilde{\mathcal{LI}}$ est le compos\'e
$$\xymatrix{
\uscms_{rig} \ar[r]^{ev_{[1]}} & \uscat \ar[r]^-{\mathcal{I}} & \T \ar[r]^-{\mathcal{L}} & S^{1}-\T},$$
o\`u $\mathcal{L}$ envoie un ensemble simplicial $X$ sur son espace de lacets $Map(S^{1},X)$, muni
de l'action induite par $S^{1}$.
Ces deux $\s$-foncteurs, sont des rel\`evements de $End_{-}(\mathbf{1})$
et $\mathcal{LI}$ en des foncteurs \`a valeurs dans $S^{1}-\T$. On dispose en particulier d'un 
morphisme d'oubli de l'action de $S^{1}$
$$Map(\widetilde{\mathcal{LI}},\widetilde{End}_{-}(\mathbf{1}))
 \longrightarrow Map(\mathcal{LI},End_{-}(\mathbf{1})),$$
dont la fibre homotopique en le morphisme $Tr$ de la proposition  \ref{p6} sera not\'ee
$$Map_{Tr}(\widetilde{\mathcal{LI}},\widetilde{End}_{-}(\mathbf{1})).$$

\begin{cor}\label{c4}
On a
$$Map_{Tr}(\widetilde{L\mathcal{I}},\widetilde{End}_{-}(\mathbf{1}))\simeq *.$$
\end{cor}

\textbf{Preuve --} En effet, $S^{1}$ op\`ere naturellement sur 
l'ensemble simplicial $Map(\mathcal{LI},End_{-}(\mathbf{1}))$, car il op\`ere
sur les $\s$-foncteurs $\mathcal{LI}$ et $End_{-}(\mathbf{1})$. De plus, 
le morphisme d'oubli
$$Map(\widetilde{\mathcal{LI}},\widetilde{End}_{-}(\mathbf{1}))
 \longrightarrow Map(\mathcal{LI},End_{-}(\mathbf{1})),$$
est alors isomorphe dans $\mathrm{Ho}(\mathbf{SEns})$ au morphisme naturel
$$Map(\mathcal{LI},End_{-}(\mathbf{1}))^{hS^{1}} \longrightarrow Map(\mathcal{LI},End_{-}(\mathbf{1})),$$
o\`u $Map(\mathcal{LI},End_{-}(\mathbf{1}))^{hS^{1}}$ d\'esigne les points fixes
homotopiques de $S^{1}$. Ainsi, 
$Map_{Tr}(\widetilde{\mathcal{LI}},\widetilde{End}_{-}(\mathbf{1}))$
s'identifie \`a la fibre homotopique du morphisme
$$Map_{Tr}(\mathcal{LI},End_{-}(\mathbf{1}))^{hS^{1}} \longrightarrow Map_{Tr}(\mathcal{LI},End_{-}(\mathbf{1})),$$
o\`u $Map_{Tr}$ d\'esigne la composante connexe contenant le morphisme $Tr$ de la
proposition \ref{p6}. 
Mais le th\'eor\`eme \ref{t2} implique que 
$Map_{Z}(\mathcal{LI},End_{-}(\mathbf{1}))$ est discret, et donc que le morphisme
$$Map_{Z}(\mathcal{LI},End_{-}(\mathbf{1}))^{hS^{1}} \longrightarrow Map_{Z}(\mathcal{LI},End_{-}(\mathbf{1}))$$
est automatiquement une \'equivalence. Ceci montre la contractibilit\'e 
de $Map_{Tr}(\widetilde{LI},\widetilde{End}_{-}(\mathbf{1}))$. 
\hfill $\Box$ \\

\begin{df}\label{d13}
Une transformation naturelle 
$$\widetilde{LI} \longrightarrow \widetilde{End}_{-}(\mathbf{1})$$
dans $\rh (\uscms_{rig},S^{1}-\T)$, qui rel\`eve la trace $Tr$ de la proposition \ref{p6}, 
 sera appel\'ee une \emph{trace cylique}. 
Elle sera not\'ee $Tr^{S^{1}}$. 
\end{df}

Notons que le corollaire \ref{c4} implique que l'espace des traces
cyliques est contractile, et donc que la trace cylique est unique (\`a isomorphismes
unique pr\`es, y compris sup\'erieurs). Nous parlerons donc 
de \emph{la} trace cylique. \\

\subsection{Additivit\'e des traces}\label{add} 

Dans ce paragraphe nous montrons que la trace cyclique $Tr^{S^{1}}$ de la d\'efinition
\ref{d13} est additive lorsque l'on se restreint aux $\s$-cat\'egories
mono\"\i dales sym\'etriques \emph{semi-additives}. Cette propri\'et\'e d'additivit\'e
n'es pas la plus g\'en\'erale, car l'on pourrait aussi d\'emontrer l'additivit\'e 
pour les triangles distingu\'es dans le cadre des $\s$-cat\'egories stables (voir \cite[\S 2]{lu4}). Cependant, 
cette forme plus restrictive, qui est plus facile \`a d\'emontrer, suffira pour 
les applications que nous avons en perspective. Elle sera utilis\'ee en particulier
pour comparer notre construction g\'en\'erale du caract\`ere de Chern avec le
caract\`ere de Chern usuel des fibr\'es vectoriels sur des vari\'et\'es alg\'ebriques lisses
(voir Appendice \ref{comparazione}). \\

Nous commencerons par les notions de $\s$-cat\'egories et de $\s$-CMS semi-additives. 

\begin{df}\label{dadd1}
\begin{enumerate}
\item 
Une $\s$-cat\'egorie $A$ est \emph{semi-additive} si elle v\'erifie les trois
conditions suivantes.
\begin{enumerate}
\item La $\s$-cat\'egorie $A$ poss\`ede des sommes et des produits finis. 
\item Le morphisme naturel
$$\emptyset \longrightarrow *,$$
de l'objet initial vers l'objet final, est une \'equivalence.
\item Pour toute paire d'objets $(x,y)$ dans $A$, le morphisme
$$x\coprod y \longrightarrow x\times y,$$
induit par les deux projections 
$$x\coprod y \longrightarrow x\coprod * \simeq x \qquad x\coprod y \longrightarrow *\coprod x \simeq y,$$
est une \'equivalence.
\end{enumerate}
\item Un $\s$-foncteur $f : A \longrightarrow B$ entre deux $\s$-cat\'egories est 
\emph{semi-additif} s'il commute aux limites finies. 
\end{enumerate}
\end{df}

Pour une $\s$-cat\'egorie $A$, l'$\s$-foncteur naturel $A \longrightarrow [A]$
est conservatif et commute aux sommes et produits finis. Ainsi, on voit qu'une 
$\s$-cat\'egorie $A$ est semi-additive si et seulement si sa cat\'egorie homotopique $[A]$ l'est. \\

Pour une $\s$-cat\'egorie $A$, poss\'edant des produits finis, nous pouvons introduire
la $\s$-cat\'egorie $CMon(A)$, des \emph{mono\"\i des commutatifs} dans $A$. 

\begin{df} Soit $\rh(\Gamma,A)$, l'$\s$-cat\'egorie
des $\s$-foncteurs de $\Gamma$ vers $A$. On d\'efinit $CMon(A)$ comme la sous-$\s$-cat\'egorie
pleine de $\rh(\Gamma,A)$ form\'ee des $\s$-foncteurs
$$H : \Gamma \longrightarrow A$$
tels que pour tout $n$ (y compris $n=0$) le morphisme naturel
$$H([n]) \longrightarrow H([1])^{n}$$
soit une \'equivalence dans $A$. \\ Les objets de $CMon(A)$ sont appel\'es les \emph{mono\"\i des 
commutatifs dans $A$}. 
\end{df}

On renvoie au $\S$ 2.1 pour des commentaires qui expliquent cette terminologie.\\ 

Nous pouvons maintenant caract\'eriser les $\s$-cat\'egories semi-additives comme
\'etant celles \'equivalentes \`a leurs $\s$-cat\'egories de mono\"\i des commutatifs. Cet \'enonc\'e
formalise l'id\'ee intuitive qu'une $\s$-cat\'egorie est semi-additive lorsque tous ses objets
poss\`edent une structure de mono\"\i des commutatifs naturelle. 

\begin{prop}\label{padd1}
Soit $A$ une $\s$-cat\'egorie poss\`edant des produits finis. 
Les deux assertions suivantes sont \'equivalentes.
\begin{enumerate}
\item L'$\s$-foncteur d'\'evaluation en $[1]$
$$CMon(A) \longrightarrow A$$
est une \'equivalence d'$\s$-cat\'egories.
\item L'$\s$-cat\'egorie $A$ est semi-additive.
\end{enumerate}
\end{prop}

\textbf{Preuve --} Pour montrer que $(1)$ implique $(2)$, il suffit de montrer que 
$CMon(A)$ est une $\s$-cat\'egorie semi-additive au sens de la d\'efinition \ref{dadd1}.
On consid\`ere $CMon(A)$ comme naturellement plong\'ee dans $\rh(\Gamma,A)$. C'est une sous-$\s$-cat\'egorie
pleine stable par limites, et donc stables par produits finis. On voit ainsi que 
$CMon(A)$ poss\`ede des produits finis, et que ceux-ci se calculent dans 
$\rh(\Gamma,A)$. L'objet final de $CMon(A)$ est donc le $\s$-foncteur ponctuel 
$* : \Gamma \longrightarrow A$, qui n'est autre que le $\s$-foncteur corepr\'esent\'e
par $[0] \in \Gamma$. Cet objet final est aussi initial, car on a, d'apr\`es le lemme
de Yoneda, pour tout mono\"\i de commutatif $M : \Gamma \longrightarrow A$
$$Map(*,M)\simeq M([0]) \simeq *.$$
Soient $M$ et $N$ deux mono\"\i des commutatifs dans $A$, et montrons que les morphismes
naturels 
$$u : \xymatrix{M \ar[r]^-{id\times *} & M\times N} \qquad
v : \xymatrix{N \ar[r]^-{*\times id} & M\times N}$$
font de $M\times N$ la somme de $M$ et de $N$ dans $CMon(A)$. Pour cela, il suffit de montrer que
pour tout $P \in [CMon(A)]$, l'application induite
$$\phi : \xymatrix{[M\times N,P] \ar[r]^-{u^{*}\times v^{*}} & [M,P] \times [N,P]}$$
est bijective. Pour cela, on construit une application inverse $\psi$ de la fa\c{c}on suivante.
Nous d\'efinissons un nouveau mono\"\i de commutatif $P\wedge P : \Gamma \longrightarrow A$, 
dont la valeur sur $[n] \in \Gamma$ est par d\'efinition $P([n+n])$, et sur un 
morphisme $a : [n] \rightarrow [m]$ est $P(a+a) : P([n+n]) \longrightarrow P([m+m])$
(on remarquera que $[n+n]$ est la somme de $[n]$ avec $[n]$ dans la cat\'egorie $\Gamma$). 
L'objet $P\wedge P$ vient avec une \'equivalence
$$P \wedge P \longrightarrow P\times P,$$
qui est induite par les morphismes
$$P([n+n]) \longrightarrow P([n])\times P([n])$$
provenant des deux morphismes de projection $[n+n] \longrightarrow [n]$. D'autre part, nous avons
un morphisme de mono\"\i des commutatifs
$$ + : P\wedge P \longrightarrow P$$
qui est induit par le morphisme
$$P([n+n]) \longrightarrow P([n])$$
provenant du morphisme $[n+n] \longrightarrow [n]$ \'egal \`a l'identit\'e sur chaque
copie de $[n]$ dans $[n+n]$. Ainsi, dans $[CMon(A)]$, on dispose d'un morphisme naturel
$$+ : P\times P \simeq P\wedge P \longrightarrow P.$$
On d\'efinit alors
$$\psi : [M,P] \times [N,P] \longrightarrow [M\times N,P]$$
qui \`a deux morphismes $f : M \rightarrow P$ et $g : N \rightarrow P$ associe
$$f+g:=(+)\circ (f\times g) : [M,P] \times [N,P] \rightarrow
[M\times N,P\times P] \rightarrow [M\times N,P].$$
On laisse le soin au lecteur de v\'erifier que $\psi$ ainsi construit est l'inverse
de $\phi$. \\

R\'eciproquement, supposons que $A$ soit semi-additive. Consid\'erons 
$CMon(A)$ comme plong\'ee dans $\rh_{*}(\Gamma,A)$, o\`u $\rh_{*}(\Gamma,A)$ est la 
sous-$\s$-cat\'egorie pleine de $\rh(\Gamma,A)$ form\'ee des $\s$-foncteurs
qui envoient $[0]$ sur $* \in A$. Le $\s$-foncteur 
d'\'evaluation en $[1]$
$$ev_{1} : \rh_{*}(\Gamma,A) \longrightarrow A,$$
poss\`ede un adjoint \`a gauche
$$L : A \longrightarrow \rh(\Gamma,A)$$
qui \`a $a \in A$ associe le $\s$-foncteur 
$$L(a) : [n] \mapsto \coprod_{\Gamma([1],[n])-*} \wedge a,$$
o\`u $* : [1] \rightarrow [n]$ est l'application constante.
Comme $A$ est semi-additive, on v\'erifie que $L(a)$ est
un objet de $CMon(A)$. Ainsi, le $\s$-foncteur d'\'evaluation en $1$
$$CMon(A) \longrightarrow A$$
poss\`ede $L$ comme adjoint \`a gauche. L'unit\'e de cette adjonction, pour $a\in A$, 
est le morphisme 
$$\coprod_{\Gamma([1],[1])-*}a \longrightarrow a$$
qui est une \'equivalence dans $A$ car $\Gamma([1],[1])-*\simeq *$. 
\hfill $\Box$ \\

Soit $\uscat^{ad}$ la sous-$\s$-cat\'egorie, non pleine, de $\uscat$ form\'e des 
$\s$-cat\'egories semi-additives et des $\s$-foncteurs qui commutent aux produits finis. 
On dispose d'un 
plongement naturel
$$j : \uscat^{ad} \longrightarrow \uscat.$$

\begin{prop}\label{padd2}
\begin{enumerate}
\item 
Le $\s$-foncteur $j$ ci-dessus poss\`ede un adjoint \`a gauche 
$$(-)^{ad} : \uscat \longrightarrow \uscat_{ad}.$$
\item Pour tout $A\in \uscat$, le morphisme d'ajonction
$$f : A \longrightarrow A^{ad}$$
est tel que pour tout $a$ et tout $b$ objet de $A$, on ait
$$Map(f(a),f(b)) \simeq \coprod_{n \geq 0}(Map(a,b)^{n})_{h\Sigma_{n}}.$$
En d'autres termes, $Map(f(a),f(b))$ est le mono\"\i de commutatif libre
dans $\T$ engendr\'e par $Map(a,b)$. 
\end{enumerate}
\end{prop}

\textbf{Preuve --} Nous allons d\'emontrer simultan\'ement $(1)$ et $(2)$ en construisant 
explicitement l'adjoint \`a gauche en question. Soit $A \in \uscat$ une $\s$-cat\'egorie, 
et $\widehat{A}$ l'$\s$-cat\'egorie des pr\'efaisceaux. On consid\`ere le plongement 
de Yoneda
$$h : A \hookrightarrow \widehat{A}.$$
Soit $CMon(\widehat{A})$ l'$\s$-cat\'egorie des mono\"\i des commutatifs dans
$\widehat{A}$, et $CMon(\widehat{A}) \longrightarrow \widehat{A}$ 
le foncteur d'\'evaluation en $[1]$. L'adjoint \`a gauche de ce foncteur
est not\'e
$$L : \widehat{A} \longrightarrow CMon(\widehat{A}).$$
On consid\`ere le foncteur compos\'e
$$\phi:=L\circ h : A \longrightarrow \widehat{A} \longrightarrow CMon(\widehat{A}).$$
On d\'efinit $A^{ad}$ comme \'etant la sous-$\s$-cat\'egorie pleine de 
$CMon(\widehat{A})$ des objets qui sont des produits finis d'objets de l'image essentielle
du $\s$-foncteur $\phi$. L'$\s$-foncteur $\phi$ nous donne un $\s$-foncteur
$$j : A \longrightarrow A^{ad}.$$ 
Il nous faut montrer que $A^{ad}$ est semi-additive et que le $\s$-foncteur 
$j$ induit, pour tout $B\in \uscat_{ad}$, une bijection
$$j^{*} : [A^{ad},B]^{ad} \simeq [A,B],$$
o\`u le membre de droite (resp. de gauche) d\'esigne l'ensemble des
morphismes dans la cat\'egorie $[\uscat_{ad}]$ (resp. $[\uscat]$). 

Pour montrer cela nous d\'efinissons une application inverse
$$p : [A,B] \longrightarrow [A^{ad},B]^{ad}$$
de la fa\c{c}on suivante. Soit $f : A \longrightarrow B$ un $\s$-foncteur, qui induit 
un $\s$-foncteur
$$f^{*} : CMon(\widehat{B})\simeq \rh (B^{op},CMon) \longrightarrow \rh (A^{op},CMon) 
\simeq CMon(\widehat{A}),$$
o\`u l'on note $CMon:=CMon(\T)$. L'adjoint \`a gauche de $f^{*}$ est not\'e
$$f_{!} : CMon(\widehat{A}) \longrightarrow CMon(\widehat{B}).$$
Le $\s$-foncteur $f_{!}$ est caract\'eris\'e par, d'une part la propri\'et\'e de commuter aux colimites
et d'autre part la propri\'et\'e de rendre le diagramme suivant commutatif
$$\xymatrix{
CMon(\widehat{A}) \ar[r]^-{f_{!}} & CMon(\widehat{B}) \\
A \ar[u]^-{\phi} \ar[r]_-{f} & B. \ar[u]_-{\phi} }$$
Les $\s$-cat\'egories $CMon(\widehat{A})$ et $CMon(\widehat{B})$ sont semi-additives, 
et le $\s$-foncteur $f_{!}$ commute aux sommes finies et est donc semi-additif. Il envoie donc
la sous-$\s$-cat\'egorie $A^{ad} \subset CMon(\widehat{A})$ dans 
$B^{ad} \subset CMon(\widehat{B})$, et d\'efinit ainsi un $\s$-foncteur semi-additif
$$f^{ad} : A^{ad} \longrightarrow B^{ad}.$$
Cependant, $B$ \'etant elle-m\^eme semi-additive, le $\s$-foncteur $\phi : B \longrightarrow 
CMon(\widehat{B})$ est une \'equivalence (proposition \ref{padd1}), et ainsi le $\s$-foncteur
$j : B \longrightarrow B^{ad}$ est une \'equivalence. On trouve ainsi un $\s$-foncteur semi-additif
$$p(f) : A^{ad} \longrightarrow B.$$
Par construction, il est facile de v\'erifier que $f \mapsto p(f)$ d\'efinit un inverse 
de l'application $j^{*} : [A^{ad},B]^{ad} \rightarrow [A,B]$. 

Ceci termine l'existence de l'adjoint $A \mapsto A^{ad}$, avec de plus une description 
explicite de $A^{ad}$ comme \'etant la sous-$\s$-cat\'egorie pleine de 
$\rh (A ^{op},CMon)$ form\'ee des produits finis de $\s$-foncteurs de la forme $L(h_{a})$ pour $a\in A$. 
Ainsi, pour d\'emontrer la seconde assertion de la proposition il nous suffit de 
calculer, pour deux objets $a$ et $b$ dans $A$, l'espace des morphismes
$$Map_{A^{ad}}(j(a),j(b))\simeq Map_{\rh (A^{op},CMon)}(L(h_{a}),L(h_{b})).$$
Mais on a
$$Map_{\rh (A^{op},CMon)}(L(h_{a}),L(h_{b})) \simeq 
Map_{\widehat{A}}(h_{a},L(h_{b}))\simeq L(h_{b})(a).$$
Par d\'efinition, $L(h_{b})(a)$ est le mono\"\i de commutatif dans $\T$ libre sur 
l'objet $Map_{A}(a,b)$. Ce mono\"\i de commutatif libre est 
d\'ecrit dans \cite[Prop. 3.6]{seg}, 
et est naturellement \'equivalent \`a
$$\coprod_{n\in \mathbb{N}} (Map_{A}(a,b)^{n})_{h\Sigma_{n}}.$$
\hfill $\Box$ \\

Nous allons maintenant nous int\'er\'esser au cas des $\s$-cat\'egories mono\"\i dales
sym\'etriques et de leur semi-additivisation. 

\begin{df}\label{dadd2}
\begin{enumerate}
\item Une \emph{$\s$-cat\'egorie mono\"\i dale sym\'etrique semi-additive} est une
$\s$-cat\'egorie mono\"\i dale sym\'etrique $T$ v\'erifiant les deux propri\'et\'es suivantes:
\begin{enumerate}
\item La $\s$-cat\'egorie sous-jacente \`a $T$ est semi-additive. 

\item Pour tout objet $x\in T$, le $\s$-foncteur
$$x\otimes - : T \longrightarrow T$$
est semi-additif. 
\end{enumerate}
\item Un $\s$-foncteur mono\"\i dal sym\'etrique entre $\s$-cat\'egories mono\"\i dales sym\'etriques 
semi-additives est \emph{semi-additif} s'il l'est en tant que $\s$-foncteur entre $\s$-cat\'egories
semi-additives. 
\end{enumerate}
\end{df}

Nous noterons $\uscms_{ad}$ la sous-$\s$-cat\'egorie (non pleine) de $\uscms$ form\'ee
des $\s$-cat\'egories mono\"\i dales sym\'etriques 
semi-additives et des $\s$-foncteurs mono\"\i daux sym\'etriques semi-additifs. On dispose de $\s$-foncteurs
d'oubli
$$\uscms_{ad} \longrightarrow \uscat_{ad} \qquad \uscms \longrightarrow \uscat,$$
qui commutent aux inclusions $\uscms_{add} \subset \uscms$ et 
$\uscat_{ad} \subset \uscat$. 

\begin{prop}\label{padd3}
\begin{enumerate}
\item Le $\s$-foncteur d'inclusion $\uscms_{ad} \hookrightarrow \uscms$ poss\`ede un 
adjoint \`a gauche
$$(-)^{\otimes-ad} : \uscms \longrightarrow \uscms_{ad}.$$
\item Le diagramme suivant
$$\xymatrix{
\uscms \ar[r]^-{(-)^{\otimes-ad}} \ar[d] & \uscms_{ad} \ar[d] \\
\uscat_{ad} \ar[r] & \uscat,}$$
commute \`a \'equivalence pr\`es. 
\end{enumerate}
\end{prop}

\textbf{Preuve --} Nous allons d\'ecrire explicitement l'adjoint \`a gauche $(-)^{\otimes-ad}$. Pour cela,
soit $T : \Gamma \longrightarrow \s-Cat$ une $\s$-cat\'egorie mono\"\i dale sym\'etrique. Notons
$$p : \mathcal{T}:=\int_{\Gamma}T \longrightarrow \Gamma$$
l'$\s$-cat\'egorie cofibr\'ee sur $\Gamma$ correspondante. Pour simplifier les notations nous supposerons
que $T$ est \`a valeur dans les $\s$-cat\'egories strictes.

Nous commen\c{c}ons par d\'efinir une $\s$-cat\'egorie $\mathcal{C}$, bi-fibr\'ee au-dessus
de $\Gamma$, de la fa\c{c}on suivante. Les objets de $\mathcal{C}$ sont des couples
$(I,x)$ o\`u:

\begin{itemize}

\item $I$ est un objet de $\Gamma$

\item $x$ est un objet de l'$\s$-cat\'egorie
$\rh(T([1])^{op},CMon)^{I^{0}}$, o\`u $CMon$ est l'$\s$-cat\'egorie
des mono\"\i des commutatifs dans $\mathbb{T}$, et $I^{0}$ est l'ensemble
$I$ priv\'e de son point de base. L'objet $x$ sera repr\'esent\'e par 
une famille de $\s$-foncteur $x_{i} : T([1])^{op} \rightarrow CMon$ pour $i\in I^{0}$. 

\end{itemize}

Ci-dessus, nous prendrons un mod\`ele strict bien pr\'ecis pour $\rh(T([1])^{op},CMon)$.
Nous noterons $\Gamma-SEns_{sp}^{f}$ la cat\'egorie des $\Gamma$-espaces
fibrants et sp\'eciaux (voir Prop. \ref{p2}, dans laquelle il faut remplacer la cat\'egorie
de mod\`eles des cat\'egories de Segal par $SEns$). La cat\'egorie
$\Gamma-SEns_{sp}^{f}$ est telle que sa localis\'ee de Dwyer-Kan v\'erifie
$$L(\Gamma-SEns_{sp}^{f}) \simeq CMon.$$
Ainsi, un mod\`ele strict explicite pour $\rh(T([1])^{op},CMon)$ est 
$$L((\Gamma-SEns_{sp}^{f})^{T[1]^{op}})\simeq \rh(T([1])^{op},CMon).$$
Les objets de $\rh(T([1])^{op},CMon)$ seront donc identifi\'es \`a des
$\s$-foncteurs stricts
$$T([1])^{op} \longrightarrow \Gamma-SEns_{sp}^{f}.$$

Les espaces de morphismes dans $\mathcal{C}$ sont d\'efinis de la fa\c{c}on suivante. Soient
$(I,x)$ et $(J,y)$ deux objets, on pose
$$\mathcal{C}((I,x),(J,y)):=
\coprod_{u : I \rightarrow J}Map_{\rh(T([I])^{op},CMon)}(s_{I}^{*}(x),u^{*}s_{J}^{*}(y)).$$
Dans cette expression, $u$ parcourt l'ensemble des morphismes $I \rightarrow J$ dans $\Gamma$, 
et $$u^{*} : \rh(T([J])^{op},CMon) \longrightarrow \rh(T([I])^{op},CMon)$$ est le morphisme
induit. L'objet $s_{I}^{*}(x)$, qui est un $\s$-foncteur
$$T([I])^{op} \longrightarrow CMon,$$
est le compos\'e
$$\xymatrix{T([I])^{op} \ar[r] &  T([1])^{I^{0}} \ar[r]^-{\boxtimes x} & CMon,}$$
o\`u $\boxtimes x$ est le $\s$-foncteur produit des $x_{i}$. La composition des
morphismes dans $\mathcal{C}$ se fait de mani\`ere naturelle. 

On dispose d'une projection $\pi : \mathcal{C} \longrightarrow \Gamma$, et on 
v\'erifie facilement que cela fait de $\mathcal{C}$ une $\s$-cat\'egorie (stricte) bi-fibr\'ee
sur $\Gamma$. Il existe une sous-$\s$-cat\'egorie pleine $\mathcal{C}_{0}$ de 
$\mathcal{C}$ qui consiste \`a se restreindre aux objets de la forme
$(I,x)$, avec chaque $x_{i} \in \rh(T([1])^{op},CMon)$ qui est 
\'equivalent \`a un produit fini de repr\'esentables. On v\'erifie que $\mathcal{C}_{0}$ 
est encore cofibr\'ee sur $\Gamma$ (mais elle n'est plus fibr\'ee), et 
d\'efinit ainsi par la proposition \ref{pcart} un foncteur
$$\Gamma \longrightarrow \s-Cat.$$
Par construction, la valeur de se foncteur en $I$ est naturellement \'equivalente \`a 
$(T([1])^{ad})^{I^{0}}$, et on voit sans peine qu'il d\'efinit une 
$\s$-cat\'egorie mono\"\i dale sym\'etrique semi-additive que l'on notera
$T^{\otimes-ad}$. La $\s$-cat\'egorie sous-jacente est \'equivalente \`a 
$T([1])^{ad}$. Pour d\'emontrer la proposition, il nous reste donc \`a v\'erifier 
que $T \mapsto T^{\otimes-ad}$ est bien adjoint \`a gauche de l'inclusion
$\uscms_{ad} \hookrightarrow \uscms$, ce que nous laissons aux lecteurs. 
\hfill $\Box$ \\

Nous sommes maintenant en mesure de d\'emontrer le caract\`ere additif des 
traces cycliques. Pour cela, nous introduisons $\uscms_{ad,rig}$, la sous-$\s$-cat\'egorie
pleine form\'ee des objets dont la $\s$-cat\'egorie mono\"\i dale sym\'etrique est 
sous-jacente est rigide. En pr\'ecomposant par le $\s$-foncteur d'oubli de la semi-additivit\'e
$$\uscms_{ad,rig} \longrightarrow \uscms_{rig},$$
les $\s$-foncteurs $\widetilde{\mathcal{LI}}$ et $\widetilde{End}_{-}(\mathbf{1})$ du paragraphe 
\S 2.3 induisent deux nouveaux $\s$-foncteurs
$$\widetilde{\mathcal{LI}}_{ad}, \widetilde{End}_{-}(\mathbf{1})_{ad} : 
\uscms_{ad,rig} \longrightarrow  S^{1}-\mathbb{T}.$$
De m\^eme, en oubliant l'action du cercle nous avons deux $\s$-foncteurs
$$\mathcal{LI}_{ad}, End_{-}(\mathbf{1})_{ad} : 
\uscms_{ad,rig} \longrightarrow  \mathbb{T}.$$
La proposition \ref{padd3} nous dit que le $\s$-foncteur $\mathcal{LI}_{ad}$
est correpr\'esent\'e par $\mathbb{B}(S^{1})^{\otimes-ad}$, dont la $\s$-cat\'egorie
semi-additive sous-jacente est $\mathbb{B}(S^{1})^{ad}$. De m\^eme, le $\s$-foncteur
$\mathcal{LI}_{ad}\times \mathcal{LI}_{ad}$ est corepr\'esent\'e par 
$\mathbb{B}(S^{1}\coprod S^{1})^{\otimes-ad}$. 
Ainsi, l'espace des transformations naturelles
de $\mathcal{LI}_{ad}\times \mathcal{LI}_{ad}$ vers 
$End_{-}(\mathbf{1})_{ad}$ est donn\'e par
$$Map(\mathcal{LI}_{ad}^{2},End_{-}(\mathbf{1})_{ad})\simeq
End_{\mathbb{B}(S^{1}\coprod S^{1})^{ad}}(\mathbf{1}).$$
D'apr\`es la proposition \ref{padd2}, on trouve
que $Map(\mathcal{LI}_{ad}^{2},End_{-}(\mathbf{1})_{ad})$ est \'equivalent au mono\"\i de commutatif
libre sur l'espace $End_{\mathbb{B}(S^{1}\coprod S^{1})}(\mathbf{1})$. L'action de $S^{1}$ sur
l'espace $Map(\mathcal{LI}^{2}_{ad},End_{-}(\mathbf{1})_{ad})$ est quand \`a elle 
induite par l'action naturelle sur $End_{\mathbb{B}(S^{1}\coprod S^{1})}(\mathbf{1})$. 

Nous allons maintenant consid\'erer deux transformations naturelles $S^{1}$-\'equivariantes
$$f,g : \mathcal{LI}_{ad}^{2} \Rightarrow End_{-}(\mathbf{1})_{ad},$$
et montrer qu'elles sont \'equivalentes, et ce de mani\`ere $S^{1}$-\'equivariante. La premi\`ere, 
$f$ est le compos\'e
$$\xymatrix{
\mathcal{LI}_{ad}^{2} \ar[r]^-{\times} & \mathcal{LI}_{ad} \ar[r]^-{Tr} & 
End_{-}(\mathbf{1})_{ad}},$$
o\`u $Tr$ est la transformation naturelle de la proposition \ref{p6} restreinte
\`a la sous-$\s$-cat\'egorie $\uscms_{ad,rig}$. La transformation naturelle
$\times$ est quand \`a elle induite par le produit direct dans la $\s$-cat\'egorie
sous-jacente, qui existe par hypoth\`ese de semi-additivit\'e. En utilisant 
le corepr\'esentabilit\'e de $\mathcal{LI}_{ad}^{2}$, la transformation 
naturelle $\times$ correspond \`a un couple $(x,\beta)$, form\'e d'un objet $x$ muni d'une auto-\'equivalence
$\beta$ dans $\mathbb{B}(S^{1}\coprod S^{1})_{ad}$. En identifiant 
$[\mathbb{B}(S^{1}\coprod S^{1})_{ad}]$ avec la cat\'egorie semi-additive engendr\'ee
par $[\mathbf{Bord}_{1}^{S^{1}\coprod S^{1}}]$, la cat\'egorie des $1$-bordismes au-dessus
de $S^{1}\coprod S^{1}$ (voir la preuve du th\'eor\`eme \ref{t2}), 
le couple $(x,\beta)$ s'identifie \`a la somme directe
$(x_{1}^{+},u_{1})$, $(x_{2}^{+},u_{2})$, o\`u $x^{+}_{i}$ est un point, orient\'e 
positivement, au-dessus de la $i$-\`eme copie de $S^{1}$, et $u_{i}$ est l'automorphisme
canonique de $x_{i}^{+}$ qui couvre la $i$-\`eme copie de $S^{1}$
$$(x,u)=(x_{1}^{+},u_{1}) \oplus (x_{2}^{+},u_{2}).$$
Ainsi, la transformation naturelle $f$ est repr\'esent\'e par un 
\'el\'ement $Z(f)$ de 
$$End_{\mathbb{B}(S^{1}\coprod S^{1})_{ad}}(\mathbf{1}),$$
qui est la trace de $(x_{1}^{+},u_{1}) \oplus (x_{2}^{+},u_{2})$. Notons
$X:=End_{\mathbb{B}(S^{1}\coprod S^{1})}(\mathbf{1})$, de sorte que d'apr\`es la proposition 
\ref{padd2} on ait
$$End_{\mathbb{B}(S^{1}\coprod S^{1})_{ad}}(\mathbf{1}) \simeq \coprod_{n\geq 0}
X^{n}_{h\Sigma_{n}}.$$
Dans $X$, nous disposons des \'el\'ements $Z_{1}$ et $Z_{2}$, qui sont les
endomorphismes de de l'unit\'e $\emptyset \rightarrow S^{1}\coprod S^{1}$
repr\'esent\'es par les inclusions canoniques $S^{1}$ dans l'un des deux
facteurs. L'\'el\'ement $Z(f)$ est alors le couple
$$(Z_{1},Z_{2}) \in X^{2}_{h\Sigma_{2}} \subset End_{\mathbb{B}(S^{1}\coprod S^{1})_{ad}}(\mathbf{1}).$$
Cet \'el\'ement est de plus naturellement $S^{1}$-\'equivariant, car $Z_{1}$ et $Z_{2}$
le sont en tant qu'endomorphismes dans $\mathbb{B}(S^{1}\coprod S^{1})$. On peut-\^etre plus
pr\'ecis ici, car l'action de $S^{1}$ sur la composante connexe de $X^{2}_{h\Sigma_{2}}$
qui contient $(Z_{1},Z_{2})$ est en r\'ealit\'e triviale d'apr\`es le th\'eor\`eme \ref{t2}, et 
la structure de $S^{1}$-point fixe est alors la structure canonique. 

Passons maintenant \`a la transformation naturelle $g$. Elle est d\'efinie comme le compos\'e
$$\xymatrix{
\mathcal{LI}_{ad}^{2} \ar[r]^-{Tr\times Tr} & End_{-}(\mathbf{1})_{ad}^{2} \ar[r]^-{\times} & 
End_{-}(\mathbf{1})_{ad}},$$
o\`u $\times$ est encore induite par le caract\`ere semi-additif qui assure que 
l'objet $\mathbb{1}$ est canoniquement un mono\"\i de commutatif, ce qui induit
une structure de mono\"\i de commutatif naturelle sur $End_{-}(\mathbf{1})_{ad}$.
Si l'on utilise la corepr\'esentabilit\'e de $\mathcal{LI}_{ad}^{2}$ par 
$\mathbb{B}(S^{1}\coprod S^{1})_{ad}$, on trouve que la transformation naturelle 
$g$ est elle aussi donn\'ee par le couple 
$$(Z_{1},Z_{2}) \in X^{2}_{h\Sigma_{2}} \subset End_{\mathbb{B}(S^{1}\coprod S^{1})_{ad}}(\mathbf{1}),$$
qui est aussi $S^{1}$-\'equivariant car l'action de $S^{1}$ est triviale sur la composante
contenant ce couple. 

La conclusion de cette discussion est le corollaire suivant, qui est la forme de l'additivit\'e
des traces cycliques dont nous aurons besoin par la suite (mais qui n'est pas la forme la plus
g\'en\'erale qui existe). 

\begin{cor}\label{tracad}
Les deux transformations naturelles d\'efinies pr\'ec\'edemment
$$f,g : \mathcal{LI}_{ad}^{2} \Rightarrow End_{-}(\mathbf{1})_{ad}$$
sont \'equivalentes en tant que transformations naturelles $S^{1}$-\'equivariantes. 
\end{cor}

En particulier, on trouve le corollaire suivant. 

\begin{cor}\label{tracad'}
Soit $T$ une $\s$-cat\'egorie mono\"\i dale sym\'etrique semi-additive. Alors le morphisme
de trace cyclique
$$Tr^{S^{1}} : \pi_{0}(\mathcal{LI}(T)) \longrightarrow \pi_{0}(End_{T}(1)^{hS^{1}})$$
est un morphisme de mono\"\i des commutatifs, pour les lois de mono\"\i des
induites par le caract\`ere semi-additif de $T$. 
\end{cor}

\subsection{Muliplicativit\'e des traces}\label{mult} 

Terminons par la propri\'et\'e de multiplicativit\'e de la trace et de la trace
cylique.  Pour cela nous allons commencer par construire un 
$\s$-foncteur
$$\mathcal{M} : \uscms \longrightarrow \rh (\Gamma,\uscms).$$
On commence par construire un foncteur au niveau des cat\'egories de mod\`eles
$$\scps \longrightarrow (\scps)^{\Gamma},$$
qui envoie un objet $A$ sur le foncteur
$$\mathcal{M}(A) : \Gamma \longrightarrow \scps$$
d\'efini par 
$$\begin{array}{cccc}
\mathcal{M}(A)([n]) : & \Gamma & \longrightarrow &  \s -\mathbf{Cat}^{pr} \\
 & [m] & \mapsto & A([nm]).
 \end{array}$$
Notons que l'on utilise ici le fait que $\Gamma$ soit naturellement munie d'une
structure de cat\'egorie mono\"\i dale sym\'etrique dont la loi mono\"\i dale
est d\'efini par $[n]\otimes [m]=[nm]$. Les coh\'erences d'associativit\'e
et d'unit\'e sont alors des identit\'es. Les coh\'erences de commutativit\'e
sont elles induites par l'automorphisme $\alpha_{n,m} : [nm] \simeq [nm]$ 
qui est un ''shuffle'' de type $(n,m)$. On v\'erifie alors que 
si $A$ est une $\s$-CMS il en est de m\^eme de $\mathcal{A}([n])$ pour tout
$[n] \in \Gamma$. On a ainsi construit un foncteur
$$\scms \longrightarrow (\scms)^{\Gamma}$$
qui par localisation induit un $\s$-foncteur
$$\mathcal{M} : \uscms \longrightarrow \rh (\Gamma,\uscms).$$

En composant ce foncteur avec $\widetilde{\mathcal{LI}}$ et 
$\widetilde{End}_{-}(\mathbf{1})$ en trouve deux $\s$-foncteurs
$$\widetilde{\mathcal{LI}}\circ \mathcal{M}, \widetilde{End}_{-}(\mathbf{1}) \circ \mathcal{M} :
\uscms \longrightarrow \rh (\Gamma,S^{1}-\T).$$
Comme $\widetilde{\mathcal{LI}}$ et 
$\widetilde{End}_{-}(\mathbf{1})$ commutent avec les produits finis il est facile de voir que 
ces les deux foncteurs $\widetilde{\mathcal{LI}}\circ  \mathcal{M}$ et 
$\widetilde{End}_{-}(\mathbf{1})\circ \mathcal{M}$
se factorisent \`a travers la sous-$\s$-cat\'egorie 
pleine de $\rh (\Gamma,S^{1}-\T)$ form\'ee des $\s$-foncteurs $F$ v\'erifiant que pour tout
$n$ le morphisme 
$$F([n]) \longrightarrow F([1])^{n}$$
est une \'equivalence dans $S^{1}-\T$. En d'autres termes, ces deux $\s$-foncteurs 
sont \`a valeurs dans la $\s$-cat\'egorie des \emph{mono\"\i des commutatifs $S^{1}$-\'equivariants
dans $\T$}. Il suffit alors de pr\'ecomposer la transformation naturelle
$Tr^{S^{1}}$ avec le $\s$-foncteur $\mathcal{M}$ pour trouver une transformation 
naturelle
$$Tr_{\otimes}^{S^{1}} : \widetilde{\mathcal{LI}}\circ \mathcal{M}  \Rightarrow 
\widetilde{End}_{-}(\mathbf{1})\circ \mathcal{M}.$$
On note qu'\'evalu\'ee en $[1] \in \Gamma$ cette transformation naturelle 
n'est autre que $Tr^{S^{1}}$. On peut v\'erifier (nous  ne le ferons pas) que 
$Tr_{\otimes}^{S^{1}}$ est l'unique
extension de $Tr^{S^{1}}$ en une transformation naturelle entre
$\widetilde{\mathcal{LI}}\circ \mathcal{M}$ et $\widetilde{End}_{-}(\mathbf{1})\circ \mathcal{M}$. 

\begin{df}\label{d14}
La \emph{trace cyclique multiplicative} est la transformation naturelle
$$Tr_{\otimes}^{S^{1}} : \widetilde{\mathcal{LI}}\circ \mathcal{M} \longrightarrow
 \widetilde{End}_{-}(\mathbf{1})\circ \mathcal{M}$$
 d\'efinie ci-dessus.
\end{df}

\section{$\s$-Topos mono\"\i del\'es et catannel\'es}

Le but de cette section est de construire le caract\`ere de Chern
\`a proprement parl\'e. Pour cela nous commecerons par 
d\'ecrire la notion de \emph{$\s$-topos catannel\'e}, qui sont 
aux $\s$-cat\'egories mono\"\i dales sym\'etriques ce que sont
les topos annel\'es aux anneaux. En d'autres termes, un 
$\s$-topos catannel\'e est la donn\'e d'un couple
$(T,\mathcal{A})$, form\'e d'un $\s$-topos $T$ et d'un
champ en $\s$-cat\'egories mono\"\i dales sym\'etriques $\mathcal{A}$ sur $T$.
Un tel couple sera dit \emph{rigide} si les valeurs de
$\mathcal{A}$ sont rigides. Nous d\'efinissons pour tout
$\s$-topos catannel\'e $(T,\mathcal{A})$ une notion 
d'homologie de Hochschild et d'homologie cyclique relativement \`a
$\mathcal{A}$. Lorsque $(T,\mathcal{A})$ est rigide nous construisons 
alors quatres transformations naturelles: une pr\'e-trace, une trace, 
un pr\'e-caract\`ere de Chern et un caract\`ere de Chern. Les deux morphismes
de traces sont construits en utilisant l'existence de traces pr\'edit
par la proposition \ref{p6} et sont \`a valeurs dans l'homologie
de Hochschild: il s'agit d'une version extr\`emement g\'en\'erale
du morphisme Trace de Dennis utilis\'e en K-th\'eorie (voir \cite[8.4]{lo}). Les deux
morphismes de caract\`ere de Chern sont construits \`a l'aide 
de l'existence de traces cycliques de \ref{d13} et sont 
\`a valeurs dans l'homologie cyclique: il s'agit 
de g\'en\'eralisation du caract\`ere de Chern \`a valeurs dans
l'homologie cyclique n\'egative utilis\'e en K-th\'eorie (\cite[11.4]{lo}).

\subsection{La $\s$-cat\'egorie des $\s$-Topos structur\'es}

Soit $A \in \uscat$ une $\s$-cat\'egorie fibrante. Nous allons
commencer par d\'efinir la $\s$-cat\'egorie $\uscat_{//A}$ des $\s$-cat\'egories fibrantes au-dessus 
de $A$. Pour cela notons $\underline{Hom}^{\Delta}$ les 
Hom simpliciaux de la cat\'egorie de mod\`eles $\s -\mathbf{Cat}^{pr}$ (voir 
\cite{hisi}), et 
$\Delta^{n}$ la cat\'egorie classifiant les chaines de $n$ morphismes composables. 
Rappelons que la structure simpliciale sur $\s -\mathbf{Cat}^{pr}$ n'est pas donn\'ee
par $n \mapsto \Delta^{n}$, mais pas sa compl\'etion en groupo\"\i des
$n \mapsto \overline{\Delta}^{n}$. Ainsi, pour deux pr\'e-cat\'egories 
de Segal $A$ et $B$, l'ensemble des $n$-simplexes de $\underline{Hom}^{\Delta}(A,B)$
est l'ensemble des morphismes 
$$\underline{\Delta}^{n} \times A \longrightarrow B.$$
On a ainsi
$$\mathbb{R}\underline{Hom}^{\Delta}(A,B) \simeq \mathcal{I}(\mathbb{R}\underline{Hom}(A,B)).$$

Les objets
de $\uscat_{//A}$ sont par d\'efinition les morphismes de $\s$-cat\'egories fibrantes
$f : B \longrightarrow A$. Pour $f : B \longrightarrow A$ et $g : C \longrightarrow A$
deux tels objets on d\'efinit
$$Map_{\uscat//A}(f,g):=\underline{Hom}^{\Delta}(B\times \Delta^{1},C) 
\times_{\underline{Hom}^{\Delta}(A,C)\times \underline{Hom}^{\Delta}(A,C)}
\underline{Hom}^{\Delta}(B,C),$$
o\`u d'une part le morphisme 
$$\underline{Hom}^{\Delta}(B\times \Delta^{1},C) \longrightarrow
\underline{Hom}^{\Delta}(A,C)\times \underline{Hom}^{\Delta}(A,C)$$ 
est induit par la 
restriction \`a $0\coprod 1 \subset \Delta^{1}$, et d'autre part le morphisme
$$\underline{Hom}^{\Delta}(B,C) \longrightarrow \underline{Hom}^{\Delta}(A,C)\times \underline{Hom}^{\Delta}(A,C)$$
est le compos\'e
$$\xymatrix{
\underline{Hom}^{\Delta}(B,C) \ar[r]^-{g\circ } & \underline{Hom}^{\Delta}(A,C) \simeq
\{f\} \times \underline{Hom}^{\Delta}(A,C) \ar[r] & \underline{Hom}^{\Delta}(A,C) \times \underline{Hom}^{\Delta}(A,C).}
$$
Plus g\'en\'eralement, si $\{f_{i} : B_{i} \longrightarrow A\}_{2\leq i\leq n}$
est une famille finie d'objets de $\uscat/A$, on d\'efinit 
$(\uscat_{//A})(f_{1}, \dots, f_{n})$ comme \'etant le produit fibr\'e suivant
$$\underline{Hom}^{\Delta}(B\times \Delta^{n},C) 
\times_{\underline{Hom}^{\Delta}(B,C)^{n+1}}(\underline{Hom}^{\Delta}(B_{1},B_{2})\times 
\underline{Hom}^{\Delta}(B_{2},B_{3})\times \dots \times \underline{Hom}^{\Delta}(B_{n-1},B_{n})).
$$
Les ensembles simpliciaux $(\uscat_{//A})(f_{1}, \dots, f_{n})$ s'orgnisent alors
naturellement en une $\s$-cat\'egorie que nous noterons
$\uscat_{//A}$. Il faut garder \`a l'esprit que 
les morphismes entre $f : B \longrightarrow A$ et $g : C \longrightarrow A$
dans $\uscat_{//A}$ sont des paires $(u,h)$, o\`u
$u : B \longrightarrow C$ est un $\s$-foncteur et $h$ est une transformation naturelle
de $g\circ u$ dans $f$. La $\s$-cat\'egorie $\uscat_{//A}$ n'est donc pas
la $\s$-cat\'egorie comma $\uscat_{/A}$ au sens usuel. Cette derni\`ere
s'identifie en r\'ealit\'e \`a la sous-$\s$-cat\'egorie (non pleine) 
de $\uscat_{//A}$ form\'ee des morphismes $(u,h)$ comme ci-dessus avec 
$h$ une \'equivalence. \\

\begin{df}\label{d15}
Soit $A$ une $\s$-cat\'egorie fibrante poss\'edant des $\mathbb{U}$-limites.
La \emph{$\s$-cat\'egorie des $\mathbb{U}-\s$-topos $A$-pr\'e-structur\'es}
(resp. \emph{$A$-structur\'es}) est 
la sous-$\s$-cat\'egorie (non pleine) de $\s -\mathbf{Cat}_{//A}$ dont
\begin{enumerate}
\item les objets sont les
$f : T^{op} \longrightarrow A$ avec $T$ un $\mathbb{U}-\s$-topos et 
$f$ un $\s$-morphisme (resp. $\s$-morphisme qui commute aux $\mathbb{U}$-limites),
\item les morphismes sont ceux induisant des $\s$-foncteurs
g\'eom\'etriques sur les $\mathbb{U}-\s$-topos sous-jacents.
\end{enumerate}
Cette $\s$-cat\'egorie sera not\'ee $\ustop^{pr,\mathbb{U}}_{/A}$
(resp. $\ustop^{\mathbb{U}}_{/A}$)
(ou $\ustop^{pr}_{/A}$ (resp. $\ustop_{/A}$) lorsque l'univers $\mathbb{U}$ est sous-entendu).
\end{df}

Revenons sur la d\'efinition des espaces de morphismes dans $\ustop_{/A}$. 
Il s'agit, pour deux $\s$-topos $A$-structur\'es $f : T_{1}^{op} \rightarrow A$ et 
$g : T_{2}^{op} \rightarrow A$
 de consid\'erer la projection naturelle
$$Map_{\s -\mathbf{Cat}//A}(f,g) \longrightarrow Map_{\s -\mathbf{Cat}}(T_{1}^{op},T_{2}^{op})\simeq
Map_{\s -\mathbf{Cat}}(T_{1},T_{2}).$$
Par d\'efinition l'ensemble simplicial $Map_{\ustop_{/A}}(f,g)$ est 
le produit fibr\'e
$$\xymatrix{
Map_{\s -\mathbf{Cat}//A}(f,g) \ar[r] & Map_{\s -\mathbf{Cat}}(T_{1},T_{2})  \\
Map_{\ustop_{/A}}(f,g) \ar[r] \ar[u] & Map^{geom}_{\s -\mathbf{Cat}}(T_{1},T_{2}) \ar[u],}$$
o\`u $Map^{geom}_{\s -\mathbf{Cat}}(T_{1},T_{2})$ d\'esigne la r\'eunion des composantes
connexes correspondantes aux morphismes g\'eom\'etriques. Ainsi, 
$Map_{\ustop_{/A}}(f,g)$ est aussi une r\'eunion de composante connexe
dans $Map_{\s -\mathbf{Cat}//A}(f,g)$. \\

Il existe un $\s$-foncteur d'oubli
$$\ustop_{/A} \longrightarrow \ustop$$
qui \`a un objet $f : T^{op} \longrightarrow A$
associe le $\s$-topos $T$. La fibre homotopique de ce $\s$-foncteur d'oubli, 
pris en un $\s$-topos fix\'e $T$, est naturellement \'equivalente \`a la 
$\s$-cat\'egorie $Ch(T,A)$, des champs sur $T$ \`a valeurs dans $A$ (voir d\'efinition \ref{d5}). \\

Nous allons maintenant sp\'ecifier deux $\s$-cat\'egories $A$ et obtenir
ainsi deux notions de $\s$-topos structur\'es qui seront pour nous
fondamentaux. 

Nous avons d\'ej\`a d\'efini la $\s$-cat\'egorie $\uscms$ des $\s$-cat\'egories
mono\"\i dales sym\'etriques, ainsi que la sous-$\s$-cat\'egorie
pleine $\uscms_{rig} \subset \uscms$ des $\s$-CMS rigides. Si on fixe
un univers $\mathbb{U}$ et que l'on se restreint aux $\s$-CMS $\mathbb{U}$-petites, 
alors $\uscms_{rig}$ et  $\uscms$ sont toutes deux $\mathbb{U}$-localement pr\'esentables, et 
en particulier poss\`ede des $\mathbb{U}$-limites. 

Nous d\'efinissons aussi une $\s$-cat\'egorie $\s -\underline{CMon}$, de la fa\c{c}on suivante. 
Notons $\mathbf{SEns}^{\Gamma}$ la cat\'egorie des foncteurs de $\Gamma$ dans $\mathbf{SEns}$. 
On consid\`ere la sous-cat\'egorie pleine $\mathbf{SEns}^{\Gamma}_{sp}$ des
foncteurs sp\'eciaux, c'est \`a dire des foncteurs $F$ qui v\'erifient que 
pour tout $n$ le morphisme
$$F([n]) \longrightarrow F([1])^{n}$$
est une \'equivalence. On pose alors 
$$\s -\underline{CMon} := L(\mathbf{SEns}^{\Gamma}_{sp}),$$
o\`u la localisation est effectu\'ee le long des \'equivalences niveaux par niveaux
dans la cat\'egorie des foncteurs $\mathbf{SEns}^{\Gamma}$. Tout comme nous l'avons fait 
pour $\uscms$ dans \S 2.1 (voir la d\'efinition \ref{d8}) 
on montre que $\s -\underline{CMon}$ est une $\s$-cat\'egorie
$\mathbb{U}$-localement pr\'esentable (si l'on se restreint aux foncteurs
$\Gamma \longrightarrow \mathbf{SEns}_{\mathbb{U}}$).

\begin{df}\label{d16}
\begin{enumerate}
\item La $\s$-cat\'egorie des \emph{$\mathbb{U}-\s$-topos pr\'e-catannel\'es} (resp.
\emph{catannel\'es})
est 
$$\ustop^{pr,\mathbb{U}}_{\otimes}:=\ustop^{pr,\mathbb{U}}_{/\uscms_{\mathbb{U}}}$$
$$(resp. \; \ustop^{\mathbb{U}}_{\otimes}:=\ustop^{\mathbb{U}}_{/\uscms_{\mathbb{U}}}.$$
\item La $\s$-cat\'egorie des \emph{$\s$-topos pr\'e-catannel\'es rigides}
(resp. \emph{catannel\'es rigides})
est 
$$\ustop^{pr,\mathbb{U}}_{\otimes,rig}:=\ustop^{pr,\mathbb{U}}_{/\uscms_{\mathbb{U},rig}}$$
$$(resp. \; \ustop^{\mathbb{U}}_{\otimes,rig}:=\ustop^{\mathbb{U}}_{/\uscms_{\mathbb{U},rig}}).$$

\item La $\s$-cat\'egorie des \emph{$\mathbb{U}-\s$-topos pr\'e-mono\"\i del\'es}
(resp. \emph{mono\"\i d\'es})
est 
$$\ustop^{pr,\mathbb{U}}_{Mon}:=\ustop^{pr,\mathbb{U}}_{/\usmon_{\mathbb{U}}}$$
$$(resp. \; \ustop^{\mathbb{U}}_{Mon}:=\ustop^{\mathbb{U}}_{/\usmon_{\mathbb{U}}}).$$
\end{enumerate}
Comme nous en avons l'habitude nous noterons simplement 
$$\ustop^{pr}_{\otimes} \qquad \ustop^{pr}_{\otimes,rig} \qquad 
\ustop^{pr}_{Mon}$$
$$\ustop_{\otimes} \qquad \ustop_{\otimes,rig} \qquad 
\ustop_{Mon}$$
lorsqu'il n'est pas n\'ecessaire de sp\'ecifier l'univers $\mathbb{U}$. 
\end{df}

Afin de rappeler les notations utilis\'ees pour la notion d'espaces annel\'es les objets
d'une des $\s$-cat\'egories pr\'ec\'edentes seront not\'es
sous forme de couple $(T,\mathcal{A})$. Une telle notation d\'esigne un 
$\mathbb{U}-\s$-topos $T$, fibrant comme $\s$-cat\'egorie, 
et un $\s$-foncteur 
$$\mathcal{A} : T^{op} \longrightarrow A$$
(qui \'eventuellement commute aux $\mathbb{U}$-limites), o\`u $A$ est l'une des $\s$-cat\'egories
$\uscms_{\mathbb{U}}$, $\uscms_{\mathbb{U},rig}$, $\usmon$. \\

Nous utiliserons la proposition suivante qui permet de construire 
des exemples de $\s$-topos structur\'es par le proc\'ed\'e de 
\emph{champs associ\'e} \`a partir de $\s$-topos pr\'e-structur\'es. 

\begin{prop}\label{p7}
Le $\s$-foncteur d'inclusion
$$i : \ustop_{/A} \longrightarrow \ustop^{pr}_{/A}$$
poss\`ede un adjoint \`a gauche
$$a : \ustop^{pr}_{/A} \longrightarrow \ustop_{/A}.$$
\end{prop}

\textbf{Preuve --} Soit $(T,\mathcal{A})$ un objet de $\ustop^{pr}_{/A}$.
On consid\`ere le $\s$-foncteur $\mathcal{A} : T^{op} \longrightarrow A$
comme un objet de $\rh (T^{op},A)$. On consid\`ere alors son champ
associ\'e
$a(\mathcal{A}) \in Ch(T,A)$ (voir la d\'efinition \ref{d5}), ce qui fournit un 
nouveau $\s$-foncteur $a(\mathcal{A}) : T^{op} \longrightarrow A$
muni d'une transformation naturelle $h : \mathcal{A} \rightarrow a(\mathcal{A})$. 
Ceci d\'efinit un morphisme $h : (T,\mathcal{A}) \longrightarrow (T,a(\mathcal{A}))$ 
dans $\ustop^{pr}_{/A}$. 

Soit maintenant $(T',\mathcal{A}')$ un objet de $\ustop_{/A}$ et consid\'erons 
le morphisme
$$h^{*} : Map((T,a(\mathcal{A})),(T',\mathcal{A}')) \longrightarrow 
Map((T,\mathcal{A}),(T',\mathcal{A}')).$$
Il s'agit de montrer que le morphisme $h^{*}$ est un isomorphisme dans $\mathrm{Ho}(\mathbf{SEns})$. En effet, 
ceci montrera que pour tout $(T,\mathcal{A})$, le $\s$-foncteur 
$$Map((T,\mathcal{A}),-) : \ustop_{/A} \longrightarrow \T$$
est corepr\'esentable, et donc que $i$ poss\`ede un adjoint \`a gauche. 

On dispose d'un diagramme commutatif dans $\mathrm{Ho}(\mathbf{SEns})$
$$\xymatrix{
Map((T,a(\mathcal{A})),(T',\mathcal{A}')) \ar[r]^-{h^{*}} \ar[rd] & Map((T,\mathcal{A}),(T',\mathcal{A}')) \ar[d] \\
 & Map^{geom}(T,T'). }$$
Fixons $u : T \longrightarrow T'$ un morphisme g\'eom\'etrique. Le morphisme induit par $h^{*}$
sur les fibres homotopiques en $u$ s'identifie au morphism
$$h^{*} : Map_{\rh (T^{op},A)}(a(\mathcal{A}),u^{*}(\mathcal{A}')) \longrightarrow 
Map_{\rh (T^{op},A)}(\mathcal{A},u^{*}(\mathcal{A}')).$$
Ce morphisme est bien un isomorphisme dans $\mathrm{Ho}(\mathbf{SEns})$ par la propri\'et\'e universelle
de $a(\mathcal{A})$ comme objet dans $\rh (T^{op},A)$.  \hfill $\Box$ \\

Dans la suite de ce travail nous nous int\'eresserons 
uniquement aux cat\'egories $[\ustop_{/A}^{pr}]$ et $[\ustop_{/A}]$. Cela simplifiera largement les
d\'etails techniques, bien que nous n\'egligerons de la sorte 
des donn\'ees de fonctorialit\'e sup\'erieures. Ces cat\'egories poss\`edent une
description explicite relativement simple. Donnons la description 
de $[\ustop_{/A}]$ \`a titre d'exemples, le cas des $\s$-topos
pr\'e-structur\'es se d\'ecrivant de fa\c{c}on similaire en om\'etant 
la condition de commutation aux limites.
Les objets de $[\ustop_{/A}]$ sont 
les couples $(T,\mathcal{A})$,  consistant en 
un ($\mathbb{U}-$)$\s$-topos $T$ fibrant comme $\s$-cat\'egorie, et 
$\mathcal{A} : T^{op} \rightarrow A$, un $\s$-foncteur qui commute aux ($\mathbb{U}$-)limites. 
L'ensemble des morphismes $[(T,\mathcal{A}),(T',\mathcal{A}')]$ se 
d\'ecrit comme suit. Ses \'el\'ements sont les classes
d'\'equivalence de couples $(u,h)$, avec 
$u : T^{op} \longrightarrow (T')^{op}$ un morphisme g\'eom\'etrique et 
$h : T^{op}\times \Delta^{1} \longrightarrow A$ est un $\s$-foncteur
avec $h_{0}=\mathcal{A}'\circ u$ et $h_{1}=\mathcal{A}$. Deux tels couples
$(u,h)$ et $(v,k)$ sont \'equivalents s'il existe deux $\s$-foncteurs
$$\psi : T^{op}\times \overline{\Delta}^{1} \longrightarrow (T')^{op} \qquad
\phi : T\times \Delta^{2} \longrightarrow A$$
tels que 
$$\psi_{0}=u \qquad \psi_{1}=v$$
et $\phi$ est un \'el\'ement dans $\underline{Hom}(T^{op},A)([2])$ dont 
les trois projections dans $\underline{Hom}(T^{op},A)([1])$
sont \'egales \`a $h$, $k$ et $\mathcal{A}'\circ \psi$. En d'autres termes
$\phi$ est une $2$-cellule faisant commuter le diagramme suivant 
dans $\underline{Hom}(T^{op},A)$
$$\xymatrix{
\mathcal{A}'\circ u \ar[r]^-{h} \ar[d]_-{\mathcal{A}'\circ \psi} & \mathcal{A} \\
\mathcal{A}'\circ v. \ar[ru]_-{k} & }$$
La composition des morphismes dans $[\ustop_{/A}]$ 
s'effectue de la fa\c{c}on suivante. Pour $(u,h) : (T,\mathcal{A}) \longrightarrow 
(T',\mathcal{A}')$ et $(v,k) : (T,\mathcal{A}) \longrightarrow 
(T',\mathcal{A}')$ deux tels morphismes on consid\`ere
$$(v,k)\circ (u,h):=(v\circ u,l) : (T,\mathcal{A}) \longrightarrow (T',\mathcal{A}')$$
o\`u $l : T^{op}\times \Delta^{1} \longrightarrow A$ est un choix pour une composition
$k\circ h$ comme morphisme dans $\underline{Hom}(T^{op},A)$. On v\'erifie que le morphisme
$(v\circ u,l)$ ainsi d\'efini est ind\'ependant du choix de $l$ et que cela d\'efinit 
une composition associative. 

\subsection{Homologie de Hochschild et homologie cyclique}

Soit $(T,\mathcal{A})$ un $\s$-topos catannel\'e. Nous allons
construire deux $\s$-topos mono\"\i del\'es $(T,HH^{\mathcal{A}})$ 
et $(T,HC^{\mathcal{A}})$, d\'eduits de $(T,\mathcal{A})$ 
par les proc\'ed\'es suivants. 

Nous commen\c{c}ons par consid\'erer $\underline{Hom}(T,T)$, la $\s$-cat\'egorie
des endofoncteurs de $T$ (rappelons que par d\'efinition 
$T$ est une $\s$-cat\'egorie fibrante). La $\s$-cat\'egorie $T$ poss\'edant
des $\mathbb{U}$-limites il en est de m\^eme de $\underline{Hom}(T,T)$. Ainsi, 
pour tout ensemble simplicial $K$ et tout objet $f \in \underline{Hom}(T,T)$ 
on peut former l'objet $f^{K}\in \underline{Hom}(T,T)$. Rappelons que cet objet vient avec un 
morphisme d'ensembles simpliciaux $K \longrightarrow Map(f^{K},f)$, tel que
pour tout $g \in \underline{Hom}(T,T)$ le morphisme induit
$$Map(g,f^{K})\longrightarrow Map_{\mathbf{SEns}}(Map(f^{K},f),Map(g,f)) \longrightarrow Map_{\mathbf{SEns}}(K,Map(g,f))$$
soit un isomorphisme dans $\mathrm{Ho}(\mathbf{SEns})$. On applique cette construction 
avec $f=id$ et $K=S^{1}=B\mathbb{Z}$, et on obtient ainsi un $\s$-foncteur $id^{S^{1}} : T \longrightarrow T$, 
muni d'un morphisme $S^{1} \longrightarrow Map(id^{S^{1}},id)$. 
D'autre part on dispose du $\s$-foncteur
$$\mathcal{M}\circ \widetilde{\mathcal{LI}} : \uscms \longrightarrow \usmon$$
d\'efini en fin de \S 2.3. Ainsi, on peut construire un nouvel 
objet $(T,\mathcal{M}\circ \widetilde{\mathcal{LI}}\mathcal{A}\circ id^{K})$ dans
$[\ustop_{Mon}^{pr}]$.

\begin{df}\label{d17}
Pour un $\s$-topos catannel\'e $(T,\mathcal{A})$, le $\s$-topos pr\'e-mono\"\i del\'e
$(T,HH_{pr}^{\mathcal{A}})$ est d\'efini par 
$$(T,HH_{pr}^{\mathcal{A}}):=(T,\mathcal{M}\circ \widetilde{\mathcal{LI}}\mathcal{A}\circ id^{S^{1}}).$$
Nous noterons 
$$(T,HH^{\mathcal{A}}):=(T,a(HH_{pr}^{\mathcal{A}}))$$
le $\s$-topos mono\"\i del\'e associ\'e.
\end{df}

En termes plus concrets le $\s$-foncteur
$$HH_{pr}^{\mathcal{A}} : T^{op} \longrightarrow \usmon$$
envoie un objet $X\in T$ sur le $\s$-mono\"\i de 
$End_{\mathcal{A}(X^{S^{1}})}(\mathbf{1})$, des endomorphismes de l'unit\'e
dans la $\s$-CMS $\mathcal{A}(X^{S^{1}})$. Ce nouveau $\s$-foncteur n'est en g\'en\'eral
plus un champ et son champ associ\'e est $HH^{\mathcal{A}}$. \\

Revenons au $\s$-endofoncteur $id^{S^{1}}$. On peut promouvoir cette construction en un 
$\s$-foncteur
$$\T^{op} \longrightarrow \underline{Hom}(T,T)$$
qui commute aux limites et qui envoie $*$ sur l'objet $id$. En effet, 
un tel $\s$-foncteur existe et est unique dans $[\uscat]$ car $\T=\widehat{*}$ (voir \S 1.4).
On trouve en particulier un morphisme bien d\'efini 
$$BS^{1}=K(\mathbb{Z},2) \longrightarrow \T \longrightarrow \underline{Hom}(T,T),$$
qui envoie le point de base de $BS^{1}=K(\mathbb{Z},2)$ sur 
l'objet $id^{S^{1}}$. En d'autre termes, on dispose d'une action naturelle
du groupe simplicial $S^{1}$ sur le $\s$-foncteur $id^{S^{1}}$. Cette action est 
bien \'evidemment induite par l'action de $S^{1}$ sur lui-m\^eme par translation. 
Ainsi, le $\s$-foncteur compos\'e
$$HH_{pr}^{\mathcal{A}} : T^{op} \longrightarrow \usmon$$
poss\`ede un rel\`evement naturel
$$\widetilde{HH_{pr}^{\mathcal{A}}} : T^{op} \longrightarrow
\rh (BS^{1},\usmon)=:S^{1}-\usmon,$$
en un $\s$-foncteur vers la $\s$-cat\'egorie des $\s$-mono\"\i des
$S^{1}$-\'equivariants. En composant avec le $\s$-foncteur 
de champ associ\'e, on trouve un rel\`evement
$$\widetilde{HH^{\mathcal{A}}} : T^{op} \longrightarrow
\rh (BS^{1},\usmon)=:S^{1}-\usmon,$$
du $\s$-foncteur $HH^{\mathcal{A}}$.
En composant
avec le $\s$-foncteurs
des points fixes homotopiques (qui est le $\s$-foncteur $lim_{BS^{1}}$)
$$(-)^{hS^{1}} : S^{1}-\usmon \longrightarrow \usmon ,$$
on trouve deux nouveaux $\s$-foncteur
$$HC_{pr}^{\mathcal{A}} : T^{op} \longrightarrow \usmon$$
$$HC^{\mathcal{A}} : T^{op} \longrightarrow \usmon.$$
Par d\'efinition on a pour tout objet $X \in T$
$$HC_{pr}^{\mathcal{A}}(X)\simeq 
\widetilde{HH}_{pr}^{\mathcal{A}}(X)^{hS^{1}}\qquad HC^{\mathcal{A}}(X)\simeq 
\widetilde{HH}^{\mathcal{A}}(X)^{hS^{1}}.$$

\begin{df}\label{d18}
Pour un $\s$-topos catannel\'e $(T,\mathcal{A})$, le $\s$-topos pr\'e-mono\"\i del\'e
$(T,HC_{pr}^{\mathcal{A}})$ est d\'efini par 
$$(T,HC_{pr}^{\mathcal{A}}):=(T,(\widetilde{HH}_{pr}^{\mathcal{A}})^{hS^{1}}).$$
\end{df}

Notons que la projection naturelle des points fixes sur le $\s$-foncteur d'oubli
induit des morphismes dans $[\ustop_{Mon}]$
$$(T,HH_{pr}^{\mathcal{A}}) \longrightarrow (T,HC_{pr}^{\mathcal{A}}) \qquad 
(T,HH^{\mathcal{A}}) \longrightarrow (T,HC^{\mathcal{A}}).$$
Pour $X \in T$ peut \'ecrire
$$HH_{pr}^{\mathcal{A}}\simeq End_{\mathcal{A}}(X^{S^{1}}) \qquad
HC_{pr}^{\mathcal{A}}\simeq End_{\mathcal{A}}(X^{S^{1}})^{hS^{1}}.$$
Plus g\'en\'eralement, il n'est pas difficile de voir que 
$(T,\mathcal{A}) \mapsto (T,HH_{pr}^{\mathcal{A}})$
et $(T,\mathcal{A}) \mapsto (T,HC_{pr}^{\mathcal{A}})$
d\'efinissent des foncteur
$$[\ustop^{pr}_{\otimes,rig}] \longrightarrow [\ustop^{pr}_{Mon}].$$
De m\^eme, on dispose de deux foncteurs
$$[\ustop_{\otimes,rig}] \longrightarrow [\ustop_{Mon}].$$
Ces foncteurs sont reli\'es par les morphismes
$$(T,HH_{pr}^{\mathcal{A}}) \longrightarrow (T,HC_{pr}^{\mathcal{A}}) \qquad 
(T,HH^{\mathcal{A}}) \longrightarrow (T,HC^{\mathcal{A}})$$
qui induisent des transformations naturelles. 

\begin{df}\label{d19}
Pour $(T,\mathcal{A})$ un $\s$-topos catannel\'e et $X \in T$ un objet.
\begin{enumerate}
\item Le \emph{$\s$-mono\"\i de de pr\'e-homologie de Hochschild de $X$}
(resp. \emph{homologie de Hochschild}) est
$$HH_{pr}^{\mathcal{A}}(X) \qquad (resp. \; HH^{\mathcal{A}}(X)).$$
\item Le \emph{$\s$-mono\"\i de de pr\'e-homologie cyclique de $X$}
(resp. \emph{homologie de cyclique}) est
$$HC_{pr}^{\mathcal{A}}(X) \qquad (resp. \; HC^{\mathcal{A}}(X)).$$
\end{enumerate}
\end{df}

Nous expliquerons au \S 4.2 que notre d\'efinition de l'homologie cyclique est en r\'ealit\'e 
une g\'en\'eralisation de l'homologie cyclique \emph{n\'egative} et non 
de l'homologie cyclique au sens propre. Notre caract\`ere de Chern sera construit \`a valeurs
dans $HC_{pr}$ et nous n'aurons jamais \`a consid\'erer un analogue de l'homologie
cyclique non n\'egative. Cette derni\`ere pourrait \^etre d\'efinie en rempla\c{c}ant 
les points fixes homotopiques par $S^{1}$ par des coinvariants homotopiques (i.e. 
$colim_{BS^{1}}$ au lieu de $lim_{BS^{1}}$). Nous ne voyons cependant pas
la pertinence d'une telle d\'efinition dans le contexte g\'en\'eral
expos\'e ici. Ainsi, pour nous, \emph{homologie cyclique} sera toujours
une expression faisant r\'ef\'erence \`a l'homologie cyclique n\'egative.

\subsection{Caract\`ere de Chern}

Soit $(T,\mathcal{A})$ un $\s$-topos catannel\'e rigide. Nous avons d\'ej\`a vu
qu'on pouvait lui associer deux $\s$-topos pr\'e-mono\"\i del\'es
$(T,HH_{pr}^{\mathcal{A}})$ et $(T,HC_{pr}^{\mathcal{A}})$. 
Nous allons maintenant en d\'efinir un troisi\`eme. Il s'agit 
de composer $\mathcal{A}$ avec le foncteur 
$$\mathcal{I}\circ \mathcal{M} : \uscms \longrightarrow \usmon,$$
pour obtenir un nouveau $\s$-foncteur
$$\mathcal{I}\circ \mathcal{M} \circ \mathcal{A} : T^{op} \longrightarrow \usmon$$
et donc un nouveau $\s$-topos pr\'e-mono\"\i del\'e que nous noterons
$(T,|\mathcal{A}|)$. 

Notons que $|\mathcal{A}| : T^{op} \longrightarrow \usmon$ 
commute aux limites et donc est repr\'esentable par un objet
$\mathcal{E} \in \rh(\Gamma,T)$. De m\^eme, le $\s$-foncteur
$$\widetilde{\mathcal{LI}}\circ  \mathcal{M}\circ \mathcal{A} : T^{op} \longrightarrow S^{1}-\usmon
\subset \rh (\Gamma \times BS^{1},\T)$$
est repr\'esentable par l'objet 
$$\mathcal{E}^{S^{1}} \in \rh(\Gamma \times BS^{1},T).$$ 
Enfin, 
le morphisme de trace cyclique multiplicative $Tr^{S^{1}}_{\otimes}$ (voir la d\'efinition \ref{d14})
induit une transformation naturelle
$$\widetilde{\mathcal{LI}} \circ \mathcal{M} \circ \mathcal{A}  \longrightarrow 
\widetilde{End}_{\mathcal{A}}(\mathbf{1}).$$
Cette transformation naturelle correspond \`a l'aide du lemme de Yoneda \`a un
\'el\'ement
$$Tr^{S^{1}}_{\otimes} \in 
\pi_{0}(End_{\mathcal{A}(\mathcal{E}^{S^{1}})}(\mathbf{1})^{hS^{1}})=
\pi_{0}(HC_{pr}^{\mathcal{A}}(\mathcal{E})).$$
Par le lemme de Yoneda cet \'el\'ement est lui-m\^eme donn\'e par une
transformation naturelle
$$Ch^{pr} : |\mathcal{A}| \longrightarrow HC_{pr}^{\mathcal{A}}$$
bien d\'efinie dans $[\rh (T^{op},\usmon)]$. 

\begin{df}\label{d20}
Soit $(T,\mathcal{A})$ un $\s$-topos catannel\'e rigide.
\begin{enumerate}
\item
Le \emph{pr\'e-caract\`ere de Chern du $\s$-topos catannel\'e rigide $(T,\mathcal{A})$}
est le morphisme
$$Ch^{pr} : |\mathcal{A}| \longrightarrow HC_{pr}^{\mathcal{A}}$$
dans $[\rh (T^{op},\usmon)]$.
\item Le \emph{caract\`ere de Chern du $\s$-topos catannel\'e rigide $(T,\mathcal{A})$}
est le morphisme compos\'e
$$Ch : \xymatrix{
|\mathcal{A}| \ar[r]^-{Ch^{pr}} & HC_{pr}^{\mathcal{A}} \ar[r] & 
HC^{\mathcal{A}}.
}$$
\item La \emph{pr\'e-trace du $\s$-topos catannel\'e rigide $(T,\mathcal{A})$}
est le morphisme compos\'e
$$Tr^{pr} : \xymatrix{
|\mathcal{A}| \ar[r]^-{Ch^{pr}} & HC_{pr}^{\mathcal{A}} \ar[r] & 
HH_{pr}^{\mathcal{A}}.
}$$
\item La \emph{trace du $\s$-topos catannel\'e rigide $(T,\mathcal{A})$}
est le morphisme compos\'e
$$Tr : \xymatrix{
|\mathcal{A}| \ar[r]^-{Ch} & HC^{\mathcal{A}} \ar[r] & 
HH^{\mathcal{A}}.
}$$
\end{enumerate}
\end{df}

On v\'erifie sans peine que le caract\`ere de Chern induit une transformation naturelle 
entre les deux foncteurs
$$[\ustop_{\uscms_{rig}}] \longrightarrow [\ustop_{\usmon}]$$
qui envoient $(T,\mathcal{A})$ sur $(T,|\mathcal{A}|)$ et 
sur $(T,HC^{\mathcal{A}})$ respectivement. En d'autres termes 
la construction de $Ch$ est fonctorielle d'une part en 
$\mathcal{A}$ et d'autre part en $T$ (remarquons que la fonctorialit\'e
en $T$ utilise l'exactitude \`a gauche des adjoints des morphismes g\'eom\'etriques
afin de comparer les objets de la forme $X^{S^{1}}$). Nous laissons le soins
au lecteur d'\'ecrire les d\'etails de ces fonctorialit\'es.

\section{Exemples et applications}

Dans cette derni\`ere section nous donnons quatres exemples
de contextes d'applications du caract\`ere de Chern construit
dans la section pr\'ec\'edente. Le premier est relativement formel
et consiste \`a remarquer qu'un $\s$-topos annel\'e donne lieu \`a deux
$\s$-topos catannel\'es, obtenus en consid\'erant les groupes additifs et multiplicatifs
sous-jacents au faisceau structural.
Dans le second exemple nous revenons sur le caract\`ere de Chern
des fibr\'es vectoriels alg\'ebriques et plus g\'en\'eralement
des complexes parfaits. Nous donnons la preuve d\'etaill\'e de la comparaison
avec le caract\`ere de Chern usuel pour les fibr\'es vectoriels.
Les deux derniers exemples
sont \`a nos yeux les plus importants car ils concernent 
le caract\`ere de Chern d'une famille de dg-cat\'egories 
dont l'existence avait \'et\'e avanc\'ee dans \cite{tove}.
Nous montrons tout d'abord comment on peut associer \`a toute
famille alg\'ebrique de dg-cat\'egories compactement engendr\'ees
sur une base $X$ (disons un sch\'ema ou plus g\'en\'eralement 
un champs alg\'ebrique) un complexe quasi-coh\'erent et 
$S^{1}$-\'equivariant sur l'objet des lacets $LX$. Cet 
objet donne lieu \`a une sorte de $\mathcal{D}_{X}$-module
filtr\'e qu'il 
faut voir comme la variation de structures de Hodge
non-commutatives induite sur l'homologie cyclique de la famille
(au sens de \cite{kakopa}). Cette variation comprend en particulier la donn\'ee
d'une connexion de Gauss-Manin, ce qui nous permet de g\'en\'eraliser
les constructions de \cite{ge,dotats} au cas des familles de dg-cat\'egories sur des
bases tr\`es g\'en\'erales. L\`a encore, de nombreux \'enonc\'es de
comparaison, par exemple entre quasi-coh\'erents $S^{1}$-\'equivariants
sur $LX$ et $\mathcal{D}_{X}$-modules, ne sont pas d\'etaill\'es: nous
reviendrons sur ces points, qui m\'eritent un traitement
s\'epar\'e, dans un travail ult\'erieur.
Enfin, lorsque la famille de dg-cat\'egories en question 
est \`a fibres satur\'ees nous expliquons comment on peut lui
associer un caract\`ere de Chern secondaire \`a valeurs
dans une nouvelle th\'eorie homologique que nous introduisons: l'homologie
cyclique secondaire. Une \'etude plus approfondie de ce
caract\`ere de Chern secondaire fera l'objet d'un travail futur. \\

Avant de d\'ecrire ces exemples nous pr\'esentons notre $\s$-topos
favori, \`a savoir celui des champs d\'eriv\'es, et qui sera utilis\'e
pour tous nos exemples. Comme nous le signalons \`a plusieurs reprises on 
aurait tout aussi bien pu consid\'erer les versions diff\'erentielles, ou
complexes analytiques, de ce $\s$-topos (voir \cite[\S 4.4]{lu3}).
Le cadre alg\'ebrique est cependant plus simple et de plus rend compte
de tous les ph\'enom\`enes int\'eressants.

Fixons nous $k$ un anneau commutatif. On consid\`ere 
$sk-Comm$, la cat\'egorie des $k$-alg\`ebres commutatives simpliciales.
En localisant les \'equivalences faibles, qui sont 
les morphismes induisant des \'equivalences faibles
sur les ensembles simpliciaux sous-jacents, on trouve
une $\s$-cat\'egorie $L(sk-Comm)$. On pose alors
$$dAff_{k}:=L(sk-Comm)^{op},$$
que nous appellerons l'$\s$-cat\'egorie des
k-sch\'emas d\'eriv\'es affines. Il existe sur
la cat\'egorie $[dAff_{k}]$ une topologie de Grothendieck, 
d\'efinie \`a l'aide des recouvrements \'etales (voir
par exemple \cite[\S 2.2]{hagII}). Cette topologie d\'efinit une
topologie sur l'$\s$-cat\'egorie $dAff_{k}$ et donne ainsi lieu 
\`a un $\s$-topos $dAff_{k}^{\sim,et}$, appel\'e l'\emph{$\s$-topos
des k-champs d\'eriv\'es} (il est not\'e $\mathbf{dSt}(k)$ dans
\cite[\S 4]{to3}). Notons que l'on a
$$dAff_{k}^{\sim,et} \simeq L(k-D^{-}Aff^{\sim,et}),$$
o\`u $k-D^{-}Aff^{\sim,et}$ est le topos de mod\`eles
des champs pour la topologie \'etale sur $k-D^{-}Aff=(sk-Comm)^{op}$.

L'$\s$-cat\'egorie $dAff_{k}^{\sim,et}$ poss\`ede 
la cat\'egorie des $k$-sch\'emas, et la $2$-cat\'egorie
des $k$-champs alg\'ebriques, comme sous-$\s$-cat\'egories pleines (voir
\cite[\S 2.2.4]{hagII}). Ainsi, nous consid\'ererons toujours
les $k$-sch\'emas et les k-champs alg\'ebriques comme des
objets de $dAff_{k}^{\sim,et}$. Notons cependant que ces sous-$\s$-cat\'egories
ne sont pas stables par limites, et en particulier par produit fibr\'e. 
Tous les produits fibr\'es que nous consid\'ererons
seront construits dans $dAff_{k}^{\sim,et}$. 

\subsection{$\s$-Topos annel\'es et $\s$-Topos catannel\'es}

Consid\'erons $Sp^{\Sigma}$ la cat\'egorie de mod\`eles des spectres sym\'etriques
de \cite{hss}, munie de sa structure positive de \cite{shi}. Les mono\"\i des commutatifs
dans cette cat\'egorie de mod\`eles forment une cat\'egorie $Comm(Sp^{\Sigma})$, 
dont les objets sont les anneaux en spectres sym\'etriques commutatifs. 
Cette cat\'egorie h\'erite d'une structure de mod\`eles pour la quelle les
fibrations et les \'equivalences sont d\'efinies sur les objets sous-jacents
dans $Sp^{\Sigma}$. Nous noterons 
$\uscomm$ la $\s$-cat\'egorie obtenue par localisation
$$\uscomm:=L(Comm(Sp^{\Sigma})).$$
Par d\'efinition un $\s$-topos annel\'e est un $\s$-topos 
structur\'e dans $\uscomm$. La cat\'egorie
des $\s$-topos annel\'es est donc $[\ustop_{\uscomm}]$.
Nous allons  construire deux $\s$-foncteurs
$$\mathbb{G}_{a} : \uscomm \longrightarrow \uscms_{rig} \qquad 
\mathbb{G}_{m} : \uscomm \longrightarrow \uscms_{rig},$$
qui de plus commutent aux limites. Ils induiront ainsi deux foncteurs
$$[\ustop_{\uscomm}] \longrightarrow [\ustop_{\uscms_{rig}}]$$
qui envoient respectivement $(T,\mathcal{A})$ sur 
$(T,\mathbb{G}_{a}\circ \mathcal{A})$ et $(T,\mathbb{G}_{m}\circ \mathcal{A})$.

Pour ce qui est du $\s$-foncteur $\mathbb{G}_{a}$ (groupe additif sous-jacent), 
on commence par consid\'erer le foncteur d'oubli $Comm(Sp^{\Sigma}) \longrightarrow
Sp^{\Sigma}$ qui induit un $\s$-foncteur sur les localis\'es 
$$\uscomm \longrightarrow L(Sp^{\Sigma}).$$
On compose alors avec l'adjonction de Quillen 
$$\mathbf{SEns}^{\Gamma} \longleftrightarrow Sp^{\Sigma}$$
entre $\Gamma$-espaces et spectres sym\'etriques (voir \cite{mmss}). 
En composant avec $\Pi_{\infty} : \T \longrightarrow \uscat$
on trouve finalement un $\s$-foncteur
$$\mathbb{G}_{a} : 
\xymatrix{uscomm \ar[r] & L(Sp^{\Sigma}) \ar[r] & L(\mathbf{SEns}^{\Gamma})\simeq \rh (\Gamma,\T) \ar[r]^-{\Pi_{\infty}} & 
\rh (\Gamma,\uscat).}$$
Ce foncteur se factorise par la sous-$\s$-cat\'egorie $\uscms \subset \rh (\Gamma,\uscat)$. 
De plus, $\mathbb{G}_{a}(R)$ est une $\s$-CMS dont la cat\'egorie homotopique $[\mathbb{G}_{a}(R)]$
est un groupo\"\i de muni d'une loi mono\"\i dale pour la quelle tous les objets sont 
inversibles. La cat\'egorie mono\"\i dale sym\'etrique $[\mathbb{G}_{a}(R)]$ est donc toujours rigide, 
et on obtient ainsi un $\s$-foncteur
$$\mathbb{G}_{a} : \uscomm \longrightarrow \uscms_{rig}.$$
Une autre fa\c{c}on de voir ce $\s$-foncteur est en utilisant l'\'equivalence entre
spectres connectifs et $\s$-CMS dont la loi mono\"\i dale est inversible mentionn\'ee au \S 2.1. 
Pour $R \in \uscomm$, la $\s$-CMS $\mathbb{G}_{a}(R)$ est celle correspondant 
au rev\`etement connectif du spectre sous-jacent \`a $R$. Finalement, on voit par construction que
$\mathbb{G}_{a}$ est obtenu par composition de $\s$-foncteurs qui  
commutent tous aux limites, il commute ainsi lui aussi avec les limites.

Construisons maintenant le $\s$-foncteur $\mathbb{G}_{m}$. Pour cela, on consid\`ere pour tout
$R \in Comm(Sp^{\Sigma})$ la cat\'egorie $R-Mod^{cof}$ des $R$-modules cofirants
(toujours pour la structure positive de \cite{shi}). 
La cat\'egorie $R-Mod^{cof}$ est une cat\'egorie mono\"\i dale 
sym\'etrique munie d'une notion compatible d'\'equivalence faible $W$. On obtient ainsi
un pseudo-foncteur $Comm(Sp^{\Sigma}) \longrightarrow \mathbf{Cat}^{\otimes}$, qui envoie
$R$ sur $R-Mod^{cof}$ et un morphisme $R \longrightarrow R'$ sur le changement de base
$R'\wedge_{R}-$. On rectifie ce pseudo-foncteur en le rempla\c{c}ant par le 
foncteur de ses sections cocart\'esiennes (voir \S 1.3 par exemple) pour obtenir un vrai foncteur
$$Comm(Sp^{\Sigma}) \longrightarrow \mathbf{Cat}^{\otimes}.$$
Ce foncteur vient avec une notion d'\'equivalence faible induite, encore not\'ee $W$, et 
en localisant on obtient ainsi un foncteur
$$Comm(Sp^{\Sigma}) \longrightarrow \scms.$$
Ce foncteur envoie lui-m\^eme \'equivalences dans $Comm(Sp^{\Sigma})$ sur
des \'equivalence dans $\s -\mathbf{Cat}$, car le changement de base par une \'equivalence
induit une \'equivalence de Quillen entre cat\'egorie de modules. 
On localise alors ce foncteur pour obtenir un $\s$-foncteur
$$\mathcal{M}od(-) : \uscomm \longrightarrow \uscms,$$
qui, \`a \'equivalence pr\`es, envoie $R \in \uscomm$ sur $L^{\otimes}(R-Mod^{cof})=:\mathcal{M}od(R)$
(voir d\'efinition \ref{d9}). On compose alors avec 
le $\s$-foncteur
$$End_{-}(\mathbf{1})\circ \mathcal{M} : \uscms \longrightarrow \usmon$$
puis avec 
$$(-)^{inv} : \usmon \longrightarrow \usmon$$
qui \`a un $\s$-monod\"\i de commutatif $E$ associe son sous-$\s$-monoide des 
\'el\'ements inversibles d\'efini par le produit fibr\'e suivant dans $\mathbf{SEns}^{\Gamma}$
$$\xymatrix{
E^{inv} \ar[r] \ar[d]  & E \ar[d] \\
\pi_{0}(E)^{inv} \ar[r] & \pi_{0}(E).}$$ 
Ceci permet d'obtenir le $\s$-foncteur cherch\'e
$$\mathbb{G}_{m} : \uscomm \longrightarrow \usmon.$$
En termes imag\'es ce $\s$-foncteur envoie $R$ sur le sous-mono\"\i de simplicial
de $\mathbb{R}\underline{End}_{R-Mod}(R,R)$ form\'e des auto\'equivalences. 
La construction pr\'ec\'edente explique
comment ce mono\"\i de est commutatif et comment il peut \^etre rendu fonctoriel en $R$
(il faut noter ici que $R$ n'est pas cofibrant comme $R$-module en g\'en\'eral, voir \cite{shi}). \\

Soit maintenant $(T,\mathcal{A})$ un $\s$-topos annel\'e auquel nous associons les deux
$\s$-topos catannel\'es rigides construits pr\'ec\'edemment
$$(T,\mathbb{G}_{a}(\mathcal{A})) \qquad (T,\mathbb{G}_{m}(\mathcal{A})).$$
On peut ainsi consid\'erer leur pr\'e-trace que nous noterons respectivement
$$d : |\mathbb{G}_{a}(\mathcal{A})| \longrightarrow HH_{pr}^{\mathbb{G}_{a}(\mathcal{A})}$$
$$dlog : |\mathbb{G}_{m}(\mathcal{A})| \longrightarrow HH_{pr}^{\mathbb{G}_{m}(\mathcal{A})}.$$
Pour un objet $X \in T$, ces morphismes \'evalu\'es en $X$ donnent
$$d : |\mathbb{G}_{a}(\mathcal{A})|(X) \longrightarrow 
HH_{pr}^{\mathbb{G}_{a}(\mathcal{A})}(X) \simeq \Omega_{*}|\mathbb{G}_{a}(\mathcal{A})|(X^{S^{1}})$$
$$dlog : |\mathbb{G}_{m}(\mathcal{A})|(X) \longrightarrow 
HH_{pr}^{\mathbb{G}_{m}(\mathcal{A})}(X) \simeq \Omega_{*}|\mathbb{G}_{m}(\mathcal{A})|(X^{S^{1}}),$$
o\`u $\Omega_{*}$ d\'esigne l'espace des lacets pris en \'el\'ement neutre des mono\"\i des correspondants.
Ici $|\mathbb{G}_{a}(\mathcal{A})|(X)$ et $|\mathbb{G}_{m}(\mathcal{A})|(X)$ ne sont autre
que les mono\"\i des simpliciaux additifs et multiplicatifs sous-jacents \`a l'anneau
en spectres $\mathcal{A}(X)$. Nous allons voir maintenant, sans donner les d\'etails cependant, 
que $\Omega_{*}|\mathbb{G}_{a}(\mathcal{A})|(X^{S^{1}})$ et 
$\Omega_{*}|\mathbb{G}_{m}(\mathcal{A})|(X^{S^{1}})$
sont des g\'en\'eralisations de l'espace des formes diff\'erentielles et formes diff\'erentielles logarithmiques, 
ce qui explique les notations $d$ et $dlog$ pour rappeler que ces morphismes sont 
des analogues des diff\'erentielles de de Rham et de Rham logarithmique de fonctions.

Soit $k$ un anneau commutatif.
Consid\'erons $sk-Comm$ la cat\'egorie des $k$-alg\`ebres commutatives simpliciales
(pour bien faire les choses il faut fixer deux univers $\mathbb{U} \in \mathbb{V}$
et ne consid\'erer que les $k$-alg\`ebres commutatives simplicials $\mathbb{U}$-petites. 
La cat\'egorie $sk-Comm$ est alors $\mathbb{V}$-petite). On dispose d'un 
foncteur 
$$sk-Comm \longrightarrow Comm(Sp^{\Sigma})$$
qui envoie une $k$-alg\'ebre commutative simpliciale $R$ sur $HR$, 
le spectre d'Eleinberg-McLane correspondant
(voir \cite{hss}). Ce foncteur pr\'eserve les \'equivalences et induit donc un $\s$-foncteur
$$L(sk-Comm) \longrightarrow \uscomm.$$
Ceci d\'efinit un pr\'echamp sur $dAff_{k}=L(sk-Comm)^{op}$
\`a valeurs dans $\uscomm$. Ce pr\'echamp est de plus un champ pour la topologie
\'etale sur $dAff_{k}$ (voir \cite[\S 2.2]{hagII} par exemple) et fournit donc un $\s$-foncteur
$$\mathcal{A} : (dAff_{k}^{\sim,et})^{op} \longrightarrow \uscomm,$$
qui fait du $\s$-topos $dAff_{k}^{\sim,et}$ un $\s$-topos annel\'e
(le champ $\mathcal{A}$ est pr\'ecis\'ement \`a valeurs dans
$\uscomm_{\mathbb{V}}$, et $dAff_{k}^{\sim,et}$ est un $\mathbb{V}$-$\s$-topos). Notons que 
le $\s$-topos $dAff_{k}^{\sim,et}$ n'est autre que le $\s$-topos des $k$-champs 
d\'eriv\'es au sens de \cite{to3}, et que $\mathcal{A}$ est son faisceau structural 
(repr\'esent\'e par l'objet en anneau $\mathbb{A}^{1}$). Soit maintenant $X$ un 
$k$-sch\'ema, consid\'er\'e comme un objet $X \in dAff_{k}^{\sim,et}$ (voir \cite[\S 2.2.4]{hagII}). 
L'objet $X^{S^{1}}$ est alors \'equivalent au produit fibr\'e $X\times_{X\times X}X$. 
On a donc un isomorphisme naturel
$$\pi_{0}(\Omega_{*}|\mathbb{G}_{a}(\mathcal{A})|(X^{S^{1}}))\simeq
\pi_{1}(|\mathbb{G}_{a}(\mathcal{A})|(X^{S^{1}}))\simeq Tor_{1}^{\mathcal{O}_{X\times X}}(\mathcal{O}_{X},\mathcal{O}_{X})\simeq \Omega^{1}_{X}.$$
Le morphisme construit ci-dessus
$$d : \pi_{0}(|\mathbb{G}_{a}(\mathcal{A})|(X))\simeq \mathcal{O}(X) \longrightarrow
\pi_{0}(\Omega_{*}|\mathbb{G}_{a}(\mathcal{A})|(X^{S^{1}}))\simeq \Omega_{X}^{1}
$$
s'identifie alors \`a la diff\'erentielle de de Rham. De plus, 
il existe des isomorphismes naturels
$$\pi_{i}(|\mathbb{G}_{a}(\mathcal{A})|(X^{S^{1}}))\simeq \pi_{i}(|\mathbb{G}_{m}(\mathcal{A})|(X^{S^{1}}))$$
pour tout $i>0$. A travers ces isomorphismes on peut montrer que le morphisme
$$dlog : \pi_{0}(|\mathbb{G}_{m}(\mathcal{A})|(X))\simeq \mathcal{O}(X)^{*} \longrightarrow
\pi_{0}(\Omega_{*}|\mathbb{G}_{m}(\mathcal{A})|(X^{S^{1}}))\simeq \Omega_{X}^{1}
$$
s'identifie \`a $f \mapsto \frac{df}{f}$, la diff\'erentielle de de Rham logarithmique. 

Notons pour terminer que l'on pourrait aussi bien remplacer le $\s$-topos
$dAff_{k}^{\sim,et}$ par le $\s$-topos des vari\'et\'es $\mathcal{C}^{\infty}$ 
d\'eriv\'ees (voir \cite{sp}), ou encore celui des espaces analytiques complexes
d\'eriv\'es (voir \cite[\S 4.4]{lu3}). Ces deux $\s$-topos sont naturellement annel\'es par
leur faisceaux structuraux respectifs (repr\'esent\'es dans les deux cas par 
la droite affine), et l'interpr\'etation que l'on vient de donner 
des morphsmes $d$ et $dlog$ en termes de diff\'erentielle de de Rham 
reste valable. 

\subsection{Fibr\'es vectoriels et complexes parfaits}

Dans cette section nous d\`ecrivons notre caract\'ere de Chern dans le cas o\`u le $\s$-topos 
$dAff_{k}^{\sim,et}$, dej\`a introduit au num\'ero pr\'ec\'edent, est annel\'e par le champ (en $\s$-CMS 
rigides) des complexes parfaits. En particulier on obtient une caract\'ere de Chern defini pour tous les 
champs d'Artin qui gen\'eralise celui connu pour les sch\'emas et les champs de Deligne-Mumford.  Quand \`a la 
comparaison avec le caract\'ere de Chern classique, nous montrons (Appendice B, Th\'eor\`eme \ref{compfibvect}) 
l'\'equivalence dans le cas des fibr\'es vectoriels sur des variet\'es lisses sur un corps de 
caract\'eristique zero; la comparaison pour les complexes parfaits (qui peut bien se deduire du cas des 
fibr\'es vectoriels) sera preuv\'e dans un travail ult\'erieur.\\

Consid\'erons de nouveau le $\s$-topos $dAff_{k}^{\sim,et}$, des
$k$-sch\'emas d\'eriv\'es. Il peut \^etre muni d'un champ
en $\s$-CMS rigides $\underline{Parf}$, des complexes parfaits, 
de la fa\c{c}on suivante. Notons comme au paragraphe pr\'ec\'edent
$sk-Comm$ la cat\'egorie des $k$-alg\`ebres commutatives simpliciales. 
Pour $R \in sk-Comm$ on dispose d'une $k$-dg-alg\`ebre commutative
normalis\'ee $N(R)$, ainsi que la cat\'egorie mono\"\i dale sym\'etrique
$N(R)-Mod$ des $N(R)$-dg-modules.
On consid\`ere
la sous-cat\'egorie mono\"\i dale pleine $Parf(R)$ form\'es des $N(R)$-dg-modules
qui sont d'une part cofibrants et d'autre part qui sont 
parfaits (c'est \`a dire dualisables dans $\mathrm{Ho}(N(R)-Mod)$, voir \cite{tova2}). 
Pour un morphisme $R \rightarrow R'$ dans
$sk-Comm$ on dispose d'un foncteur de changement de bases
$R'\otimes_{R} - : Parf(R) \longrightarrow Parf(R')$. Ceci 
d\'efinit un foncteur faible 
$$Parf : sk-Comm \longrightarrow \mathbf{Cat}^{\otimes}$$
de $sk-Comm$ vers la $2$-cat\'egorie des cat\'egories mono\"\i dales sym\'etriques. 
De plus, ce foncteur est muni du sous-foncteur form\'e des quasi-isomorphismes
de $N(R)$-dg-modules. Ainsi, en rectifiant ce foncteur faible, puis en 
localisant, on obtient un nouveau foncteur
$$sk-Comm \longrightarrow \scms.$$
Ce foncteur envoie \'equivalences de $k$-alg\`ebres simpliciales sur \'equivalences
de $\uscms$ et donne ainsi lieu \`a un $\s$-foncteur apr\`es localisation
$$\underline{Parf} : L(sk-Comm)=dAff_{k}^{op} \longrightarrow \uscms,$$
qui est tel que $\underline{Parf}(R)$ soit \'equivalent \`a 
$L^{\otimes}(Parf(R))$ (voir d\'efinition \ref{d9}). Le $\s$-foncteur satisfait de plus \`a la condition
de descente pour la topologie \'etale (voir \cite[\S 2.2]{hagII}) et donc fournit un 
champ 
$$\underline{Parf} : dAff_{k}^{\sim,et} \longrightarrow \uscms.$$
On obtient ainsi un $\s$-topos catannel\'e $(dAff_{k}^{\sim,et},\underline{Parf})$. 
Ce $\s$-topos catannel\'e est de plus rigide, car d'une part, par d\'efinition
tout objet de $\underline{Parf}(R)$ est rigide, et d'autre part, 
tout objet $X \in dAff_{k}^{\sim,et}$ est une colimite 
d'objet repr\'esentable. La $\s$-cat\'egorie
$\underline{Parf}(X)$ peut ainsi s'\'ecrire comme une limite (dans $\uscms$)
de $\s$-CMS rigides et est donc elle-m\^eme rigide d'apr\`es le th\'eor\`eme \ref{t1}. \\

\noindent \textbf{Cas des sch\'emas.} Soit maintenant $X$ un $k$-sch\'ema, consid\'er\'e comme
un objet de $dAff_{k}^{\sim,et}$. On a 
$$X^{S^{1}}\simeq X \times_{X\times X} X,$$
o\`u ce produit fibr\'e est pris dans les k-sch\'emas d\'eriv\'es.
Ceci implique en particulier que l'on a 
$$HH^{\underline{Parf}}_{pr}(X)\simeq \mathcal{O}(X^{S^{1}})\simeq
\mathbb{HH}(X),$$
o\`u $\mathbb{HH}(X)$ est un mod\`ele simplicial au complexe
d'homologie de Hochschild de $X$ (obtenu par exemple
par la correspondance de Dold-Kan \`a partir du complexe
de Hochschild de \cite{ke}). En particulier on a
$$\pi_{i}(HH^{\underline{Parf}}_{pr}(X))\simeq HH_{-i}(X)$$
pour tout $i\geq 0$. La pr\'e-trace de la d\'efinition \ref{d20} est donc un morphisme
$$Tr_{X} : [\underline{Parf}(X)]/iso \longrightarrow
HH_{0}(X).$$
Lorsque $X$ est lisse sur $k$ de caract\'eristique nulle on a 
$HH_{0}(X)\simeq \oplus_{p}H^{p}(X,\Omega_{X}^{p})$ et
on peut montrer que, 
pour un fibr\'e vectoriel $V$ sur $X$, 
$Tr_{X}(V)$ n'est autre que le caract\`ere de Chern usuel (\cite[Exp. V]{sga6}) de $V$ \`a valeurs
en cohomologie de Hodge. Nous montrerons ci-dessous un r\'esultat plus fine.\\

Revenons donc au cas g\'en\'eral d'un $k$-sch\'ema $X$. Le groupe simplicial 
$S^{1}$ op\`ere naturellement sur la $k$-alg\`ebre simpliciale
$\mathcal{O}(X^{S^{1}})\simeq \mathbb{HH}(X)$, et on peut montrer que
cette action est donn\'ee par l'op\'erateur de Connes sur le complexe
d'homologie de Hochschild. En particulier, 
on trouve un isomorphisme naturel
$$\pi_{0}(HC_{pr}^{\underline{Parf}}(X))\simeq 
HC_{0}^{-}(X),$$
o\`u $HC_{0}^{-}(X)$ est l'homologie cyclique n\'egative du sch\'ema $X$ d\'efinie par exemple 
dans \cite{ke} (voir \cite[Cor. 4.2 (4)]{rhamloop} pour 
cette comparaison dans le cas lisse et de caract\'eristique nulle). On peut montrer alors
que le caract\`ere de Chern
$$Ch_{X} : [\underline{Parf}(X)]/iso \longrightarrow
HC^{-}_{0}(X)$$
de la d\'efinition \ref{d20} s'identifie au caract\`ere de Chern \`a valeurs dans
la cohomologie cyclique n\'egative (voir par exemple \cite{ke}). 
Nous prouvons ici (Appendice B, Th\'eor\`eme \ref{compfibvect}) ce r\'esultat \emph{dans le cas des fibr\'es vectoriels sur des vari\'et\'es quasi-projectives et lisses sur un corps de caract\'eristique $0$} et reviendrons sur la comparaison g\'enerale dans un travail ult\'erieur. Nous avons choisi de mettre en appendice ce th\'eor\`eme de comparaison, dont la preuve est plut\^ot longue, pour ne pas trop couper le fil de l'exposition. \\

\noindent \textbf{Cas des champs alg\'ebriques.} Supposons maintenant que $X$ soit un $(1-)$champ alg\'ebrique, disons au sens
d'Artin. Le pr\'e-caract\`ere de Chern fournit ainsi un morphisme
$$Ch^{pr} : [\underline{Parf}(R)]/iso \longrightarrow HC^{pr}_{0}(X).$$
Contrairement au cas des sch\'emas il n'est plus vrai que
$HC^{pr}(X)\simeq HC(X)$, et le pr\'e-caract\`ere de Chern devient alors
un rafinement non trivial du caract\`ere de Chern. Par exemple, lorsque 
$X$ est un champ de Deligne-Mumford lisse sur $k$ de caract\'eristique nulle, alors
$HC_{0}^{pr}(X)$ est isomorphe \`a la cohmologie de de Rham (paire) de $IX$, le champ 
d'inertie de $X$, alors que $HC_{0}(X)$ n'est que la cohomologie 
de de Rham de $X$. Le pr\'e-caract\`ere de Chern $Ch^{pr}$ devient alors
le caract\`ere de Chern \`a coefficients dans les repr\'esentations 
de \cite{to}. Dans le cas g\'en\'eral d'un champ d'Artin $X$ sur un anneau de base quelconque $k$, 
$HC_{*}^{pr}(X)$ forme une th\'eorie cohomologique digne d'int\'er\^et. Le caract\`ere de
Chern 
$$Ch^{pr} : [\underline{Parf}(X)]/iso \longrightarrow HC_{0}^{pr}(X)$$
est alors une g\'en\'eralisation du caract\`ere de Chern construit dans \cite{to} au cadre
beaucoup plus g\'en\'eral des champs d'Artin. Notons que lorsque $X=BG$ est le champ
classifiant d'un sch\'ema en groupe lisse sur $k$ alors ce caract\`ere de Chern
$$Ch^{pr} : [\underline{Parf}(BG)]/iso \longrightarrow HC_{0}^{pr}(BG)\simeq \mathcal{O}(G)^{G},$$
n'est autre que le morphisme de trace, qui \`a un complexe parfait $E$ muni d'une
action de $G$ associe la fonction $Tr_{E}$, qui envoie $g\in G$ sur la trace
de l'endomorphisme $g$ de $E$. Dans le cas plus g\'en\'eral d'un champ
quotient $[X/G]$, $Ch^{pr}$ est une combinaison relativement subtile 
du caract\`ere de Chern pour le sch\'ema des points fixes $X^{fix}=\{(g,x)/g.x=x\} \subset G\times X$ et 
du morphisme de trace pour $G$. Nous espr\'erons revenir sur 
une description plus pr\'ecise de la th\'eorie cohomologique $X \mapsto HC_{0}^{pr}(X)$ 
pour les champs d'Artin, ainsi que du caract\`ere de Chern qui lui est associ\'e, dans
un travail ult\'erieur.  \\

Tout comme dans notre exemple pr\'ec\'edent nous aurions pu remplacer
le $\s$-topos $dAff^{\sim,et}$ par le $\s$-topos des vari\'et\'es
diff\'erentielles ou analytiques d\'eriv\'ees, munis de leurs champs
des complexes parfaits. La comparaison entre notre caract\`ere de Chern et 
le caract\`ere de Chern usuel des fibr\'es $\mathcal{C}^{\infty}$ ou holomorphes \`a valeurs
dans la cohomologie de de Rham ($\mathcal{C}^{\infty}$ ou holomorphe) resterait alors
valable. Nous reviendrons aussi sur ce point dans un travail ult\'erieur.

\subsection{Familles de dg-cat\'egories compactement engendr\'ees}

Comme dans le paragraphe pr\'ec\'edent nous consid\'erons
$dAff_{k}^{\sim,et}$, le $\s$-topos des $k$-sch\'emas d\'eriv\'es. 
Pour $R \in sk-Comm$ une $k$-alg\`ebre commutative simpliciale nous
consid\'erons $N(R)$ la dg-alg\`ebre normalis\'ee associ\'ee. On 
note alors $dg-cat_{R}$ la cat\'egorie des cat\'egories enrichies
dans $N(R)-Mod$, munie de sa structure de mod\`eles
de \cite{tab}. Pour pr\'eciser les univers rappelons que 
$sk-Comm$ d\'esigne la cat\'egorie des $k$-alg\`ebres simpliciales $\mathbb{U}$-petites, et 
de m\^eme $N(R)-Mod$ d\'esigne les $N(R)$-dg-modules $\mathbb{U}$-petits. 
Soit $\mathbb{V}$ un univers avec $\mathbb{U} \in \mathbb{V}$. Alors, par d\'efinition, 
$dg-cat_{R}$ consistera pr\'ecis\'ement en les $N(R)-Mod$-cat\'egories
qui sont $\mathbb{V}$-petites. Pour tout $R$ soit 
$dg-cat_{R}^{c}$ la cat\'egorie des objets cofibrants dans $dg-cat_{R}$. Le
produit tensoriel au-dessus de $N(R)$ fait de
$dg-cat_{R}^{c}$ une cat\'egorie sym\'etrique mono\"\i dale, munie d'une notion
d'\'equivalences qui est compatible. Ainsi, en rectifiant le pseudo-foncteur
$$R \mapsto dg-cat_{R}^{c} \qquad (R\rightarrow R') \mapsto N(R')\otimes_{N(R)}- ,$$
puis en utilisant la version mono\"\i dale de la localisation $L^{\otimes}$ 
(voir d\'efinition \ref{d9}),
 on trouve un 
foncteur
$$D : sk-Comm \longrightarrow \scms.$$
Pour bien faire les choses ici
il nous faut fixer 
$\mathbb{W}$ un troisi\`eme univers
avec $\mathbb{V} \in \mathbb{W}$. On dispose alors pr\'ecis\'ement d'un foncteur
$$D : dAff_{k}^{op} \longrightarrow \scms_{\mathbb{W}}.$$
Ce foncteur est tel que 
$D(R)$ est \'equivalente \`a 
$L^{\otimes}(dg-cat_{R}^{c})$, la $\s$-CMS  des $R$-dg-cat\'egories qui sont $\mathbb{V}$-petites. 
Ce foncteur $D$ peut aussi \'etre consid\'er\'e comme un foncteur
$$D : sk-Comm \times \Gamma \longrightarrow \s -\mathbf{Cat}.$$
Nous utilisons maintenant la proposition  \ref{pcart} pour obtenir une $\s$-cat\'egorie cofibr\'ee 
$$\int_{sk-Comm \times \Gamma}D \longrightarrow sk-Comm \times \Gamma.$$

Nous allons maintenant d\'efinir une
sous-$\s$-cat\'egorie cofibr\'ee 
de $\int_{sk-Comm\times \Gamma}D$. 
Pour cela, notons que $[\int_{sk-Comm\times \Gamma}D]$ se d\'ecrit de la fa\c{c}on suivante. 
Ses objets sont des triplets $(R,[n],A)$, avec $R\in sk-Comm$, 
$[n] \in \Gamma$, et 
$A=(A_{1}, \dots, A_{n})$ une famille de dg-cat\'egories sur $N(R)$ qui sont
$\mathbb{V}$-petites. Un morphisme $(R,[n],A) \longrightarrow (S,[m],B)$ est la donn\'ee 
d'un morphisme $R \rightarrow S$, d'un morphisme $u : [n] \rightarrow [m]$, et pour tout
$j \in \{1, \dots, n\}$
d'un morphisme dans $\mathrm{Ho}(dg-cat_{S})$
$$S\otimes_{R}^{\mathbb{L}}\left(\otimes^{\mathbb{L}}_{i\in u^{-1}(j)}A_{i}\right) \longrightarrow B_{j}.$$
Rappelons alors qu'une dg-cat\'egorie $A$ sur $N(R)$ est 
$\mathbb{U}$-compactement engendr\'ee si elle est isomorphe dans $\mathrm{Ho}(dg-cat_{R})$
\`a une dg-cat\'egorie de la forme 
$$\widehat{A_{0}}=\rh (A_{0}^{op},\widehat{\mathbf{1}}),$$ 
avec $A_{0}$ une dg-cat\'egorie $\mathbb{U}$-petite (ou encore 
\`a la dg-cat\'egorie des $\mathbb{U}$-dg-modules cofibrants sur une dg-cat\'egorie $\mathbb{U}$-petite, voir
\cite{to1}). Nous d\'efinissions une sous-cat\'egorie $C$, non pleine,  de $[\int_{sk-Comm\times \Gamma}D]$
de la fa\c{c}on suivante.

\begin{enumerate}
\item Les objets de $C$ sont les $(R,[n],A)$ tels que chaque $A_{i}$ soit
$\mathbb{U}$-compactement engendr\'ee.
\item Les morphismes $(R,[n],A) \rightarrow (S,[m],B)$ dans $C$ sont 
ceux pour les quels pour tout $1\leq j\leq m$ le morphisme de $N(R)$-dg-cat\'egories
$$\otimes^{\mathbb{L}}_{i\in u^{-1}(j)}A_{i} \longrightarrow B_{j}$$
soit \emph{multi-continu}, c'est \`a dire commute aux sommes dans chacune
des variables $A_{i}$ ind\'ependamment.
\end{enumerate}

On d\'efinit alors la sous-$\s$-cat\'egorie $\mathcal{C} \subset \int_{sk-Comm\times \Gamma}D$
par le carr\'e (homotopiquement) cart\'esien suivant
$$\xymatrix{
\mathcal{C} \ar[r] \ar[d] & \int_{sk-Comm\times \Gamma}D \ar[d] \\
C \ar[r] & [\int_{sk-Comm\times \Gamma}D].}$$
Le point crucial ici est que la sous-$\s$-cat\'egorie $\mathcal{C}$
est encore cofibr\'ee au sens de la d\'efinition \ref{dcart}. En effet, en d\'eroulant les d\'efinitions
cela se d\'eduit des r\'esultats de \cite{to1} de la fa\c{c}on suivante. 
Soit $(R,[n],A)$ un objet de $C$ et $R \rightarrow S$, $u : [n] \rightarrow [m]$
un morphisme dans $sk-Comm \times \Gamma$. Choisissons
des dg-cat\'egories $\mathbb{U}$-petites $A_{i,0}$ avec 
$A_{i}\simeq \widehat{A_{i,0}}$. On consid\`ere, pour tout $1\leq j\leq m$
le morphisme naturel
$$\otimes^{\mathbb{L}}_{i\in u^{-1}(j)}A_{i} \longrightarrow 
\widehat{\otimes^{\mathbb{L}}_{i\in u^{-1}(j)}A_{i,0}},$$
induit par le plongement de Yoneda restreint \`a 
$$\otimes^{\mathbb{L}}_{i\in u^{-1}(j)}A_{i,0} \hookrightarrow \otimes^{\mathbb{L}}_{i\in u^{-1}(j)}A_{i}.$$
Ainsi, si l'on pose 
$$B_{j}:=\widehat{\otimes^{\mathbb{L}}_{i\in u^{-1}(j)}A_{i,0}},$$
on dispose d'un morphisme $(R,[n],A) \rightarrow (S,[m],B)$, que l'on voit facilement 
\^etre dans la sous-$\s$-cat\'egorie $C$. De plus, le th\'eor\`eme \cite[Thm. 7.2]{to1} implique 
qu'un rel\`evement de ce morphisme en un morphisme 
de $\mathcal{C}$ est cocart\'esien dans $\mathcal{C}$. Ceci termine de montrer que 
$\mathcal{C}$ est fibr\'ee. 

Finalement, on utilise de nouveau la proposition \ref{pcart} en appliquant cette fois le foncteur
$\mathbb{R}Se_{sk-Comm\times \Gamma}$ pour obtenir un nouveau foncteur
$$sk-Comm \times \Gamma \longrightarrow \s -\mathbf{Cat}_{\mathbb{W}},$$
o\`u encore un foncteur
$$sk-Comm \longrightarrow \s -\mathbf{Cat}_{\mathbb{W}}^{\Gamma}.$$
Il n'est pas difficile de voir que ce foncteur se factorise en
$$sk-Comm \longrightarrow \scms_{\mathbb{W}}.$$
Il envoie \'equivalences sur \'equivalences et donc induit un $\s$-foncteur 
$$\mathbb{D}g : dAff_{k} \longrightarrow \uscms_{\mathbb{W}}.$$
Pour $R \in sk-Comm$ la cat\'egorie
sym\'etrique mono\"\i dale $[\mathbb{D}g(R)]$ est naturellement \'equivalente
\`a $\mathrm{Ho}(dg-cat_{R})^{ct}$, la sous-cat\'egorie des dg-cat\'egories
$\mathbb{U}$-compactement engendr\'ees et dg-foncteurs continus. Sa structure
mono\"\i dale $\otimes^{ct}$ est donn\'ee par la formule
$$\widehat{A_{0}} \otimes^{ct} \widehat{B_{0}} \simeq \widehat{A_{0}\otimes^{\mathbb{L}}B_{0}}.$$
De m\^eme, le foncteur de changement de bases pour $R \rightarrow S$ est donn\'e par 
$$S\otimes^{ct}_{R} \widehat{A_{0}} \simeq \widehat{S\otimes^{\mathbb{L}}_{R}A_{0}}.$$
Le th\'eor\`eme \cite[Thm. 7.2]{to1} implique alors que $[\mathbb{D}g(R)]$ est rigide pour
tout $R$: en effet, le dual de $\widehat{A_{0}}$ n'est autre
que $\widehat{A_{0}^{op}}$, et le morphisme d'unit\'e 
$$\widehat{\mathbf{1}} \longrightarrow \widehat{A_{0}\otimes^{\mathbb{L}}A_{0}^{op}}$$
n'est autre que le $A_{0}$-bi-dg-module $(a,b) \mapsto A(b,a)$.

De plus, d'apr\`es  \cite[Thm. 0.2]{taz}, 
une forme tordue de dg-cat\'egories compactement engendr\'ees est encore
compactement engendr\'ee, ce qui implique que le pr\'echamp 
$\mathbb{D}g$ est en r\'ealit\'e un champ pour la topologie \'etale (et 
m\^eme fppf). Comme chaque $\mathbb{D}g(R)$ est 
une $\s$-CMS rigide, nous venons ainsi de construire 
un $\s$-topos catannel\'e rigide $(dAff_{k}^{\sim,et},\mathbb{D}g)$. \\

Appliquons la construction du pr\'e-caract\`ere de Chern de la d\'efinition \ref{d20} au 
$\s$-topos catannel\'e rigide $(dAff_{k}^{\sim,et},\mathbb{D}g)$. 
Il s'agit, pour $X \in dAff_{k}^{\sim,et}$ un $k$-champ d\'eriv\'e, d'une
application
$$Ch^{pr} : [\mathbb{D}g(X)]/iso \longrightarrow 
HC^{pr}_{0}(X).$$
La cat\'egorie $[\mathbb{D}g(X)]$ est la cat\'egorie homotopique
des champs en dg-cat\'egories compactement engendr\'ees parametr\'es par $X$. 
D'autre part, le $\s$-foncteur $End_{\mathbb{D}g}(\mathbf{1})$ s'identifie 
au pr\'e-champ des complexes quasi-coh\'erents sur $dAff_{k}^{\sim,et}$, et nous
savons que ce pr\'echamp est un champ (voir  \cite[\S 2.2]{hagII}). Ainsi, on a 
$$HC^{pr}_{0}(X) \simeq \pi_{0}(End_{\mathbb{D}g}(\mathbf{1})(X^{S^{1}})^{hS^{1}}) \simeq
D_{qcoh}^{S^{1}}(X^{S^{1}})/iso,$$
o\`u $D_{qcoh}^{S^{1}}(X^{S^{1}})$ est la cat\'egorie d\'eriv\'ee des 
complexes quasi-coh\'erents $S^{1}$-\'equivariants sur $X^{S^{1}}$
(c'est \`a dire des complexes quasi-coh\'erents sur le champ
quotient $[X^{S^{1}}/S^{1}]$). Ainsi, le pr\'e-caract\`ere de Chern
d'un champ $A$ en dg-cat\'egories compactement engendr\'ees sur $X$ est 
un complexe quasi-coh\'erent $S^{1}$-\'equivariant
$$Ch^{pr}(A) \in D_{qcoh}^{S^{1}}(X^{S^{1}})/iso.$$
La fibre de $Ch^{pr}(A)$ en un point 
$\gamma : S^{1} \longrightarrow X$, est l'homologie de
Hochschild de la fibre $A_{\gamma(*)}$ \`a coefficient dans l'endomorphisme de
monodromie de $\gamma$ (qui est une auto\'equivalence de 
$A_{\gamma(*)}$). Le (complexe de) faisceau $Ch^{pr}(A)$ peut \^etre vu
comme le \emph{faisceau d'anomalies de $A$}, au sens o\`u ce terme
est utilis\'e dans \cite[\S 6.2]{bry}. \\

Nous allons terminer ce parapraphe par deux cas particulier 
du pr\'e-caract\`ere de Chern ci-dessus. \\

\textbf{Connexion de Gauss-Manin non-commutative.} Supposons que $X$ soit un sch\'ema
lisse sur $k$ de caract\'eristique nulle. Soit $\mathcal{D}_{X}$ le faisceau
des op\'erateurs diff\'erentiels sur $X$. On construit un faisceau de $k$-dg-alg\`ebres
sur $X$ en posant
$$\mathcal{R}_{X}:=\left(\oplus_{i\in \mathbb{N}}\mathcal{D}_{X}^{\leq i}[-2i]\right),$$
o\`u $\mathcal{D}_{X}^{\leq i}$ est le sous-faisceau des op\'erateurs diff\'erentiels
de degr\'e inf\'erieur \`a $i$. C'est un faisceau de $k[u]$-dg-alg\`ebres sur $X$, o\`u 
$u$ une variable plac\'ee en degr\'e $2$ et op\'erant par les inclusions naturelles
$$u. : \mathcal{D}_{X}^{\leq i} \hookrightarrow \mathcal{D}_{X}^{\leq i+1}.$$
On peut v\'erifier (voir par exemple \cite{bena}) que le faisceau 
$\mathcal{R}_{X}$ est la dg-alg\`ebre des endomorphismes
de l'objet $\mathcal{O}_{X}$ vu dans $D_{qcoh}^{S^{1}}(X^{S^{1}})$ par 
le foncteur d'image directes le long de l'inclusion naturelle $X \hookrightarrow X^{S^{1}}$.
Ainsi, 
\`a l'aide du foncteur $E \mapsto \rh (\mathcal{O}_{X},E)$ on construit un foncteur
$$\phi : D_{qcoh}^{S^{1}}(X^{S^{1}}) \longrightarrow D_{qcoh}(\mathcal{R}_{X}).$$
Ce foncteur, restreint aux objets satisfaisant certaines conditions de finitudes convenables est de plus
une \'equivalence de cat\'egories, mais nous nous contenterons ici de son existence.

Soit maintenant $A \in \mathbb{D}g(X)$ un
champ en dg-cat\'egories compactement engendr\'ees sur $X$. Son
pr\'e-caract\`ere de Chern $Ch^{pr}(A)$ est un objet dans $D_{qcoh}^{S^{1}}(X^{S^{1}})$, qui, 
\`a l'aide du foncteur $\phi$, peut \^etre vu comme un objet
$$\phi(Ch^{pr}(A)) \in D_{qcoh}(\mathcal{R}_{X}).$$
Cet objet donne lieu d'une part \`a un $\mathcal{D}_{X}$-module
$\mathbb{Z}/2$-gradu\'e
$$GM(A):=\phi(Ch^{pr}(A))[u^{-1}] \in D_{qcoh}(\mathcal{R}_{X}[u^{-1}])\simeq
D_{qcoh}^{\mathbb{Z}/2}(\mathcal{D}_{X}),$$
et d'autre part \`a un \emph{gradu\'e associ\'e}
$$Char(A):=\phi(Ch^{pr}(A)[u=0]) \in D_{qcoh}(\mathcal{R}_{X}\otimes_{k[u]}^{\mathbb{L}}k)\simeq
D_{qcoh}(Sym_{\mathcal{O}_{X}}(\mathbb{T}_{X}[-2])).$$
L'objet $GM(A)$ n'est autre que le complexe d'homologie p\'eriodique de $A$ relativement 
\`a $X$, muni d'une connexion plate: la connexion de Gauss-Manin. Nous ne comparerons pas ici cette
connexion avec les connexions de Gauss-Manin apparaissant dans d'autres contextes. D'autre part, 
$Char(A)$ est le \emph{module caract\'eristique} du $\mathcal{D}_{X}$-module $GM(A)$, associ\'e
\`a une certaine filtration: la filtration de Hodge. De m\^eme, nous ne comparerons pas ici cette
filtration avec les connexions de Gauss-Manin apparaissant dans d'autres contextes. 
L'objet global $\phi(Ch^{pr}(A))$ contient ainsi la connexion de Gauss-Manin et 
la filtration de Hodge, toutes deux d\'efinies sur le complexe d'homologie p\'eriodique
de $A$. Il faut donc comprendre l'objet $\phi(Ch^{pr}(A))$, et donc  l'objet
$Ch^{pr}(A)$ car $\phi$ n'est pas bien loin d'\^etre une \'equivalence de cat\'egorie, 
comme la \emph{partie alg\'ebrique de la variation de structures de Hodge associ\'ee
\`a $A$, vu comme une famille de vari\'et\'es non-commutatives param\'etr\'ee par $X$}
(au sens de \cite{kakopa}). \\
 
\textbf{Faisceaux des caract\`eres d'une repr\'esentation dg-cat\'egorique.} 
Soit $G$ un groupe alg\'ebrique sur un corps $k$ et $X=BG$. Soit $A$ un
champ en dg-cat\'egories compactement engendr\'ees sur $X$. Un tel objet 
est une \emph{repr\'esentation dg-cat\'egorique de $G$}. On peut en construire 
en prenant par exemple une dg-cat\'egorie $T$ dont le groupe des auto\'equivalences 
forme un groupe alg\'ebrique (par exemple lorsque $T$ est satur\'ee, voir \cite[Cor. 3.26]{tova2}), et 
en consid\'erant une repr\'esentation de $G$ dans ce groupe. De telles repr\'esentations
cat\'egoriques apparaissent aussi dans le contexte de la correspondance de Langlands
g\'eom\'etrique (voir par exemple \cite[\S 21]{fregai}). Le pr\'e-caract\`ere de Chern pour $A$ fournit 
alors un complexe quasi-coh\'erent  sur $G$, qui est 
d'une part \'equivariant pour l'action de $G$ sur lui-m\^eme par conjugaison, et d'autre
part muni d'une action compatible de $S^{1}$ 
$$Ch^{pr}(A) \in D_{qcoh}^{S^{1}}([G/G])
$$
(remarquer que $S^{1}$ op\`ere
naturellement sur $X^{S^{1}}\simeq [G/G]$). Par d\'efinition, 
ce complexe quasi-coh\'erent est le \emph{(complexe de) faisceau des caract\`eres
de la repr\'esentation de $G$}. La restriction de ce complexe dans un voisinage 
formel de $e\in G$ fournit de nouveau un complexe de $\mathcal{R}_{X}$-modules, 
de nouveau correspondant \`a la variation de structures de Hodge
associ\'ee \`a $A$ (voir \cite[Cor. 5.6]{bena}). Cependant, l'objet $Ch^{pr}(A)$ contient plus
d'informations que cette variation car il n'est plus vrai que 
le foncteur $D_{qcoh}^{S^{1}}(X^{S^{1}}) \longrightarrow D_{qcoh}(\mathcal{R}_{X})$ est 
une \'equivalence lorsque $X$ est un champ (m\^eme en imposant 
des conditions de finitudes). Par exemple, lorsque $G$ est fini, 
le complexe $Ch^{pr}(A)$ se souvient d'une information non triviale de la repr\'esentation 
$G$ alors que le $\mathcal{R}_{X}$-module correspondant ne voit que 
l'action de $G$ induite sur l'homologie cylique de $A$. La situation 
est en r\'ealit\'e tout \`a fait similaire \`a celle pour le 
caract\`ere de Chern des complexes parfaits que nous avons
d\'ecrite dans le paragraphe \S 4.2. 

\subsection{Familles de dg-cat\'egories satur\'ees et homologie cyclique secondaire}

Nous terminerons cette section par un rafinement des constructions
du paragraphe pr\'ec\'edent obtenu en rajoutant une condition de finitude
sur les dg-cat\'egories consid\'er\'ees. Pour cela, nous revenons
au $\s$-topos catannel\'e $(dAff_{k}^{\sim,et},\mathbb{D}g)$. 
Pour tout $R\in sk-Comm$ la cat\'egorie
$[\mathbb{D}g(R)]$ s'identifie \`a la cat\'egorie homotopique des 
dg-cat\'egories sur $R$, compactement engendr\'ees, et des morphismes
continus. Nous dirons qu'un morphisme dans $[\mathbb{D}g(R)]$ est 
\emph{compact}, s'il correspond \`a un dg-foncteur 
$A \longrightarrow B$ dont le foncteur induit
$[f] : [A] \longrightarrow [B]$ pr\'eserve les objets compacts
(rappelons ici que les cat\'egories $[A]$ et $[B]$ sont 
naturellement triangul\'ees, et qu'un objet $x$ d'une cat\'egorie triangul\'ee
est compact si $[x,-]$ commute aux sommes). De mani\`ere 
\'equivalente, $f$ est compact si et seulement si son adjoint \`a droite
(au sens dg-cat\'egorique) est un dg-foncteur continu. On dispose ainsi d'une
sous-cat\'egorie, non pleine, $[\mathbb{D}g^{c}(R)] \subset [\mathbb{D}g(R)]$ 
form\'ee des morphismes compacts. Il n'est pas difficile de voir que ces 
sous-cat\'egories sont stables par la structure mono\"\i dale $\otimes^{ct}$, ainsi 
que par les changement de bases par des morphismes $R \rightarrow S$ 
dans $sk-Comm$. Ainsi, si pour $R \in sk-Comm$, on d\'efinit
$\mathbb{D}g^{c}(R) \subset \mathbb{D}g(R)$ par la carr\'e cart\'esien suivant
$$\xymatrix{
  \mathbb{D}g^{c}(R) \ar[r] \ar[d] & \mathbb{D}g(R) \ar[d] \\
  [\mathbb{D}g^{c}(R)] \ar[r] & [\mathbb{D}g(R)],}$$
alors $R \mapsto \mathbb{D}g^{c}(R)$ d\'efinit un sous-$\s$-foncteur 
de $\mathbb{D}g(R)$. Cependant, pour une dg-cat\'egorie $A \in [\mathbb{D}g(R)]$, 
les morphismes d'unit\'e et de co-unit\'e
$$\widehat{\mathbf{1}} \longrightarrow A \otimes^{ct}A^{\vee} \qquad
A \otimes^{ct}A^{\vee} \longrightarrow \widehat{\mathbf{1}}$$
ne sont pas des morphismes compacts en g\'en\'eral. Cela implique que
les objets de $\mathbb{D}g^{c}(R)$ ne sont plus rigides en g\'en\'eral. 

Par d\'efinition, une dg-cat\'egorie $A \in \mathbb{D}g(R)$ sera dite
\emph{satur\'ee} si elle est rigide en tant qu'objet de la $\s$-CMS
$\mathbb{D}g^{c}(R)$. Nous noterons $\mathbb{D}g^{sat}(R)$ la sous-$\s$-CMS
pleine form\'ee des objets rigides dans $\mathbb{D}g^{c}(R)$. Cela d\'efinit un nouveau
$\s$-foncteur
$$\mathbb{D}g^{sat} : (dAff_{k}^{\sim,et})^{op} \longrightarrow \uscms_{rig}.$$
Le pr\'e-champ $\mathbb{D}g^{sat}$ est un sous-champ 
$\mathbb{D}g$, car on v\'erifie qu'\^etre satur\'e est une condition
locale pour la topologie \'etale (voir \cite{taz} pour plus de d\'etails). On dipose ainsi 
d'un $\s$-topos catannel\'e rigide $(dAff_{k}^{\sim,et},\mathbb{D}g^{sat})$. 

Soit maintenant $X$ un champ alg\'ebrique au sens d'Artin. On a
$$HC_{0}^{pr}(X)\simeq D_{parf}^{S^{1}}(X^{S^{1}})/iso,$$
o\`u $D_{parf}^{S^{1}}(X^{S^{1}})$ est la sous-cat\'egorie pleine de
$D_{qcoh}^{S^{1}}(X^{S^{1}})$ form\'ee des objets dont le complexe quasi-coh\'erent
sous-jacent est parfait sur $X^{S^{1}}$. Le pr\'e-caract\`ere de Chern
induit ainsi une application
$$Ch^{pr} : [\mathbb{D}g^{sat}(X)]/iso \longrightarrow
D_{parf}^{S^{1}}(X^{S^{1}})/iso.$$
Cette application est de plus compatible avec le pr\'e-caract\`ere de Chern 
du $\s$-topos catannel\'e $(dAff_{k}^{\sim,et},\mathbb{D}g)$, ce
qui fournit un diagramme commutatif
$$\xymatrix{
[\mathbb{D}g^{sat}(X)]/iso \ar[r] \ar[d] & D_{parf}^{S^{1}}(X^{S^{1}})/iso \ar[d] \\
[\mathbb{D}g(X)]/iso \ar[r] & D_{qcoh}^{S^{1}}(X^{S^{1}})/iso.}$$
On d\'eduit de cela et des consid\'erations sur l'existence de la connexion 
de Gauss-Manin sur l'homologie p\'eriodique du paragraphe pr\'ec\'edent le fait 
important suivant.

\begin{cor}\label{c5}
Soit $X$ un sch\'ema lisse sur $k$ de caract\'eristique nulle et 
$A \in [\mathbb{D}g^{sat}(X)]$ une famille de dg-cat\'egories
satur\'ees param\'etr\'ee par $X$. Alors les faisceaux 
d'homologie p\'eriodique $HP_{i}(A)$ sont des fibr\'es vectoriels sur $X$.
\end{cor}

Revenons au cas o\`u $X$ est un champ alg\'ebrique sur $k$ quelconque et au morphisme
$$Ch^{pr} : 
[\mathbb{D}g^{sat}(X)]/iso \longrightarrow D_{parf}^{S^{1}}(X^{S^{1}})/iso.$$
Soit $A \in [\mathbb{D}g^{sat}(X)]$. Le complexe
parfait $S^{1}$-\'equivariant $Ch^{pr}(A)$  poss\`ede lui-m\^eme un 
pr\'e-caract\`ere de Chern d\'ecrit dans le paragraphe \S 4.2, qui est un \'el\'ement dans
$$Ch^{pr}(Ch^{pr}(A)) \in \pi_{0}(\mathcal{O}((X^{S^{1}})^{S^{1}})^{hS^{1}}).$$
Or, 
$(X^{S^{1}})^{S^{1}}\simeq X^{S^{1}\times S^{1}}$, et ainsi
le fait que $Ch^{pr}(A)$ soit $S^{1}$-\'equivariant implique que 
$Ch^{pr}(Ch^{pr}(A))$ est un \'el\'ement dans 
$$\pi_{0}(\mathcal{O}(X^{S^{1} \times S^{1}})^{hS^{1}\times S^{1}}).$$
Ceci d\'efinit une nouvelle application
$$Ch^{pr,(2)} : [\mathbb{D}g^{sat}(X)]/iso \longrightarrow
\pi_{0}(\mathcal{O}(X^{S^{1}\times S^{1}})^{hS^{1}\times S^{1}}).$$
Le membre de droite de cette application est, par d\'efinition, la pr\'e-homologie cylique
secondaire de $X$ et sera not\'ee
$$HC_{i}^{pr,(2)}(X):=\pi_{i}(\mathcal{O}(X^{S^{1} \times S^{1}})^{hS^{1}\times S^{1}}).$$
Il faut penser \`a $HC_{i}^{pr,(2)}(X)$ comme \`a \emph{l'homologie cyclique
$S^{1}$-\'equivariante de $X^{S^{1}}$}. Lorsque 
$X$ est un sch\'ema lisse sur $k$ de caract\'eristique nulle nous pensons que
$HC_{i}^{pr,(2)}(X)$ poss\`ede une description en termes de
\emph{complexe de de Rham secondaire}. Nous esp\'erons ainsi pouvoir 
d\'ecrire, au moins partiellement, $Ch^{pr,(2)}(A)$ en termes
de la variation de structures de Hodge $HP_{*}(A)$ induite sur 
$X$. Nous reviendrons sur cette question dans un travail ult\'erieur. \\ 

\begin{appendix}

\section{Sur les cat\'egories de simplexes}

Dans cet appendice nous revenons sur deux r\'esultats concernant
la cat\'egorie de simplexes $\Delta(I)$ d'une cat\'egorie $I$ donn\'ee. 
Ces deux r\'esultats sont utilis\'es lors de la preuve du lemme
\ref{lcart}. \\

Soit donc $I$ une cat\'egorie et notons $\Delta(I)$ sa cat\'egorie
des simplexes. Les objets de $\Delta(I)$ sont les couples
$([n],u)$, form\'es d'un objet $[n]\in \Delta$ et d'un
foncteur $u : \Delta^{m} \longrightarrow I$. Un morphisme
$([n],u) \rightarrow ([m],v)$ est la donn\'ee d'un morphisme
$f : [m] \rightarrow [n]$ dans $\Delta$ tel que le diagramme suivant
$$\xymatrix{
\Delta^{n} \ar[r]^-{f} \ar[d]_-{u} & \Delta^{m} \ar[dl]_-{v} \\
I & }$$
commute (strictement). On peut aussi \'ecrire
$$\Delta(I):=\int_{\Delta^{op}}N(I),$$
o\`u $N(I) : \Delta^{op} \longrightarrow Ens$
est le nerf de $C$.

On dispose d'une projection naturelle
$$\pi : \Delta(I) \longrightarrow I$$
d\'efinie sur les objets par $\pi([n],u):=u(0)$. Pour un morphisme
$([n],u) \rightarrow ([m],v)$ comme ci-dessus, on dispose d'un morphisme
naturel 
$v(0) \rightarrow v(f(0))=u(0)$, ce qui d\'efinit $\pi$ sur les morphismes. 
Nous dirons alors qu'un morphisme $([n],u) \rightarrow ([m],v)$ est 
\emph{vertical} si son image par $\pi$ est une identit\'e (ou 
de mani\`ere \'equivalente si le morphisme $[m] \rightarrow [n]$
envoie $0$ sur $0$).

\begin{prop}\label{pa1}
Soit $W$ l'ensemble des morphismes verticaux de $\Delta(I)$. Alors le
$\s$-foncteur induit par $\pi$
$$p : L_{W}\Delta(I) \longrightarrow I$$
est une \'equivalence. 
\end{prop}

\textbf{Preuve --}  Le foncteur $\pi$ \'etant essentiellement surjectif il en est
de m\^eme de $p : L_{W}\Delta(I) \longrightarrow I$. Pour montrer qu'il est 
pleinement fid\`ele, il suffit de montrer que le $\s$-foncteur induit
$$\mathbb{L}p_{!} : \widehat{L_{W}\Delta(I)} \longrightarrow \widehat{I}$$
est pleinement fid\`ele. Cela est aussi \'equivalent au fait que le $\s$-foncteur
$$\pi^{*} : \widehat{I} \longrightarrow \widehat{\Delta(I)}$$
est pleinement fid\`ele et que son image essentielle consiste en 
les $\s$-foncteurs $F : \Delta(I)^{op} \longrightarrow \T$ 
qui envoient les morphismes verticaux sur des \'equivalences. En traduisant 
cela en termes de cat\'egories de mod\`eles on voit qu'il faut montrer que l'adjonction 
de Quillen induite sur les cat\'egorie de pr\'efaisceaux simpliciaux
$$\pi_{!} : SPr(\Delta(I)) \longleftrightarrow SPr(I) : \pi^{*}$$
induit un foncteur pleinement fid\`ele
$$\pi^{*} :  \mathrm{Ho}(SPr(I)) \longrightarrow \mathrm{Ho}(SPr(\Delta(I)))$$
d'image essentielle consistant en les foncteurs qui envoient les
morphismes verticaux sur des \'equivalences. Comme le foncteur 
$\pi$ est surjectif sur les ensembles d'objets le foncteur $\pi^{*}$ est 
conservatif. Ainsi, il nous suffit de montrer que pour 
tout foncteur $F : \Delta(I)^{op} \longrightarrow \mathbf{SEns}$, qui envoie
morphismes verticaux sur \'equivalences, le morphisme d'adjonction
$$\pi^{*}\mathbb{L}\pi_{!}(F) \longrightarrow F$$
est un isomorphisme dans $\mathrm{Ho}(SPr(\Delta(I)))$.

Soit alors $i\in I$, et consid\'erons le diagramme
commutatif
$$\xymatrix{
\Delta(I) \ar[r]^-{\pi} & I \\
\pi^{-1}(i) \ar[u]^-{v} \ar[r]_-{q} & \{i\}=*. \ar[u]_-{u}}$$
Ce diagramme induit une transformation naturelle de changement de base
$$\mathbb{L}q_{!}v^{*} \Rightarrow u^{*}\mathbb{L}\pi_{!}.$$
Cette transformation naturelle est en  r\'ealit\'e un isomorphisme. En effet, 
comme tous les foncteurs en question commutent aux colimites homotopiques, et 
que tout pr\'efaisceau simplicial est obtenu par colimites homotopiques \`a partir de pr\'efaisceaux
repr\'esentables, il nous faut montrer que pour un objet $j\in \Delta(I)$
le morphisme naturel
$$\mathbb{L}q_{!}v^{*}(h_{j}) \longrightarrow u^{*}\mathbb{L}\pi_{!}(h_{j})$$
est un isomorphisme (rappelons que $h_{j}:=Hom(-,j)$ est le pr\'efaisceau
repr\'esent\'e par $j$).  Or, nous avons $\mathbb{L}\pi_{!}(h_{j})\simeq h_{\pi(j)}$, et 
le morphisme ci-dessus est donc isomorphe, dans $\mathrm{Ho}(\mathbf{SEns})$, au morphisme naturel
$$\mathrm{Hocolim}_{k\in \pi^{-1}(i)}Hom(k,j) \longrightarrow Hom(i,\pi(j)).$$
On remarque alors que le foncteur $\pi : \Delta(I) \longrightarrow I$
est fibr\'e (c'est \`a dire que $\pi^{op}$ est cofibr\'e au sens de notre
d\'efinition \ref{dcart}). Ceci implique que la fibre du morphisme ci-dessus, prise
en $u \in Hom(i,\pi(j))$ est isomorphe \`a 
$$\mathrm{Hocolim}_{k\in \pi^{-1}(i)}Hom^{u}(k,j)\simeq \mathrm{Hocolim}_{k\in \pi^{-1}(i)}Hom^{id}(k,u^{*}(j)),$$
o\`u $u^{*}(j) \rightarrow j$ est un rel\`evement cart\'esien de $u$, et o\`u 
$Hom^{u}$ (resp. $Hom^{id}$) d\'esigne les sous-ensembles de morphismes
dont l'images par $\pi$ est \'egale \`a $u$ (resp. \`a $id$). Or, on a 
$$\mathrm{Hocolim}_{k\in \pi^{-1}(i)}Hom^{id}(k,u^{*}(j)) \simeq
\mathbb{L}q_{!}(h_{u^{*}(j)}) \simeq h_{\{i\}}\simeq *.$$
Ceci montre bien que la transformation naturelle
$$\mathbb{L}q_{!}v^{*} \Rightarrow u^{*}\mathbb{L}\pi_{!}$$
est un isomorphisme. En particulier, si $F \in \mathrm{Ho}(SPr(\Delta(I))$, le morphisme d'adjonction
$$\pi^{*}\mathbb{L}\pi_{!}(F) \longrightarrow F,$$
\'evalu\'e en $j\in \Delta(I)$, s'\'ecrit
$$\mathbb{L}\pi_{!}(F)(\pi(j))\simeq \mathbb{L}q_{!}v^{*}(F)
 \longrightarrow F(j),$$
o\`u $q : \pi^{-1}(\pi(j)) \longrightarrow *$, et $v^{*}$ est le foncteur
de retriction \`a $\pi^{-1}(\pi(j)) \subset \Delta(I)$. Ce morphisme s'\'ecrit
aussi
$$\mathbb{L}q_{!}v^{*}(F) \simeq \mathrm{Hocolim}_{k \in \pi^{-1}(\pi(j))}F(k) \longrightarrow F(j).$$
Supposons maintenant que $F$ envoie les morphismes verticaux sur des 
\'equivalences. Comme la cat\'egorie $\pi^{-1}(\pi(j))$ poss\`ede un objet initial
(qui est le foncteur constant $[0] \rightarrow I$ \'egal \`a $\pi(j)$, not\'e
simplement $\pi(j)$), le pr\'efaisceau
$v^{*}(F)$ restrint \`a $\pi^{-1}(\pi(j))$ est \'equivalent au pr\'efaisceau constant
de valeurs $F(\pi(j))=X$. Ainsi, il nous reste \`a montrer que le morphisme naturel
$$\mathrm{Hocolim}_{k \in \pi^{-1}(\pi(j))}X \longrightarrow X\simeq F(j)$$
est un isomorphisme. Ceci est finalement \'equivalent au fait que 
$$\mathrm{Hocolim}_{k \in \pi^{-1}(\pi(j))}* \simeq N(\pi^{-1}(\pi(j))) \simeq *.$$
Or, ceci est vrai car la cat\'egorie $ \pi^{-1}(\pi(j))$ poss\`ede un objet final et poss\`ede
donc un nerf contractile. \hfill $\Box$ \\

Revenons maintenant \`a la cat\'egorie $\Delta(I)$. On dispose d'un foncteur
$$\Delta(I)^{op} \longrightarrow \mathbf{Cat}/I\subset \s -\mathbf{Cat}/I,$$
qui \`a un objet $([n],u)$ associe la cat\'egorie $\Delta^{n}$, muni de son foncteur
$u : \Delta^{n} \longrightarrow I$. Ce foncteur induit alors un morphisme bien d\'efini
dans $\mathrm{Ho}(\s -\mathbf{Cat})$
$$\mathrm{Hocolim}_{([n],u)\in \Delta(I)^{op}}\Delta^{n} \longrightarrow I.$$

\begin{prop}\label{pa2}
Le morphisme ci-dessus est un isomorphisme dans $\mathrm{Ho}(\s -\mathbf{Cat})$. 
\end{prop}

\textbf{Preuve --} On commence par r\'e\'ecrire la colomite homotopique en question
de la fa\c{c}on suivante
$$\mathrm{Hocolim}_{([n],u)\in \Delta(I)^{op}}\Delta^{n} \simeq
\mathrm{Hocolim}_{n\in \Delta^{op}} \left( \coprod_{[m_{1}] \rightarrow 
[m_{2}] \rightarrow \dots \rightarrow [m_{n}]} Hom(\Delta^{m_{n}},I)\times \Delta^{m_{1}} \right).$$
Le membre de droite n'est autre que le coend homotopique du diagramme
$$(p,q)\in \Delta^{op} \times \Delta \mapsto \left( Hom(\Delta^{p},I)\times \Delta^{q} \right).$$
Plus g\'en\'eralement, pour un objet simplicial $X_{*} : \Delta^{op} \longrightarrow \s -\mathbf{Cat}^{pr}$, 
nous noterons 
$$|X_{*}|_{\Delta}:=Colim \left( 
\coprod_{[n]} X_{n}\times \Delta^{n} \leftleftarrows \coprod_{[p] \rightarrow [q]}
X_{q} \times \Delta^{p}\right).$$
Ceci d\'efinit un foncteur
$$|.|_{\Delta} : (\s -\mathbf{Cat}^{pr})^{\Delta^{op}} \longrightarrow
\s -\mathbf{Cat}^{pr}.$$
Ce foncteur poss\`ede un adjoint \`a droite 
$$\s -\mathbf{Cat}^{pr} \longrightarrow (\s -\mathbf{Cat}^{pr})^{\Delta^{op}}$$
qui envoie $A$ sur 
$$\begin{array}{cccc}
\underline{Hom}(\Delta^{*},A) : & \Delta^{op} & \longrightarrow &  \s -\mathbf{Cat}^{pr} \\
 & [n] & \mapsto & \underline{Hom}(\Delta^{n},A).
 \end{array}
$$
Nous munissons $(\s -\mathbf{Cat}^{pr})^{\Delta^{op}}$ de la structure de mod\`eles de Reedy de 
\cite[]{ho}, pour laquelle le foncteur $A \mapsto \underline{Hom}(\Delta^{*},A)$
est de Quillen \`a droite. De plus, comme tout objet de $\s -\mathbf{Cat}^{pr}$ est cofibrant on voit 
que les cofibrations de $\s -\mathbf{Cat}^{pr}$ ne sont autre que les monomorphismes (voir par exemple \cite[]{hir}). 
Ainsi, tout objet $X_{*}$ de $\s -\mathbf{Cat}^{pr}$ est cofibrant, ce qui implique que le morphisme naturel
$$\mathbb{L}|X_{*}|_{\Delta} \longrightarrow |X_{*}|_{\Delta}$$
est un isomorphisme dans $\mathrm{Ho}(\s -\mathbf{Cat}^{pr})$. Or, on a
$$\mathbb{L}|X_{*}|_{\Delta}\simeq 
\mathrm{Hocolim}_{n\in \Delta^{op}} \left( \coprod_{[m_{1}] \rightarrow 
[m_{2}] \rightarrow \dots \rightarrow [m_{n}]} X_{m_{n}})\times \Delta^{m_{1}} \right).$$
Cela implique en particulier que l'on a un isomorphisme naturel
$$\mathrm{Hocolim}_{([n],u)\in \Delta(I)^{op}}\Delta^{n} \simeq
|N(I)|_{\Delta},$$
o\`u $N(I)$ est le nerf de $I$ consid\'er\'e comme objet simplicial
$$N(I) : \Delta^{op} \longrightarrow Ens \subset \s -\mathbf{Cat}^{pr}.$$
Ainsi, il nous reste \`a montrer que le morphisme naturel
$$|N(I)|_{\Delta} \longrightarrow I$$
est un isomorphisme dans $\s -\mathbf{Cat}^{pr}$, ce qui est une cons\'equence directe du fait 
que pour un ensemble simplicial $X_{*}$ le morphisme naturel
$$Colim \left( 
\coprod_{[n]} X_{n}\times \Delta^{n} \leftleftarrows \coprod_{[p] \rightarrow [q]}
X_{q} \times \Delta^{p}\right) \longrightarrow 
X_{*}$$
soit un isomorphisme. 
\hfill $\Box$ \\

\section{Comparaison avec le caract\`ere de Chern usuel}\label{comparazione}

Dans cet appendice nous montrons comment le caract\`ere de Chern introduit dans \S 4.2 induit une 
transformation naturelle $Ch : K_0 \rightarrow H^{ev}_{\mathrm{dR}}$ de foncteurs de la cat\'egorie des 
vari\'et\'es quasi projectives et lisses sur un corps $k$ de caract\'eristique nulle, et \`a 
valeurs dans la cat\'egorie des anneaux commutatifs. Nous montrerons alors que cette transformation
naturelle est \'egal au caract\`ere de Chern usuel. \\

\noindent \textbf{Espace des lacets d'un groupe alg\'ebrique.}
Rappelons que pour tout champ $X$ sur $k$, $\mathcal{L}X:=X^{S^{1}} \in \mathbf{dSt}_{k}$ est 
l'espace des lacets deriv\'ees de $X$.  
Soit $G$ un champ alg\'ebrique en groupes commutatifs sur le corps $k$ et $\mathrm{e}:\mathrm{Spec}\, k \longrightarrow G$ sa section unit\'e.
L'espace $\Omega_{\mathrm{e}}G$ des lacets point\'es en $\mathrm{e}$ est defini par le carr\'e homotopiquement cart\'esien dans $\mathbf{dSt}_{k}$ $$\xymatrix{\Omega_{\mathrm{e}}G \ar[r]^-{j} \ar[d] & \mathcal{L}G \ar[d]^-{\textrm{ev}_{0}} \\ 
\mathrm{Spec}\, k \ar[r]_-{\textrm{e}} & G.}$$

\noindent La composition $$\xymatrix{\mathcal{L}(G) \ar[r]^-{\Delta} & \mathcal{L}(G)\times \mathcal{L}(G) \ar[r]^-{\textrm{ev}_{0} \times \textrm{id}} & G \times \mathcal{L}(G) \ar[r]^-{(-)^{-1}\times \textrm{id}} & G \times \mathcal{L}(G) \ar[r]^-{\sigma} & \mathcal{L}(G) ,}$$ $\sigma$ etant induite par la multiplication \`a gauche de $G$ sur lui-meme, induit un morphisme $$\gamma_{\mathrm{e}}: \mathcal{L}G\longrightarrow \Omega_{\mathrm{e}}G$$

Nous aurons besoin du r\'esultat suivant

\begin{prop} Soit $G$ un champ alg\'ebrique en groupes commutatifs sur le corps $k$. Le morphisme $$(\gamma_{\mathrm{e}}, \mathrm{ev}_{0}): \mathcal{L}G \longrightarrow \Omega_{\mathrm{e}}G \times G$$ est un isomorphisme dans $\mathrm{Ho}(\mathbf{dSt}_{k}/G)$.
\end{prop}

\noindent \textbf{Preuve --} L'inverse est donn\'e par la composition $$\xymatrix{G\times \Omega_{\mathrm{e}}G \ar[r]^{\textrm{id}\times j} & G\times \mathcal{L}G \ar[r]^-{\sigma} & \mathcal{L}G  }.$$ 
\hfill $\Box$ \\

\begin{cor} \label{decomp} Si $\mathbb{G}_{m}$ est le sch\'ema en groupes multiplicatifs sur $k$, on a des isomorphismes canoniques dans $\mathrm{Ho}(\textbf{dSt}_{k})$ $$ \mathcal{L}(\mathbb{G}_{m})\simeq \mathbb{G}_{m}\times \Omega_{1}\mathbb{G}_{m},$$ $$\mathcal{L}(B\mathbb{G}_{m})\simeq B\mathbb{G}_{m}\times \mathbb{G}_{m}.$$ 
\end{cor}

\noindent \textbf{Preuve --} Il suffit de remarquer que l'espace des lacets de $B\mathbb{G}_{m}$ point\'es en $1:=\mathrm{e}$ est canoniquement isomorphe \`a $\mathbb{G}_{m}$.
\hfill $\Box$ \\

Comme $\mathcal{L}(B\mathbb{G}_{m})$ est isomorphe  \`a un espace non deriv\'e, l'action de $S^1$ sur $\mathcal{O}(\mathcal{L}(B\mathbb{G}_{m}))$ est triviale, $\mathcal{O}(\mathcal{L}(B\mathbb{G}_{m}))^{\mathrm{h}S^1}\simeq \mathcal{O}(\mathcal{L}(B\mathbb{G}_{m}))\simeq \mathcal{O}(\mathbb{G}_{m})$.\\

\noindent \textbf{Le caract\`ere $Ch:K_0 \rightarrow H^{ev}_{\mathrm{dR}}$.} 
Soit $\mathbf{QProjLis}_{k}$ la cat\'egorie des vari\'et\'es quasi-projectives et lisses sur un corps $k$ de caract\'eristique zero, $\mathbf{ChQuotLis}_{k}$ la sous-cat\'egorie pleine des 1-champs d'Artin sur $k$ qui sont des quotients des vari\'et\'es quasi-projectives lisses par l'action d'un groupe $GL_n$, pour  $n$ quelconque, et $\mathbf{Comm}$ la cat\'egorie des anneaux commutatifs.
On utilisera la correspondence de Dold-Kan (\cite[\S 8.4]{we}) implicitement, pour identifier les complexes cohomologiques en degr\'es non-positifs avec les modules simpliciaux.\\

En utilisant \cite[Thm. 4.1]{rhamloop}, on peut definir le champ non deriv\'e d'\emph{homologie cyclique periodique}
$$\mathrm{HP}:\mathrm{alg}_{k} \longrightarrow \mathrm{SEns}:\,\, R \longmapsto \mathcal{O}(\mathcal{L}(R))^{hS^1}[u^{-1}]$$
o\`u $\mathcal{L}(R)= L(\mathbb{R}\mathrm{Spec}(R))$ dans les notations de \cite{rhamloop}.
Par construction, le caract\`ere de Chern du paragraphe \S 4.2 donne un morphisme des champs non 
deriv\'es 
$$\underline{Ch}:\mathrm{Vect} \longrightarrow \mathrm{HP}$$ 
qu'on peut \'etendre de mani\`ere canonique \`a la cat\'egorie $\mathbf{QProjLis}_{k}$. Pour chaque $F\in \mathbf{ChQuotLis}_{k}$, on obtient donc un application d'ensembles 
$$[F, \mathrm{Vect}] \longrightarrow [F, \mathrm{HP}]$$ 
o\`u $[\,,\,]$ designe l'ensemble des morphismes dans la cat\'egorie homotopique des champs sur $k$, pour la topologie \'etale (\cite[2.1.1]{hagII}). Mais 
$$[F, \mathrm{Vect}]\simeq \mathrm{Vect}(F)/\mathrm{iso}$$ et $[F, \mathrm{HP}]$ est identifi\'e par le lemme suivant

\begin{lem} On a un isomorphisme canonique $[F, \mathrm{HP}]\simeq H_{\textrm{dR}}^{\mathrm{ev}}(F)$
\end{lem}
\noindent \textbf{Preuve --} Soit $U_{\bullet}\rightarrow F$ un atlas simplicial de $F$ form\'e d'objets dans $\mathbf{QProjLis}_{k}$. D'apr\`es  \cite[Cor. 4.2]{rhamloop} on a
$$\mathbb{R}\underline{Hom}(F,\mathrm{HP}) \simeq \mathrm{holim}_{\Delta}\mathbb{R}\underline{Hom}(U_{\bullet},\mathrm{HP})\simeq \mathrm{holim}_{\Delta}C_{\textrm{dR}}(U_{\bullet})[u,u^{-1}]=:C_{\textrm{dR}}(F)[u,u^{-1}],$$
o\`u $C_{\textrm{dR}}(F)[u,u^{-1}]$ est ici le complexe de de Rham $\mathbb{Z}/2$-p\'eriodique de $F$.
 En prenant le $\pi_{0}$ on obtient la bijection $[F, \mathrm{HP}]\simeq H_{\textrm{dR}}^{\mathrm{ev}}(F)$ 
\hfill $\Box$ \\

On a donc une application $Ch^{\textrm{iso}}:\mathrm{Vect}(F)/\mathrm{iso} \longrightarrow  H_{\textrm{dR}}^{\mathrm{ev}}(F)$ fonctoriel par rapport \`a $F\in \mathbf{ChQuotLis}_{k}$. D'apr\'es 
l'additivit\'e (voir Cor. \ref{tracad'}) et la multiplicativit\'e (\S \ref{mult}) 
de la trace cylique on a une factorisation, fonctorielle en $F$

$$\xymatrix{\mathrm{Vect}(F)/\mathrm{iso} \ar[rr]^-{Ch^{\textrm{iso}}} \ar[rd] & &   H_{\textrm{dR}}^{\mathrm{ev}}(F) \\
& K_{0}^{\oplus}(F) \ar[ru]_-{Ch^{\oplus}} &  }$$ 
o\`u $K_{0}^{\oplus}(F)$ est le groupe de Grothendieck des fibr\'es vectoriels sur $F$ (pour la somme
directe, et non pour les suites exactes non-scind\'ees !) 
et $Ch^{\oplus}$ est une transformation naturelle de foncteurs de $\mathbf{ChQuotLis}_{k}$ dans la cat\'egorie $\mathbf{Comm}$ des anneaux commutatifs. 
Pour tout $F\in \mathbf{ChQuotLis}_{k}$, et tout suite exacte $(E)$ 
 de fibr\'es vectoriels sur $F$,  il existe un 
champ $F'$ et une  \'equivalence $\mathbb{A}^{1}$-locale (\cite[Def. 2.1, p. 106]{mv}) $F'\rightarrow F$ telle 
que la suite $(E)$ se scinde sur $F'$ (on peut prendre pour 
$F'$ un torseur sous un fibr\'e vectoriel sur $F$). Comme $F \mapsto H^{ev}_{DR}(F)$ est 
$\mathbb{A}^{1}$-invariant, on en d\'eduit une factorization

$$\xymatrix{K_{0}^{\oplus}(F) \ar[rr]^-{Ch^{\oplus}} \ar[rd] & &   H_{\textrm{dR}}^{\mathrm{ev}}(F) \\
& K_{0}(F) \ar[ru]_-{Ch} &  }$$ 
o\`u cette fois $K_{0}(F)$ est le v\'eritable groupe de Grothendieck des fibr\'es vectoriels sur $F$. \\

Nous venons de voir que le caract\`ere de Chern de \S 4.2 induit une transformation naturelle $$Ch: K_{0} \rightarrow H^{ev}_{\mathrm{dR}}$$ 
entre foncteurs de $\mathbf{ChQuotLis}^{\textrm{op}}_{k}$ vers la cat\'egorie $\mathbf{Comm}$ des anneaux commutatifs. Nous allons maintenant comparer cette transformation naturelle avec 
le caract\`ere de Chern usuel (voir par exemple \cite{gi}) 
$$Ch^{\mathrm{cl}}: K_{0} \rightarrow H^{ev}_{\mathrm{dR}},$$
et montrer que $Ch=Ch^{cl}$ en tant que transformations naturelles. \\

\noindent \textbf{La comparaison $Ch=Ch^{\mathrm{cl}}$.} Le r\'esultat principal de cet appendice est la comparaison suivante.

\begin{thm}\label{compfibvect} Soit $\mathbf{Comm}$ la cat\'egorie des anneaux commutatifs. Notre 
transformation naturelle $Ch$ et le caract\'ere de Chern usuel $Ch^{\mathrm{cl}}$ sont des 
transformations naturelles \'egales (entre fonceturs 
$\mathbf{QProjLis}^{\mathrm{op}}_{k}\rightarrow \mathbf{Comm}$). 
\end{thm}

\noindent \textbf{Preuve --} La preuve consiste en quatre \'etapes.\\

\noindent \textsf{Pr\'emi\`ere \'etape.} Le pr\'emier pas est de caracteriser les transformations naturelles $K_{0}\rightarrow H^{ev}_{dR}$ comme foncteurs definies sur la cat\'egorie $\mathbf{ChQuotLis}_{k}$. C'est le contenu du lemme suivant.

\begin{lem}\label{nategalk} Si $\mathrm{Nat}(K_0,H_{\mathrm{dR}})$ est l'ensemble des transformations naturelles entre foncteurs $$K_0,\, H_{\mathrm{dR}}: \mathbf{ChQuotLis}_{k}^{\mathrm{op}}\rightarrow \mathbf{Comm},$$ on \`a une injection canonique 

$$\mathrm{Nat}(K_0, H_{\mathrm{dR}})\hookrightarrow Hom_{\mathbf{GrpForm}_{k}}(\hat{\mathbb{G}}_{a,\,k}, \hat{\mathbb{G}}_{m,\,k}) \simeq k,$$
o\`u  $Hom_{\mathbf{GrpForm}_{k}}(\hat{\mathbb{G}}_{a,\,k}, \hat{\mathbb{G}}_{m,\,k})$
est l'ensemble des morphismes de groupes formels entre le groupe additif et le groupe multiplicatif. 
\end{lem}

\noindent \textsl{Preuve du lemme} -- Par le \emph{splitting principle} pour $K_0$ (qui est valable pour les objets de $\mathbf{ChQuotjLis}_{k}$), l'evaluation en $B\mathbb{G}_{m}\in \mathbf{ChQuotLis}_{k}$ $$\mathrm{Nat}(K_0, H_{\mathrm{dR}})\rightarrow Hom_{\mathbf{Comm}}(K_0(B\mathbb{G}_{m}), H_{\mathrm{dR}}(B\mathbb{G}_{m}))$$ est injective. Mais, comme $B\mathbb{G}_{m}$ est un champ en groupes abeliens, chaque transformation naturelle $K_0 \rightarrow H_{dR}$ donne en fait un morphisme $$K_0 (B\mathbb{G}_{m})\otimes k \simeq k[t,t^{-1}]\rightarrow k[[v]]\simeq H_{\mathrm{dR}}(B\mathbb{G}_m)$$
d'alg\`ebres de Hopf sur $k$ (\'eventuellement compl\`etes). La propri\'et\'e universelle de la compl\'etion montre que ce morphisme est \'equivalent \`a la donn\'e d'un morphisme de groupes formels sur $k$, $\hat{\mathbb{G}}_{a,\,k}\rightarrow \hat{\mathbb{G}}_{m,\,k}$. Donc, l'\'evaluation en $B\mathbb{G}_{m}$ donne en fait une injection $$\mathrm{Nat}(K_0, H_{dR})\rightarrow Hom_{\mathbf{GrpForm}_{k}}(\hat{\mathbb{G}}_{a,\,k}, \hat{\mathbb{G}}_{m,\,k}).$$  
On conclut en remarquant le fait bien connu que $$k\simeq Hom_{\mathbf{GrpForm}_{k}}(\hat{\mathbb{G}}_{a,\,k}, \hat{\mathbb{G}}_{m,\,k})$$ par l'application $\lambda \longmapsto \exp(\lambda \cdot - ).$
\hfill $\spadesuit$ \\

\noindent \textsf{Deuxi\`eme \'etape.}  Par le lemme pr\'ec\'edent $$\mathrm{Nat}(K_0, H_{dR})\hookrightarrow Hom_{\mathbf{Comm}}(K_0(B\mathbb{G}_{m}), H_{dR}(B\mathbb{G}_{m}))\hookrightarrow Hom_{\mathbf{GrpForm}_{k}}(\hat{\mathbb{G}}_{a,\,k}, \hat{\mathbb{G}}_{m,\,k}) \simeq k,$$ donc, comme $K_{0}(B\mathbb{G}_{m})\simeq k[t,t^{-1}]$ ($t$ etant la classe du fibr\'e en droites universel sur $B\mathbb{G}_{m}$), il suffit de montrer que  $$Ch(B\mathbb{G}_{m})(t)=Ch^{\mathrm{cl}}(B\mathbb{G}_{m})(t)\equiv 1 + v +\textrm{O}(v^2 )$$ dans $H^{ev}_{\mathrm{dR}}(B\mathbb{G}_{m})\simeq k[[v]]$. 
Par definition de $Ch$, on a un diagramme commutatif
$$\xymatrix{\pi_{0}(\mathrm{Vect}(B\mathbb{G}_{m})) \ar[rr]^-{a} \ar[rd]_-{Ch} & &   \mathcal{O}(\mathcal{L}B\mathbb{G}_{m})^{hS^{1}} \ar[ld]^-{\psi} \\
& H_{\mathrm{dR}}(B\mathbb{G}_{m}).  &  }$$  

Le corollaire Cor. \ref{decomp} implique que $$\mathcal{L}B\mathbb{G}_{m} \simeq \mathbb{G}_{m}\times B\mathbb{G}_{m}$$ 
et donc $$\mathcal{O}(\mathcal{L}B\mathbb{G}_{m})^{hS^1}\simeq \pi_{0}(\mathcal{O}(\mathcal{L}B\mathbb{G}_{m})^{hS^1})\simeq \mathcal{O}(\mathbb{G}_{m})\simeq k[t, t^{-1}].$$
De plus, le fibr\'e en droites universel sur $B\mathbb{G}_{m}$ envoy\'e par le morphisme $a$ sur $t$. On 
se ram\`ene ainsi \`a calculer $\psi(t)\in H_{\mathrm{dR}}(B\mathbb{G}_{m})$.\\

\noindent \textsf{Troisi\`eme \'etape.}  Pour calculer $\psi(t)$ on utilise la r\'esolution 
simpliciale standard $$\xymatrix{\mathrm{B}^{\bullet}\mathbb{G}_{m}  \equiv  (\cdots  \ar@<1ex>[r] \ar@<-1ex>[r]^-{\cdots} & \mathbb{G}_{m} \ar@<1ex>[r] \ar@<-1ex>[r] & \mathrm{Spec}(\mathbb{C}) \ar[l]) \ar[r]_-{s} & B\mathbb{G}_{m}}.$$ On obtient le diagramme suivant

$$\xymatrix{\mathcal{O}(\mathcal{L}B\mathbb{G}_{m})^{hS^1} \ar[r]^-{\sim} \ar[d]_-{h}^-{\sim} & k[t,t^{-1}] & \\
\mathrm{Holim}_{n\in \Delta} \mathcal{O}(\mathcal{L}(\mathbb{G}_{m}^{n}))^{hS^1} \ar[r]^-{\tau_{\leq 1}} \ar[d]_-{\rho} &  \mathrm{Holim}\, \tau_{n\leq 1}(\mathcal{O}(\mathcal{L}(\mathbb{G}_{m}^{n}))^{hS^1})  \ar[r]^-{\textrm{oub-}S^{1}} \ar[d]_-{\rho_{\leq 1}} &  \mathrm{Holim}\, \tau_{n\leq 1}(\mathcal{O}(\mathcal{L}(\mathbb{G}_{m}^{n}))) \ar[d]^-{\rho_{\emptyset}} \\
\mathrm{Holim}_{n\in \Delta} \mathrm{C}_{\textrm{dR}}^{\mathbb{Z}/2}(\mathbb{G}_{m}^{n}) \ar[r]^-{\tau_{\leq 1}} \ar[d]_-{p} & \mathrm{Holim} \, \tau_{n\leq 1} (\mathrm{C}_{\textrm{dR}}^{\mathbb{Z}/2}(\mathbb{G}_{m}^{n})) \ar[r]^-{\textrm{oub-}d} \ar[d]^-{p_{\leq 1}} & \mathrm{Holim} \, \tau^{\Feyn d}_{n\leq 1} (\mathrm{C}_{\textrm{dR}}^{\mathbb{Z}/2}(\mathbb{G}_{m}^{n})) \ar[d]^-{p_{\emptyset}}\\
k[[v]]\simeq H_{\textrm{dR}}^{ev}(B\mathbb{G}_{m}) \ar[r]_{q} & k\oplus H_{\textrm{dR}}^{1}(\mathbb{G}_{m}) \ar[r]_{q'} & k[t,t^{-1}] \oplus \Omega_{\textrm{dR}}^{1}(\mathbb{G}_{m})}$$ 

o\`u \begin{itemize}
\item $\tau_{n\leq 1}(\mathcal{O}(\mathcal{L}(\mathbb{G}_{m}^{n}))^{hS^1}):=(\xymatrix{k \ar@<1ex>[r] \ar@<-1ex>[r] & \mathcal{O}(\mathcal{L}\mathbb{G}_{m})^{hS^1} \ar[l]})$ est la restriction du diagramme cosimpliciale a la sous-cat\'egorie pleine $\Delta_{\leq 1}$ de $\Delta$ des objets $[0]$ et $[1]$;
\item $\tau_{n\leq 1}(\mathcal{O}(\mathcal{L}(\mathbb{G}_{m}^{n})))$ est le diagramme (d'alg\`ebres simpliciales) sur $\Delta_{\leq 1}$, $(\xymatrix{k \ar@<1ex>[r] \ar@<-1ex>[r] & \mathcal{O}(\mathcal{L}\mathbb{G}_{m})  \ar[l]})$ ;
\item $\tau_{n\leq 1} (\mathrm{C}_{\textrm{dR}}^{\mathbb{Z}/2}(\mathbb{G}_{m}^{n})):=(\xymatrix{k \ar@<1ex>[r] \ar@<-1ex>[r] & \mathrm{C}_{\textrm{dR}}^{\mathbb{Z}/2}(\mathbb{G}_{m}) \ar[l]})$ est la restriction du diagramme cosimpliciale (de complexes de $k$-modules) a la sous-cat\'egorie pleine $\Delta_{\leq 1}$ de $\Delta$ des objets $[0]$ et $[1]$;
\item $\tau^{\Feyn d}_{n\leq 1} (\mathrm{C}_{\textrm{dR}}^{\mathbb{Z}/2}(\mathbb{G}_{m}^{n}))$ est le diagramme (de $k$-modules) sur $\Delta_{\leq 1}$, $(\xymatrix{k \ar@<1ex>[r] \ar@<-1ex>[r] & (k[t,t^{-1}]\oplus \Omega_{\textrm{dR}}^{1}(\mathbb{G}_{m}) \ar[l]})$. Notons que en fait $\mathrm{Holim} \, \tau^{\Feyn d}_{n\leq 1} (\mathrm{C}_{\textrm{dR}}^{\mathbb{Z}/2}(\mathbb{G}_{m}^{n}))= \mathrm{lim} \, \tau^{\Feyn d}_{n\leq 1} (\mathrm{C}_{\textrm{dR}}^{\mathbb{Z}/2}(\mathbb{G}_{m}^{n}))$;
\item si $X$ est un sch\'ema lisse sur $k$, $\mathrm{C}_{\textrm{dR}}^{\mathbb{Z}/2}(X)$ est le complexe $2$-periodique de de Rham alg\'ebrique de $X$, $\xymatrix{\Omega_{\textrm{dR}}^{\textrm{even}}(X)\ar[r]^-{d} & \Omega_{\textrm{dR}}^{\textrm{odd}}(X)}$;
\item les morphismes $\rho$ (deuxi\`eme ligne) r\'esultent de \cite[Th. 1.1]{rhamloop} ;
\item les morphismes $p$ (troisi\`eme ligne) d\'esignent le passage au $\pi_{0}$ (et $p_{\emptyset}$ est un isomorphisme);
\item le morphisme compos\'e $p\circ \rho \circ h$ est $\psi$ (definie \`a la 2me \'etape).
\end{itemize}
\medskip

Le lemme suivant identifie le morphisme compos\'e $q'\circ q$.

\begin{lem}\label{koszul} Le morphisme $q'\circ q : k[[v]]\simeq H_{\mathrm{dR}}^{ev}(B\mathbb{G}_{m}) \longrightarrow k[t,t^{-1}]\oplus \Omega_{\mathrm{dR}}^{1}(\mathbb{G}_{m})$ envoie $1$ sur $(1,0)$ et $v$ sur $(0,dt/t)$.
\end{lem}

\noindent \textsl{Preuve du lemme} -- On a $\mathrm{C}_{\textrm{dR}}^{\mathbb{Z}/2}(\mathbb{G}_{m})\simeq \mathcal{O}(\mathbb{G}_{m})\oplus \mathfrak{g}^{\vee}$ o\`u $\mathfrak{g}\simeq k\cdot (dt/t)$ est l'alg\`ebre de Lie de $\mathbb{G}_{m}$. Le r\'esultat suit alors par dualit\'e de Koszul.
\hfill $\spadesuit$ \\

\noindent Gr\^ace aux diagrammes et lemmes pr\'ec\'edents, il suffit donc de montrer que, dans le diagramme $$\xymatrix{k[t,t^{-1}]\simeq \mathcal{O}(\mathcal{L}B\mathbb{G}_{m})^{hS^1} \ar[r]^-{h_{\leq1}^{\Feyn S^{1}}}   &
 \mathrm{Holim}\, \tau_{n\leq 1}(\mathcal{O}(\mathcal{L}(\mathbb{G}_{m}^{n}))) \ar[r]^-{\rho_{\emptyset}} & \mathrm{Holim} \, \tau^{\Feyn d}_{n\leq 1} (\mathrm{C}_{\textrm{dR}}^{\mathbb{Z}/2}(\mathbb{G}_{m}^{n})) \ar[r]^-{p_{\emptyset}} & k[t,t^{-1}] \oplus \Omega_{\textrm{dR}}^{1}(\mathbb{G}_{m})}$$ (o\`u $h_{\leq1}^{\Feyn S^{1}}:= (\mathrm{oub-}S^1) \circ \tau_{\leq 1}\circ h$), l'image de $t\in k[t,t^{-1}]$ est $(1,dt/t)$.\\

\noindent \textsf{Quatri\`eme \'etape.} La r\'esolution standard $s$ $$\xymatrix{\mathrm{B}^{\bullet}\mathbb{G}_{m}  \equiv  (\cdots  \ar@<1ex>[r] \ar@<-1ex>[r]^-{\cdots} & \mathbb{G}_{m} \ar@<1ex>[r] \ar@<-1ex>[r] & \mathrm{Spec}(\mathbb{C}) \ar[l]) \ar[r]_-{s} & B\mathbb{G}_{m}}$$ (qui induit $h$ et donc $h_{\leq1}^{\Feyn S^{1}}$) donne une auto-homotopie canonique de la composition $$s_1:\xymatrix{\mathbb{G}_{m} \ar@<1ex>[r] \ar@<-1ex>[r] & \mathrm{Spec}(\mathbb{C}) \ar[r]_-{s} & B\mathbb{G}_{m}}.$$ Cette homotopie, par adjonction et Lemme \ref{decomp}, est donn\'e par le morphisme $$H:\mathbb{G}_{m}\longrightarrow \mathcal{L}(B\mathbb{G}_{m})\simeq \mathbb{G}_{m}\times B\mathbb{G}_{m}$$ dont la projection sur $B\mathbb{G}_{m}$ est $s_1$ et la projection sur $\mathbb{G}_{m}$ est l'identit\'e. Cette $H$ donne fonctorielment une auto-homotopie $H^{\mathcal{L}}$ de $L$ de la facon suivante.\\ Pour donner la construction, on se place dans une cadre un peu plus g\'en\'eral. Soit $T$ un objet de $\mathbf{dSt}_{k}$ et notons aussi par $T$ le foncteur $$\mathrm{Ho}(\mathbf{dSt}_{k}) \longrightarrow \mathrm{Ho}(\mathbf{dSt}_{k}) \, :\, X\longmapsto \mathbb{R}\mathrm{HOM}(T,X). $$ Si on a une morphisme $f:X\rightarrow Y$ dans $\mathrm{Ho}(\mathbf{dSt}_{k})$ et une auto-homotopie $H:X\rightarrow Y^{S^1}$ de $f$, on obtient aussi une auto-homotopie $H^{T}: T(X)\rightarrow (T(Y))^{S^1}$ de $T(f)$ en considerant la composition $$ \xymatrix{T(X) \ar[r]^-{T(H)} & T(Y^{S^1}) \simeq \mathbb{R}\mathrm{HOM}(T,Y^{S^1}) \ar[r]^-{\sim} & \mathbb{R}\mathrm{HOM}(T\times S^1,Y) \ar[r]^-{\mathrm{tw}} & \mathbb{R}\mathrm{HOM}(S^1\times T,Y) \ar[r]^-{\sim} & (T(Y))^{S^1} }$$ o\`u on a utilis\'e l'isomorphisme de twist $T\times S^1 \simeq S^1 \times T$.\\
L'application de cette construction au cas $X=\mathbb{G}_{m}$, $Y=B\mathbb{G}_{m}$ et $T=S^1$, nous donne l'auto-homotopie (de $\mathcal{L}(s_1)$) suivante $$H^{\mathcal{L}}: S^1 \times \mathcal{L}(\mathbb{G}_{m}) \longrightarrow \mathcal{L}(B\mathbb{G}_{m})\simeq \mathbb{G}_{m}\times B\mathbb{G}_{m}$$ qui 
\begin{itemize}
\item sur la composante $\mathbb{G}_{m}$ de $\mathcal{L}(B\mathbb{G}_{m})$ est donn\'e par $$\xymatrix{S^1 \times \mathcal{L}(\mathbb{G}_{m}) \ar[r]^-{\textrm{id}\times \varphi} & S^1 \times \mathcal{L}(\mathbb{G}_{m}) \ar[r]^-{\textrm{ev}} & \mathbb{G}_{m}  }$$  o\`u $\varphi$ est la composition $$\xymatrix{\mathcal{L}(\mathbb{G}_{m}) \ar[r]^-{\Delta} & \mathcal{L}(\mathbb{G}_{m})\times \mathcal{L}(\mathbb{G}_{m}) \ar[r]^-{\textrm{ev}_{0} \times \textrm{id}} & \mathbb{G}_{m} \times \mathcal{L}(\mathbb{G}_{m}) \ar[r]^-{\textrm{inv}\times \textrm{id}} & \mathbb{G}_{m} \times \mathcal{L}(\mathbb{G}_{m}) \ar[r]^-{\sigma} & \mathcal{L}(\mathbb{G}_{m}) }$$ $\sigma$ etant induite par la multiplication \`a gauche de $\mathbb{G}_{m}$ sur lui-meme;
\item sur la composante $B\mathbb{G}_{m}$ de $\mathcal{L}(B\mathbb{G}_{m})$ par la composition $$\xymatrix{S^1 \times \mathcal{L}(\mathbb{G}_{m})\ar[r]^-{\textrm{id}\times \textrm{pr}_{k}} & S^1 \times \mathrm{Spec}\, k \ar[r]^-{\textrm{id}\times 1} & S^1 \times \mathbb{G}_{m} \ar[r]^-{H^{\flat}} & B\mathbb{G}_{m} }$$  o\`u $H^{\flat}$ est adjointe \`a $H$ et $1$ la section unit\'e de $\mathbb{G}_{m}$.
\end{itemize}
L'auto-homotopie $H^{\mathcal{L}}$ donne lieu, en passant aux fonctions $\mathcal{O}$, a deux morphismes $$a_{0}:= \pi_{0}(H^{\mathcal{L}}(0,-)^{*}): k[t, t^{-1}]\simeq \pi_{0}(\mathcal{O}(\mathcal{L}(B\mathbb{G}_{m}))) \longrightarrow \pi_{0}(\mathcal{O}(\mathcal{L}(\mathbb{G}_{m})))\simeq \mathcal{O}(\mathbb{G}_{m})$$ qui envoie $t$ sur $1$,  et $$a_{1}: k[t, t^{-1}]\simeq \pi_{0}(\mathcal{O}(\mathcal{L}(B\mathbb{G}_{m}))) \longrightarrow \pi_{1}(\mathcal{O}(\mathcal{L}(\mathbb{G}_{m})))\simeq \Omega_{\textrm{dR}}^{1}(\mathbb{G}_{m})$$ qui envoie $t$ sur $dt/t$ (par definition de $\varphi$, qui donne le facteur $1/t$, et par \cite[Th. 4.1]{rhamloop}, qui montre comme l'homotopie universelle $\mathrm{ev}: S^1 \times \mathcal{L}X \longrightarrow X$ induit sur les $\mathcal{O}$ le diff\'erentiel de de Rham, pour $X$ lisse sur $k$). Puisque la composition $p_{\emptyset} \circ \rho_{\emptyset} \circ h_{\leq1}^{\Feyn S^{1}}$ est egale \`a la somme $a_0 + a_1$, on conclut que telle composition envoie $t$ sur $(1,dt/t)$.

\hfill $\Box$ \\

\begin{rmk}\label{compparf} \emph{On peut montrer que le r\'esult\`at du Th\'eor\`eme \ref{compfibvect} est aussi valable pour le complexes parfaits, et en est en fait un corollaire. Nous esquissons ici les grandes
 lignes de la d\'emonstration. D'abord, notons que le foncteur alg\`ebre de de Rham $$\mathrm{HP}: \mathbf{QProjLis}_{k}\rightarrow \mathbf{Ch}(k)$$ envoie les $\mathbb{A}^1$-equivalences (au sens de \cite{mv}) sur des quasi-isomorphismes. Mais, en $\mathbb{A}^1$-homotopie des sch\'emas on dispose de $\mathbb{A}^1$-equivalences $$\xymatrix{BGL_{\infty} \ar[r]^{\sim} & \mathrm{Vect}^{+}\ar[r]^{\sim} & \mathrm{Parf}}.$$
 Ainsi, pour le complexe parfait universel $\mathcal{E}$ sur $\mathrm{Parf}$ on aura $Ch (\mathcal{E})=Ch^{cl} (\mathcal{E})$ dans $H_{dR}^{ev}(\mathrm{Parf}):=\pi_{0}(\mathcal{O}(L\mathrm{Parf})^{hS^1})$ si et seulement si $Ch (\mathcal{E})=Ch^{cl} (\mathcal{E})$ dans $H_{dR}^{ev}(B\mathbb{G}L_{n})$ pour tout $n\geq 0$. Donc il suffit de montrer que $Ch_{B\mathbb{G}L_{n}}=Ch^{cl}_{B\mathbb{G}L_{n}}$ pour tout $n\geq 0$. Mais ce decoule du Th\'eor\`eme \ref{compfibvect} gr\^ace \`a l'observation suivante du \`a B. Totaro: si $F=[Y/\mathbb{G}L_{n}]$ est un 1-champ quotient d'une vari\'et\'e projective lisse $Y$, alors pour chaque $n_{0}$ il existe une vari\'et\'e projective lisse $X_{n_0}$ et un morphisme $p_{n_0}:X_{n_0}\rightarrow F$  tel que l'image inverse induit des isomorphismes $H^{m}_{dR}(F)\simeq H_{dR}^{m}(X_{n_0})$ pour tout $m\leq n_{0}$.}
\end{rmk}

\end{appendix}

\end{document}